\documentstyle[10pt]{amsart}
\newcommand{\R}{{\Bbb R}}
\newcommand{\Z}{{\Bbb Z}}
\newcommand{\Nb}{{\Bbb N}}

\newcommand{\C}{{\Bbb C}}

\newcommand{\half}{\frac{1}{2}}
\newcommand{\Pp}{{\cal P}}
\newcommand{\I}{{\cal I}}

\theoremstyle{plain}

\title{Wave invariants at elliptic closed geodesics}
\author{Steven Zelditch*}
\date{July 1995}
\thanks{ *Partially supported by NSF grant \#DMS-9404637.}
\address{Johns Hopkins University, Baltimore, Maryland  21218 \newline
email:{\tt zel\@chow.mat.jhu.edu}  or  {\tt zelditch\@msri.org}}

\begin{document}
\maketitle

\addtolength{\baselineskip}{1pt} 

\setcounter{section}{-1}
\section{Introduction}

The purpose of this article is to provide an effective method for calculating the wave invariants 
 associated to a non-degenerate elliptic closed geodesic $\gamma$ of a compact 
Riemannian manifold $(M,g)$: that is, the coefficients in the singularity expansion
$$Tr U(t) = e_o(t) + \sum_{\gamma}  e_{\gamma}(t)$$
$$e_{\gamma}(t) \sim a_{\gamma\;-1}(t - L_{\gamma} + i0)^{-1} + \sum_{k=0}^{\infty}
a_{\gamma\;k}(t-L_{\gamma} +i0)^k log (t - L_{\gamma} + i0)$$
of the trace of the wave group $U(t):= expit\sqrt{\Delta}$ at $t=L_{\gamma}$ (the length of $\gamma.$)
We will show that
$$a_{\gamma\;k} = \int_{\gamma} I_{\gamma\\;k}(s) ds$$
for certain  homogeneous invariant  densities $I_{\gamma\;k}(s)ds$ on $\gamma$, given by
at most 2k+1 integrals over $\gamma$ of polynomials
 in the curvature, Jacobi fields, length, inverse length and Floquet invariants
$\beta_j:=(1 - e^{i\alpha_j})^{-1}$ along $\gamma.$   
  These expressions  characterize the wave invariants
in much the same way that the heat invariants (or wave invariants at t=0) are characterized
as integrals $\int_M P_j(R, \nabla R, ...)dvol$ of homogeneous curvature 
polynomials over M [ABP] [Gi].  Moreover,  in combination with the recent inverse results
of Guillemin [G.1,2], the method produces a list of new spectral invariants
of this kind, simpler than the wave invariants themselves  (the so-called 
quantum Birkhoff normal form coefficients $B_{\gamma k;j}$.)

To state the results, we will need some notation.  We let ${\cal J}_{\gamma}^{\bot}
\otimes \C$ denote the space of complex normal Jacobi fields along $\gamma$, a symplectic
vector space of (complex) dimension 2n (n=dim M-1) with respect to the Wronskian
$$\omega(X,Y) = g(X, \frac{D}{ds}Y) - g(\frac{D}{ds}X, Y).$$
 The linear Poincare map $P_{\gamma}$ is then the linear symplectic map on ${\cal J}_{\gamma}
^{\bot} \otimes \C$ defined by $P_{\gamma} Y(t) = Y(t + L_{\gamma}).$  We will assume
$\gamma$ to be non-degenerate elliptic, i.e. that the eigenvalues of $P_{\gamma}$ are
of the form $\{ e^{\pm i \alpha_j}, j=1,...,n\}$ with (Floquet) exponents $\{\alpha_1,
...,\alpha_n\}$, together with $\pi$,  independent over ${\bf Q}$. The associated normalized eigenvectors
will be denoted $\{ Y_j, \overline{Y_j}, j=1,...,n \}$,
$$P{\gamma} Y_j = e^{i \alpha_j}Y_j \;\;\;\;\;\;P_{\gamma}\overline{Y}_j =
e^{-i\alpha_j} \overline{Y}_j \;\;\;\; \omega(Y_j, \overline{Y}_k) = \delta_{jk}$$
 and relative to a fixed parallel
normal frame $e(s):= (e_1(s),...,e_n(s))$ along $\gamma$ they will be written in the form
$Y_j(s)= \sum_{k=1}^n y_{jk}(s)e_k(s).$  The metric
coefficients $g_{ij}$ will always be taken relative to Fermi
normal coordinates $(s,y)$ along $\gamma$. The mth jet of $g$ along $\gamma$ will
be denoted by $j^m_{\gamma}g$, the curvature tensor  by
$R$ and its covariant derivatives by $\nabla^m R$.  The vector fields 
$\frac{\partial}{\partial s},\frac{\partial}{\partial y_j}$ and their real linear combinations
will be referred to as Fermi normal vector fields along $\gamma$ and contractions of
tensor products of the $\nabla^m R$'s with these vector fields will be referred to as
{\it Fermi curvature polynomials}.   Such polynomials will be called {\it invariant} if
they are invariant under the action of $O(n)$ in the normal spaces.   Invariant contractions
against $\frac{\partial}{\partial s}$ and against the Jacobi eigenfields $Y_j, \overline{Y}_j$,
with coefficients given by invariant polynomials in the components $y_{jk}$,
will be called {\it Fermi-Jacobi polynomials}.  We will alse use this term for functions
on $\gamma$ given by repeated indefinite integrals over $\gamma$ of such FJ polynomials.
Finally, FJ polynomials whose coefficients are given by polynomials in the Floquet invariants
$\beta_j = (1- e^{i\alpha_j})^{-1}$ will be called {\it Fermi-Jacobi-Floquet} polynomials.
 We give weights to
the variables $g_{ij}, D_{s,y}^{\beta} g_{ij}, L:= L_{\gamma}, \alpha_j, y_{ij}, \dot{y}_{ij}$
as follows: $wgt( D_{s,y}^{\beta} g_{ij}) = -|\beta|, wgt(L) = 1, wgt(\alpha_j) = 0,
wgt(y_{ij})= \half, wgt(\dot{y}_{ij}) = -\half.$  As will be seen, these weights reflect the
scaling of these objects under $g \rightarrow \epsilon^2 g$.  A polynomial in this data is
homogeneous of weight s if all its monomials  have weight s under this scaling.
\medskip

\noindent{\bf Theorem A} {\it  Let ${\gamma}$ be an elliptic closed geodesic
with $\{\alpha_1,...,\alpha_n, \pi \}$ independent over ${\bf Q}$.  Then 
$a_{\gamma k} = \int_{\gamma} I_{\gamma k} (s; g)ds$
where :

\noindent(i)  $I_{\gamma k }(s,g)$  is a homogeneous Fermi-Jacobi-Floquet
 polynomial  of weight -k - 1 in the data
$\{L,  y_{ij}, \dot{y}_{ij}, D^{\beta}_{s,y}g \}$ with $|\beta| \leq 2k+4$ ;

\noindent(ii) The degree of $I_{\gamma k }$ in the Jacobi field components is at most
6k+6;

\noindent(iii) At most 2k+1 indefinite integrations over $\gamma$ occur in $I_{\gamma k }$;

\noindent(iv) The degree of $I_{\gamma k }$ in the Floquet invariants $\beta_j$ is
at most k+2.}

\medskip

For instance, in dimension 2  the residual wave invariant $a_{\gamma o}$ is given by:
$$a_{\gamma o} = \frac{c_{\gamma}}{L^{\#}}[ B_{\gamma o;4} (2\beta^2 - \beta
-\frac{3}{4}) + B_{\gamma o;0} ]$$
where:

\noindent (a) $c_{\gamma}$ is the principal wave invariant
$i^{\sigma}L^{\#} |I-P_{\gamma}|^{-\half}$;

\noindent(b) $L^{\#}$ the primitive length of $\gamma$; $\sigma$ is its Morse index; $P_{\gamma}$
is its Poincare map;

\noindent(c) $B_{\gamma o; j}$ have the form:
$$B_{\gamma o;j} = \frac{1}{L^{\#}} \int_o^{L^{\#}} [a\; |\dot{Y}|^4 +
 b_1\; \tau |\dot{Y}\cdot Y|^2 + b_2\; \tau Re (\bar{Y} \dot{Y})^2 
+c \;\tau^2 |Y|^4 + d \;\tau_{\nu \nu}|Y|^4 + e\; \delta_{jo} \tau ]ds  $$
$$ +\frac{1}{L^{\#}}\sum_{0\leq m,n \leq 3; m+n=3} C_{1;mn}\;\frac{sin((n-m)\alpha)}{|(1-e^{i(m-n)\alpha})|^2}\;
|\int_o^{L^{\#}}  \tau_{\nu}(s)\bar{Y}^m\cdot Y^n(s)ds|^2 $$
$$+\frac{1}{L^{\#}}\sum_{0\leq m,n \leq 3; m+n=3} C_{2;mn}\;Im \;\{\int_o^{L^{\#}} \tau_{\nu}(s)\bar
{Y}^m\cdot Y^n(s)[\int_o^s\tau_{\nu}(t)\bar {Y}^n\cdot Y^m(t)dt] ds \}$$
for various universal (computable) coefficients. Here,

\noindent (d)  $\tau$ denotes the scalar curvature, $\tau_{\nu}$  its unit normal derivative, $\tau_{\nu \nu}$
the Hessian $Hess(\tau)(\nu,\nu)$;
 $Y$ denotes the unique
normalized Jacobi eigenfield, $\dot{Y}$ its time-derivative and $\delta_{jo}$ the Kronecker symbol
(1 if j=0 and otherwise 0.) 

We note that the residual wave invariant already saturates the description in Theorem A.
 
This characterization of the wave invariants 
makes more concrete (for Laplacians) the recent results of Guillemin [G1,2] which show that the
wave invariants may be expressed  in terms of the quantum Birkhoff normal form coefficients
of the wave operator around $\gamma.$  For instance, Guillemin shows [G.2, (8.24)]:
$$a_{\gamma o} = ic_{\gamma}\sigma[\sum_{i\not=j}\frac{\partial^2}{\partial I_i \partial I_j}
H_1(0,\sigma)(\beta_i + \half)(\beta_j + \half)$$
$$+ \half \sum_i (2\beta_i^2 + \beta_i - \frac{1}{4})\frac{\partial^2}{\partial I_i^2}H_1(0,
\sigma) -\sum_i(\beta_i + \half)\frac{\partial}{\partial I_i}H_o(0,\sigma) + 
H_{-1}(0,\sigma) ]$$
where $H_{1-r}(\sigma, I_1,\dots,I_n)$ is the term of order 1-r in the complete symbol
of the quantum normal form.  The coefficients $B_{\gamma o;j}$ above are essentially these QBNF (quantum
Birkhoff normal form) coefficients. 

Indeed, the first step in the proof of Theorem A (which we call Theorem B)  
 consists in explicitly constructing the QBNF for
$\sqrt{\Delta}$ near $\gamma$. 
 The method is different from that in [G.1,2] and effectively calculates
 the QBNF coefficients  as integrals of Fermi-Jacobi
 polynomials.  It 
 is  based to some extent on the construction of 
a complete set of quasi-modes associated to $\gamma$ as presented in Babich-Buldyrev
[B.B] although the emphasis is on   intertwining operators rather than
on quasi-modes per se.  In the construction of the intertwining operators it also employs several
ideas in Sjostrand [Sj], although it does not begin by putting the principal symbol in
Birkhoff normal form as is done in [Sj] and [G.2]. 

The second step (which we call Theorem C)  consists in calculating the wave invariants
from the normal form.  This is possible because
 the wave invariants are non-commutative residues of
the wave operator and its $t$ derivatives and hence are invariant under conjugation by unitary
Fourier Integral operators (cf. [G.1,2],
[Z.1]).  The  residues of the normal form wave group will be easily seen  
 to be polynomials in the QBNF
coefficients and in the $\beta_j$'s. 
 Some dimensional analysis of the wave coefficients then leads to the description
in the statement of Theorem A.

 Guillemin [loc.cit] has also proved the remarkable inverse result that, conversely, the QBNF
coefficients can be determined from the wave invariants associated to $\gamma$ and
its iterates, and therefore are themselves spectral invariants.  In view of Theorems A-B this
gives a  list of new spectral invariants, 
 for instance the QBNF coefficients $B_{\gamma o;j}$ above,
in the form of geodesic integrals of FJ polynomials,

Let us now describe the ingredients of the proofs more precisely. Henceforth we will reserve the
notation $\gamma$ for a primitive closed geodesic and will denote its iterates by $\gamma^m$.
 To introduce the
quantum Birkhoff normal form at $\gamma$ and its role in the calculations,   we note first that the wave
invariants associated to $\gamma$ are determined by
the microlocalization of $\Delta$ to a conic neighborhood
$$|y| < \epsilon \;\;\;\;\;\;\;\;\;\;\;\;\;\; \frac{|\eta|}{\sigma} < \epsilon \leqno(0.1)$$
of the cone $\R^{+} \gamma$ thru $\gamma$ in $T^*M - 0$.  Here, $(\sigma, \eta)$ 
  denote the symplectically dual coordinates to the Fermi normal
coordinates $(s,y)$ above.  The microlocalization
of $\Delta$ to (0.1) is then given by $$\Delta_{\psi}:=\psi(s,D_s, y, D_y)^* \Delta 
\psi(s,D_s,y,D_y)$$ where  $D_{x_j}:=\frac{\partial}{i\partial x_j}$ and where
$\psi(s, \sigma, y, \eta)$ is supported in (0.1) and identically one in some smaller
conic neighborhoood .  Often we omit explicit mention of the microlocal cut-off
$\psi$ in calculations which are valid on its microsupport.
  Under the exponential map
$$exp : N_{\gamma} \rightarrow M$$
along the normal bundle to $\gamma$,  the localization of $g$, resp. $\Delta$, to the
tubular neighborhood $|y| < \epsilon$ pulls back to a locally well-defined metric, respectively
Laplacian, on a similar neighborhood in $N_{\gamma}$.  Hence $exp$ conjugates $\Delta_{\psi}$
to an isometric microlocalized Laplacian on $N_{\gamma}$, which we continue to denote
by $\Delta_{\psi}$.  We are thus reduced to calculating the wave invariants of a Laplacian
on the model space $S^1_L \times \R^n$ at the closed geodesic $\gamma = S^1_L \times \{0\}$,
where $S^1_L := \R /L\Z.$  For the sake of simplicity we  assume the normal bundle is orientable
but note that the reduction is valid for immersed as well as embedded closed geodesics.  

We now wish  to put $\Delta_{\psi}$ into normal form,
 which is first of all to conjugate it (modulo a small error)
  into a distinguished
maximal abelian algebra ${\cal A}$ of pseudodifferential operators on the model space
 $S^1_L \times\R^n$.  Roughly speaking, ${\cal A}$ is generated by the tangential operator
$D_s:=\frac{\partial}{i \partial s}$ on $S^1_L$ together with the transverse harmonic
oscillators
$$I_j=I_j(y,D_y) := \frac{1}{2} (D_{y_j}^2 + y_j^2) .\leqno(0.2)$$
 In the construction of the normal form, a special role
will be played by the distinguished element  
$${\cal R}:=\frac{1}{L}( L D_s + H_{\alpha}) \leqno(0.3)$$
where
$$H_{\alpha}:= \frac{1}{2}\sum_{k=1}^n \alpha_k I_k \leqno(0.4)$$
and where the choice of sign in $\pm \alpha_k$ will be specified below.  This element comes
up naturally as the semi-classical parameter in the construction of quasi-modes, although
$D_s$ is more suitable for analysing the wave invariants. 
  Note that both are elliptic elements in the
conic neighborhood
$$ |I_j| < \epsilon \sigma \;\;\;\;\;\;\;\;\;\; I_j(y,\eta):= \frac{1}{2} 
(y_j^2 + \eta_j^2),\leqno(0.5)$$
which will be the image of (0.1) under the conjugation to normal form. 

The classical Birkhoff normal form theorem states roughly the following: near a non-
degenerate elliptic closed geodesic $\gamma$ the Hamiltonian
 $$H(x,\xi)=|\xi|:= \sqrt{\sum_{ij=1}^{n+1} g^{ij}\xi_i\xi_j}$$
 can
be conjugated by a homogeneous local canonical transformation $\chi$  to the normal form
$$\chi^*H \equiv 
  \sigma + \frac{1}{L} \sum_{i,j =1}^n \alpha_j I_j +\frac{p_1(I_1,...,I_n)}{\sigma}+...\;\;\;\; 
mod\;\;\;\;\; O^1_{\infty}\leqno(0.6)$$
where $p_k$ is homogeneous of order k+1 in $I_1,\dots, I_n$, and where
 $O^1_{\infty}$ is the space of germs
of functions homogeneous of degree 1 which vanish to infinite order along $\gamma.$
Note that all the terms in (0.6) are homogenous of degree 1 in $(\sigma, I_1,...,I_n)$, and
that the order of vanishing at $|I|=0$ equals one plus the order of decay in $\sigma$.  The
coefficients of the monomials in the $p_j(I_1,\dots,I_n)$ are known as the classical
Birkhoff normal form invariants.  (See the Appendix for some further details).
 
The quantum Birkhoff normal form is the more or less analogous statement on the operator level.  In
the following the symbol $\equiv$ means that the two sides agree modulo operators whose
complete symbols are of order 1 and vanish to infinite order on $\gamma.$  Also, $O_m \Psi^r$
denotes the space of pseudodifferential operators of order $r$ whose complete symbols
vanish to order $m$ at $(y,\eta)=(0,0).$
\medskip

\noindent{\bf Theorem B} {\it There exists a microlocally elliptic
  Fourier Integral operator $W$ 
from the conic neighborhood (0.1) of $\R^+\gamma$ in $T^*N_{\gamma}-0$ to the conic neighborhood
(0.5) of $T^*_{+}S^1_L$ in $T^*(S^1_L \times R^n)$ such that:

$${\cal D}:=W \sqrt{\Delta_{\psi}}W^{-1} \equiv \overline{\psi}({\cal R}, I_1, ...,I_n)\;\;
[ {\cal R} +\frac{ p_1(I_1,....,I_n)}{L{\cal R}}
+\frac{p_2(I_1,...,I_n)}{(L{\cal R})^2} + ...+\frac{p_{k+1}(I_1,\dots,I_n)}{(L{\cal R})^{k+1}}+\dots]$$
$$ \equiv D_s + \frac{1}{L}H_{\alpha} + \frac{\tilde{p}_1(I_1,\dots,I_n)}{L D_s} +
 \frac{\tilde{p}_2(I_1, \dots,I_n)}{(LD_s)^2}
+\dots+\frac{\tilde{p}_{k+1}(I_1,\dots,I_n)}{(L D_s)^{k+1}}+ \dots$$
where the numerators $p_j(I_1,...,I_n), \tilde{p}_j(I_1,...,I_n)$ are polynomials of degree
j+1 in the variables $I_1,...,I_n$, where $\overline{\psi}$ is microlocally
supported in (0.5), and where $W^{-1}$ denotes a microlocal inverse to $W$ in (0.5). The
kth remainder term lies in the space $\oplus_{j=o}^{k+2} O_{2(k+2-j)}\Psi^{1-j}$.}
\medskip

The QBNF coefficients will by definition be the coefficients of the monomials
in the classical action ($I_j$-) variables in
the complete Weyl symbols of the operators
$\tilde{p}_j(I_1,\dots,I_n)$.   As mentioned above,
the proof of Theorem B  gives an effective method for calculating them as  integrals over $\gamma$ of
FJ polynomials.

The asymptotic relation in the above expansion  may be viewed in either of two ways:
 First, as mentioned in the statement of the Theorem, the kth remainder is a sum of terms
in ${\cal A}$ of orders $1, 0, \dots, -(k+1)$ where the complete symbol of  the term
of order $1-j$ must vanish to order $2(k+2-j)$.   This characterization of the remainder
will play the key role in the calculation of the wave invariants, since terms in the
normal form with low pseudodifferential order or with high vanishing order make no
contribution to a given wave invariant.   On the other hand,
 it may be viewed as a
semi-classical asymptotic relation with ${\cal R}$ playing the role of semi-classical
parameter; thus the theorem gives
 a semi-classical expansion for ${\cal D}$ in terms of ${\cal R}$.   This point of view
comes up naturally in the theory of quasi-modes associated to $\gamma:$
 Indeed,  consider
the joint ${\cal A}$- eigenfunctions 
$$\phi_{kq}^o(s,y):= e_k(s)\otimes\gamma_q(y) \leqno(0.7)$$
  with 
$e_k(s):= e^{\frac{2 \pi}{L}iks}$  and with $\gamma_q$ 
$qth$ normalized Hermite function  ($q \in \Nb^n$).  The corresponding eigenvalues
of ${\cal D}$ 
then have the semi-classical expansions
$$\lambda_{kq} \equiv r_{kq} + \frac{p_1(q)}{r_{kq}} + \frac{p_2(q)}{r_{kq}^2} + ...\leqno(0.8)$$
where
$$r_{kq} = \frac{1}{L} (2 \pi k + \sum_{j=1}^n (q_j + \frac{1}{2}) \alpha_j)\leqno(0.9)$$
are the eigenvalues of ${\cal R}$.
Here the index $q$ is held fixed as $k \rightarrow \infty$.  We recognize in (0.8)  the
familiar form of the quasi-eigenvalues associated to $\gamma$ (cf. [B.B., ch.9]); hence the
intertwining operator $W$ is the operator taking the eigenfunctions (0.7) to quasi-modes of infinite
order for $\Delta$ at $\gamma.$

As will be seen in (\S4), Theorem B implies that the wave invariants of $\sqrt{\Delta}_{\psi}$ are the
same as the wave invariants of ${\cal D}$.
The second main step in the calculation of the wave invariants is then the use of the
non-commutative residue to connect the terms in the normal
form expansion with the terms in the singularity expansion for $Tr$ exp $it{\cal D}.$
The main point here is that 
$$a_{\gamma k} = res\ D_t^k e^{it\sqrt{\Delta_{\psi}}}: = Res_{s=0} TrD_t^k
 e^{it\sqrt{\Delta_{\psi}}} \sqrt{\Delta_{\psi}}^{-s},\leqno(0.10)$$ 
with $res$  invariant under conjugation by (microlocal) unitary operators, and
depending  on only a finite jet of the Laplacian near $\gamma.$  Hence it may
be calculated by conjugating to the normal form, and indeed will only depend on
a finite part of the normal form.  
Applying $D_t^k$ and formally exponentiating the terms of order $\leq -1$ 
in $D_s$ we get
$$ res \overline{\psi}(D_s,I_1,...,I_n) D_t^k e^{it{\cal D}}|_{t=L} =
 res \overline{\psi}(D_s ,I_1,...,I_n) e^{i  L D_s}e^{iH_{\alpha}} {\cal D}^k
 (I + iL \frac{\tilde{p}_1(I_1,...,
I_n)}{L D_s} + ...)\leqno(0.11)$$
which suggests that the wave coefficient $a_{k \gamma}$ is the regularized trace  of
the coefficient of ${D_s}^{-1}$ in (0.11).  This is not clear, even formally, since many of
the terms of negative order in $D_s$ in the exponent have overall
 order 1 as pseudodifferential operators; but
it will  prove to be the case.  Since $e^{i L D_s}= I $ on $\R/ L {\bf Z}$
 the Fourier Integral factor in (0.11) is just $e^{iH_{\alpha}}$.
Regarding the regularized traces of the coefficients, we note that   
$$Tr e^{iH_{\alpha}} =\sum_{q \in {\bf N}^n} e^{i\sum_{k=1}^n
 (q_k + \frac{1}{2})\alpha_k}\leqno(0.12)$$
is well-defined as the tempered distribution 
$$T(\alpha) =\Pi_{k=1}^n \frac{e^{\frac{i}{2}(\alpha_k+i0)}}{(1 - e^{i(\alpha_k+i0)})}
 \leqno(0.13)$$
on $\R^n_{\alpha}.$  Since its singular support is the union of the hyperplanes $\Z_{km}:=
\{(\alpha_1,\dots,\alpha_n) \in \R^n: \alpha_k = 2 \pi m\}, T$ has
smooth localization to neighborhoods of  
  Floquet exponents $(\alpha_1,\dots,\alpha_n)$ which are independent
of $\pi$ over ${\bf Q}$.  Hence the regularized trace is simply the evaluation of the
distribution trace at a regular point.
Similarly, the coefficient of ${D_s}^{-k-1}$ in (0.12) has a distribution trace of the form
$$\sum_{q \in {\bf N}^n} {\cal F}_{k,-1}(q_1 + \half,...,q_n + \half) e^{i \sum_{k=1}^n (q_k + \half)
\alpha_k} \leqno(0.14)$$ 
$$={\cal F}_{k,-1}(D_{\alpha_1} ,...,D_{\alpha_n} ) T(\alpha)$$
for a certain polynomial ${\cal F}_{k,-1}$.
Hence  this trace is also a locally smooth function in  neighborhoods of non-resonant exponents.   
\medskip

\noindent{\bf Theorem C}  {\it The wave invariants are given by 
$$a_{k \gamma} = {\cal F}_{k,-1}(D_{\alpha_1} ,...,D_{\alpha_n} ) T(\alpha).$$}

\medskip

The coefficients of the polynomials ${\cal F}_{k,-1}$ are evidently polynomials in the
QBNF coefficients and the differentiation process produces polynomials in the $\beta_j$'s.
Combined with Theorem B and a dimensional analysis, this proves Theorem A.

The proofs of Theorems A-C also lead to a somewhat simpler proof of Guillemin's inverse
theorem that the classical (in fact the full quantum) normal form is determined by
the wave trace invariants for all the iterates of $\gamma.$  We  only sketch the
proof here, assuming the reader's familiarity with the original proof of Guillemin in [G.2].
   The key point is to focus on the  Floquet invariants $\beta_j:=(1 - e^{\theta_j})^{-1}$
 for all the iterates $\gamma^m$
of $\gamma$, that is the residues (0.10) for $t=L, 2L, 3L,\dots.$  It follows from the
calculations in Theorems A-C that the wave invariants  are polynomials in the $\beta_j$'s:
more precisely, for each m,
 $a_{k \gamma^m}$ is the special value at $\theta_j = m\alpha_j$ of the fixed polynomial
$I_{\gamma;k}$ in the variables $\beta_j$.  Under the irrationality condition
above, the points $(e^{im\alpha_1},\dots,e^{im\alpha_n})$ form a dense set on the torus,
and hence the special values at these points
 determine the entire polynomial. The coefficients of the kth polynomial  $I_{\gamma; k}$
 are therefore determined by the wave invariants
for $\gamma, \gamma^2,\dots$.  By studying the relation of the coefficients of
$I_{\gamma; k}$ to the normal form invariants, Guillemin proves that all of the latter can be
determined from  the former.

Although Theorems A-C are only proved here under the hypothesis that $\gamma$ is 
non-degenerate elliptic, they have analogues for  hyperbolic and mixed
hyperbolic-elliptic geodesics, which we plan to describe in a future article [Z.3].
We also note that for closed geodesics possessing neighborhoods in
which the metric has no pairs of conjugate points for $t \leq L$, the wave invariants
can be calculated directly from a Hadamard parametrix [D] [Z.1,5].  Since a sufficiently large number of
iterates of an elliptic closed geodesics will
always contain pairs of conjugate points , but a small iterate may contain none, the 
calculation here and in [Z.1] overlap but are independent.  The calculations of [Z.1] also
apply to hyperbolic geodesics without pairs of conjugate points, showing that the form
of the wave invariants is essentially the same for the hyperbolic and elliptic cases.  However
the form resulting from the Hadamard parametrix is not immediately that of FJF polynomials,
 and it takes  considerable manipulation to show that
the formulae given here and in that paper agree.
 In the opposite extreme of Zoll manifolds, all of whose geodesics are closed and of completely
degenerate elliptic type, the
wave invariants are  calculated in [Z3,4] by yet another method. 
  
The organization of this paper is as follows:

\S 1: The models

\S 2: Semi-classical normal form of the Laplacian

\S 3: Normal form: Proof of Theorem B

\S4: Residues and wave invariants: Proof of Theorem C

\S5: Local formulae for the residues: Proof of Theorem A

\S6: Quantum Birkhoff normal form coefficients

\S7: Explicit formulae  in dimension 2

\S8: Appendix: The classical Birkhoff normal form

\S9: Index of Notation

The author wishes to thank S.Graffi for discussions of quantum Birkhoff normal forms
during a visit to the University of Bologna in June, 1994, where this work was begun.
  He also wishes to thank Y.Colin de Verdiere for his remarks on a preliminary version 
of the results which were presented during a visit to the Institut Fourier in January 1995.
  The final version was completed after the author received a copy of the article [G.2]
of Guillemin, and has benefited a good deal from the discussion there of normal forms.

\section{The models}

As mentioned above, the calculation of the wave invariants associated to a closed geodesic
$\gamma$ of a Riemannian manifold $(M,g)$ can be transplanted to the normal bundle
$N_{\gamma}$ by means of the exponential map.  Thus, the model space is the cylinder
$S^1_L \times \R^n$, where as above $S^1_L = \R / L\Z$. 
 In this section we first collect together some basic formulae and facts concerning the
``Hermite package" on this model space: that is, those aspects of analysis which come
from the representation theory of the Heisenberg and metaplectic algebras.   
We then
transfer the Hermite package to $N_{\gamma}$ in a way particularly well-adapted to
the metric along $\gamma$. 
\medskip

\noindent{\bf \S1.1: The model: ${\cal H} = H^2(S^1_L) \otimes L^2(\R^n)$}
\medskip

Since we are only concerned with the conic neighborhood
(0.2) of $\R^{+} \gamma$, we only consider the
 positive part $T^*_{+}S^1_L \times T^*\R^n$
and its quantum analogue the Hardy space $H^2(S_L^1)\otimes L^2(\R^n)$ with $k\geq 0.$

On the phase space level, the model is roughly $T^*(S^1_L \times \R^n)$. More 
precisely it is the cone (0.1) in the natural symplectic coordinates
$(s,\sigma,y,\eta)$. Since (0.1) is a conic neighborhood
 of $\R^{+} \gamma$, we will view it as a subcone of the
 positive part $T^*_{+}S^1_L \times T^*\R^n$ ($\sigma >0$).
 On the Hilbert space level the model is then ${\cal H}:= H^2(S^1_L) \otimes
L^2(\R^n)$, where $H^2(S^1_L)$ is the Hardy space, or more precisely
  the range ${\cal H}_{\psi}$ of the microlocal cutoff $\psi$ of the
introduction; generally we omit the subscript unless we need to emphasize the
role of $\psi$.
 We now introduce some distinguished algebras of operators on the
 model space.

First is 
the (complexified) Heiseberg algebra ${\bf h}_n\otimes \C$, which will be identified
with its usual realization on $L^2(\R^n)$.  It is then generated by the elements
$y_j =$ ``multiplication by $y_j$" and by $D_{y_j}=\frac{\partial}{i\partial y_j}$, or
equivalently by the
 annihilation, resp. creation, operators
$$A_j := y_j + iD_{y_j}  \;\;\;\;\;\;\;\; A_j^* = y_j - i D_{y_j}$$
which satisfy the commutation relations
$$[A_j, A_k] = [A_j^*, A_k^*] = 0 \;\;\;\;\;\; [A_j, A_k^*] = 2\delta_{ij} I.$$
The enveloping algebra of the Heisenberg algebra
$${\cal E} := <Y_1,...,Y_n, D_1,...,D_n> \leqno(1.1.1)$$ 
 is 
the algebra of partial differential operators on $\R^n$ with polynomial coefficients. We
let
${\cal E}^n$ denote the subspace of polynomials of degree n in the variables $y_j,D_{y_j}.$
In the usual isotropic Weyl algebra ${\cal W}^*$ of pseudo-differential operators on $\R^n$, the
operators $y_j, D_{y_j}$ are given the order $\frac{1}{2}$, so that
$${\cal E}^n \subset {\cal W}^{n/2}\leqno(1.1.1b)$$
$$[{\cal E}^m, {\cal E}^n] \subset {\cal E}^{m + n - 2}. $$
The symplectic algebra $sp(n,\C)$ is represented in ${\cal E}^2$ by 
homogeneous quadratic polynomials in $Y_j, D_j$, and a maximal abelian subalgebra of it
is spanned by the  harmonic oscillators (0.2).
We denote by 
$${\cal I}:= <I_1,...,I_n> \leqno(1.1.2a)$$
the (maximal abelian) subalgebra they generate in ${\cal W}$, with ${\cal I}^k:= {\cal I} \cap {\cal W}^k$,
and by
$${\cal P}_{{\cal I}} = {\cal I} \cap {\cal E}\leqno(1.1.2b)$$
the subalgebra of polynomials in the generators (0.2)), 
with ${\cal P}^k_{{\cal I}}$ the
space of polynomials of degree k.

The full pseudo-differential algebra on $S^1_L \times \R^n$ is the doubly filtered
algebra
$$\Psi^{**}(S^1_L \times \R^n) \equiv \Psi^*(S^1_L) \otimes {\cal W}^*$$ with
$$\Psi^{mn}(S^1_L \times \R^n) \equiv \Psi^m(S^1_L) \otimes {\cal W}^n.$$  A maximal abelian
subalgebra of it is given by 
$${\cal A}:= < D_s, I_1,\dots,I_n> = <{\cal R}, I_1,\dots, I_n> \leqno(1.1.3)$$
where ${\cal R}$ is the distinguished element (0.3).  It inherits
a double filtration ${\cal A}^{mn}$. As above, our interest is really in
the microlocalization of (1.1.3)    to the cone (0.2), 
i.e. the operators in (1.1.3)
will only be used in composition with the microlocal cut-off $\overline{\psi}({\cal
R}, I_1,...,I_n)$. In this cone $D_s$ and ${\cal R}$ are elliptic.  Hence 
  the subalgebra
$$<{\cal R}> \otimes {\cal P}_{{\cal I}} \leqno(1.1.4)$$
 of pseudo-differential symbols in ${\cal R}$ with coefficients
 in ${\Pp}_{\I}$ is well defined.  

 An orthonormal basis of $L^2(\R^n)$ of joint eigenfunctions of ${\cal I}$ is
 provided by the Hermite
functions $\gamma_{ q}, q \in \Nb^n.$ Here, $\gamma_o$  is the
Gaussian $\gamma_o(y)=\gamma_{iI}(y):= e^{-\half |y|^2}$.  It is
 the    unique ``vacuum state",  i.e. the state annihilated by the annihilation operators.
  The qth Hermite function is
then given by
$ \gamma_{ q}:= C_q A_1^{* q_1}...A_n^{*q_n} \gamma_o (q \in {\bf N}^n\}),$
with $C_q = (2 \pi)^{-n/2} (q!)^{-1/2}, q!=q_1!...q_n!.$  The notation
  ``$\gamma_{q}$" for ``Gaussian"
is standard in this context, see [F],  and should not be confused with the
 notation for closed geodesics.
An orthonormal basis of $H^2(S^1_L) \otimes L^2(\R^n)$ of  joint ${\cal A}$-
eigenfunctions is then furnished by
$$\phi_{kq}^o(s,y):= e_k(s)\otimes\gamma_{ q}(y),\;\;\;\;\;\;e_k(s):= e^{\frac{2 \pi}{L}iks}. \leqno(1.1.5)$$
  \medskip

\noindent{\bf \S1.2: The twisted model ${\cal H}_{\alpha}$}
\medskip

We now introduce a unitarily equivalent (twisted) version of the model, in which
the distinguished element ${\cal R}$ gets conjugated to $D_s.$  This will
eventually help to simplify the transport equations in \S2.  

The unitary equivalence will be given by conjugation with 
the unitary operator 
$$\mu(r_{\alpha}):= \int_{S^1_L}^{\oplus} \mu(r_{\alpha}(s))ds =
\int_{S^1_L}^{\oplus} e^{i \frac{s}{L} H_{\alpha}} ds\leqno(1.2.1a)$$
where $\mu$ is the metaplectic representation, and where $r_{\alpha}(s)$ is the
block diagonal orthogonal transformation on $\R^{2n}$ with blocks
$$r_{\alpha_j}(s):= \left( \begin{array}{ll} cos \alpha_j \frac{ s}{L} & sin \alpha_j
\frac{ s}{L} \\ -sin \alpha_j \frac{ s}{L} & cos \alpha_j \frac{ s}{L}
 \end{array} \right).\leqno(1.2.1b)$$
The direct integral here refers to the representation of $L^2(S^1)\otimes L^2(\R^n)$
as $\int_{S^1_L} L^2(\R^n)ds$, that is,
$$\int_{S^1_L}^{\oplus}\mu(r_{\alpha}(s))ds f(s,y) = \mu(r_{\alpha}(s))f(s,y)$$
where the right side is the application the operator in the $y$-variables with $s$
fixed. As will be seen below, $\mu(r_{\alpha})$ conjugates ${\cal R}$ to $D_s$,
and commutes with $I_1,\dots, I_n$.  Hence it  preserves the algebra (1.1.3).

On the other hand, it does not preserve the Hilbert space ${\cal H}.$
Indeed, the elements $\phi^o_{kq}$ get transformed into the elements 
$$e_{kq}(s) \otimes \gamma_q(y),\;\;\;\; \;\;\;\;\;\;\;\;\;\;\;\;
 e_{kq}(s):= e^{ir_{kq}s} \leqno(1.2.2)$$
which are not periodic in $s$.  Rather they satisfy
$$e_{kq} \otimes \gamma (s+L,y)= e^{i \kappa_q} e_{kq}
 \otimes \gamma(s,y) \leqno(1.2.3a)$$
with
$$\kappa_q = \sum_{j=1}^n (q_j + \frac{1}{2}) \alpha_j.\leqno(1.2.3b)$$
The space $H^2(S^1_L)\otimes L^2(\R^n)$ thus gets taken to the space ${\cal H}_{\alpha}$
of elements of the form 
$$f(s,y) =\sum_{k=o}^{\infty}\sum_{q \in \Nb^n}
 \hat{f}(k,q) e_{kq} \otimes \gamma_q
\leqno(1.2.4)$$
with square summable coefficients.

We can better describe this Hilbert space (and its associated phase space) in the
language of `quantized mapping cylinders'.

 On the phase space level we have the
symplectic  map $r_{\alpha}(L)$ of $T^*\R^n$ of (1.2.1b), essentially the Poincare map
of our problem. As in [F.G][G.2] we can introduce its (homogeneous) symplectic
mapping cylinder $C_{r_{\alpha}(L)}$: namely, the quotient of $T^*\R \times T^*\R^n$
under the cylic group $<R_{\alpha}(L)>$ generated by the symplectic map
$$\tilde{r}_{\alpha}(L) :  T^*\R \times T^*\R^n \rightarrow T^*\R \times T^*\R^n,
\;\;\;\;\;\;\;\; \tilde{r}_{\alpha}(L) (s,\sigma, y,\eta):= (s + L, \sigma,
r_{\alpha}(L)(y,
\eta)). \leqno (1.2.5)$$
Note that the first return time is constant, which is consistent with [F.G][G.2,
2.10] since elements of $Sp(2n,\R)$ preserve the contact form $ds - \half(\eta dy -
yd\eta).$ Also note that the mapping cylinder can be untwisted via the symplectic map
$$R_{\alpha} : T^*\R \times T^*\R^n\rightarrow T^*\R \times T^*\R^n,\;\;\;\;\;\;\;
R_{\alpha} (s, \sigma, y, \eta):= (s, \sigma + \half \sum_j \alpha_j (y_j^2 + 
\eta_j^2),r_{\alpha}(s)(y,\eta)).\leqno (1.2.7)$$
Indeed, we have
$$R_{\alpha}  (s + L, \sigma, y, \eta) = \tilde{r}_{\alpha}(L) R_{\alpha}(s,\sigma,
y,\eta)\leqno(1.2.8)$$
so that $R_{\alpha}$ induces a symplectic equivalence $C_{r_{\alpha}(L)} \sim
T^*(S^1_L)\times R^*(\R^n).$ 
  
On the quantum level, the analogue of the mapping cylinder is the Hilbert space
${\cal H}_{\alpha}$ of functions on $\R \times \R^n$ satisfying
$$f(s + L, y) = \mu(r_{\alpha}(L)) f(s,y) \leqno (1.2.6)$$
and square integrable on $[0,L)\times \R^n$.  The intertwining operator
$\mu(r_{\alpha})$ is essentially the analogue of $R_{\alpha}$.  More precisely,
we have:
\medskip

\noindent{\bf (1.2.9) Proposition}{\it

(i) $ \mu(r_{\alpha})^* D_s \mu(r_{\alpha}) = {\cal R};$

(ii) $\mu(r_{\alpha}) \phi^o_{kq}(s,y) = e^{ir_{kq}s}\gamma_q;$

(iii) $\mu(r_{\alpha}) : {\cal H} \rightarrow {\cal H}_{\alpha}$.}
\medskip

\noindent{\bf Proof}:

(i)  Follows from the fact that $\mu$ takes the jth diagonal block
 $$\left(\begin{array}{ll} 
0 & 1 \\ -1 & 0  \end{array}\right)$$ to $I_j$ and hence that $\mu(r_{\alpha}(s)) =
e^{i \frac{s}{L} H_{\alpha}}$. 

(ii) Follows from (i) and the fact that $\gamma_q$ is an eigenfunction of eigenvalue
$q_j + \half$ of $I_j$.

(iii) Follows from (ii).\qed
\medskip 

Under conjugation by $\mu(r_{\alpha})$, the sub-algebra (1.1.4) goes over to
the algebra
$$<D_s> \otimes {\cal P}_{{\cal I}} \leqno(1.2.10)$$
of pseudodifferential symbols in $D_s$ with coefficients in the $I_j$'s.  As
usual, it is understood to be microlocalized to the cone (0.2). 

We now conjugate by a further unitary equivalence to transfer this Hermite package to
$N_{\gamma}$.  It will induce a new model  better adapted to the geometry of $(M,g)$
near $\gamma;$ we call it the adapted model.

\noindent{\bf \S1.3: The adapted model}

To define it, we must discuss the Jacobi equation along $\gamma.$  Let Y(s) be a vector
field along $\gamma,$ and as above let $Y(s) = \sum_{i=1}^{n} Y_i(s) e_i(s)$ be its expression
in terms of the parallel normal frame.  The Jacobi equation is then
$$ \frac{d^2}{ds^2} Y_i + \sum_{ij=1}^{n} K_{ij} Y_j = 0.$$
Let ${\cal J}^{\bot}_{\gamma}\otimes {\bf C}$ denote the space of complex orthogonal Jacobi fields
along $\gamma.$ Equipped with the symplectic form $\omega(X,Y):=
<X, \frac{DY}{ds}>-<\frac{DX}{ds},Y>,$
it is a symplectic vector space of dimension 2n (over ${\bf C}$). 
Here $\frac{D}{ds}$ denotes
covariant differentiation along $\gamma,$ and $< , >$ is the
 inner product defined by the
metric.
Now let $Z(s) := (Y(s), \frac{DY}{ds}).$ The linear Poincare
 map $P_{\gamma}$ is by defintion
the operator on ${\cal J}^{\bot}_{\gamma}\otimes {\bf C}$ given
 by $P_{\gamma} Z(t)=Z(t+L).$  We recall that
$\gamma$ is assumed non-degenerate elliptic, hence 
$$Spec(P_{\gamma}) \subset S^1 \;\;\;\;\;\; \pm1 \notin spec (P_{\gamma}). $$
Since $P_{\gamma} \in Sp(2n, {\bf R}),$ its eigenvalues come in complex conjugate pairs
$\{ e^{i \alpha_k}, e^{-i\alpha_k} \}.$ The exponents $\alpha_k$ (called Floquet)
 are only defined up multiples of $2 \pi;$ they will be normalized below as in [B.B,
(9.3.17)]. Then there exists a basis of complex 
eigenvectors $\{Y_1,...,Y_n\}$ satisfying
$$P_{\gamma} Y_k = e^{i \alpha_k} Y_k
 \;\;\;\;\;\; Y_j \in {\cal J} \otimes {\bf C} \leqno(1.3.1)$$
$$\omega(Y_k, Y_j)= \omega(\overline{Y}_j,
 \overline{Y}_k) = 0,\;\;\;\;\;\; \omega(Y_j, \overline{Y}_k)
=\delta_{jk}.$$
with one choice of eigenvalue from each complex conjugate pair.  Equivalently,
 the span of $\{Y_1, ..., Y_n\}$ defines a $P_{\gamma}$-invariant positive
 Lagrangean subspace of ${\cal J}_{\gamma}^{\bot}\otimes {\bf C}$.

Consider now the (modified) Wronskian matrix 
$$a_s:= \left( \begin{array}{ll} Im\dot{Y}(s)^* \;\;\;&ImY(s)^*\\
Re\dot{Y}(s)^*\;\;\;&Re Y(s)^*  \end{array} \right). \leqno(1.3.2)$$
Here, with a little abuse of notation, we denote the components of $Y(s)$, resp. $\frac{DY}{ds}$
relative to the  parallel normal frame $e(s)$ by $Y(s)$, resp. $\dot{Y}(s)$ and the notation $A^*$
refers to the adjoint of a matrix $A$.
From (1.3.1) we have that $a_s \in Sp(2n, {\bf R})$ , which is equivalent to
$$\overline{Y}\frac{DY}{ds} - \overline{\frac{DY}{ds}}Y = iI.$$

As above, we let $\mu$ denote the metaplectic representation. 
  Identifying $a_s$ with one of its two possible lifts from
$Sp(2n, \R)$ to $Mp(2n, R)$ we introduce the unitary operator
$$\mu(a):=\int_{\gamma}^{\oplus}\mu(a_s)ds \leqno(1.3.3)$$
on
$\int_{S^1_L}^{\oplus} L^2(\R^n) ds$. In other words,
$$\mu(a) f(s,y) = \mu(a_s)f(s,y)$$
where the operator on the right side acts in the $y$-variables. 
 Informally, we think of the range as the
Hilbert space $\int_{\gamma}^{\oplus}
L^2(N_{\gamma(s)})ds,$ and below we will describe it more completely
in terms of quantum mapping cylinders. 
We now use $\mu(a)$ to transfer the Hermite package to $N_{\gamma}$.  We begin with
the generators of the Heisenberg algebra, and set:

$$\begin{array}{ll} P_j := \mu(a)^{*} D_j \mu(a) & Q_j:=\mu(a)^{*} Y_j \mu(a) \\

\Lambda_j:=\mu(a)^{*}A_j\mu(a) & \Lambda_j^*=
\mu(a)^{*}A^*_j\mu(a) \end{array}.$$
  
 $$I_{\gamma j}:= \frac{1}{2}\Lambda_j \Lambda_j^*.$$

We will refer to the $\Lambda_j$'s, resp. $\Lambda_j^*$'s, as the adapted annihilation
and creation operators and to the operators $I_{\gamma j}$  
in the symplectic algebra as the adapted action operators.
These adapted operators play a key role in the study of quasi-modes associated to $\gamma$. To
establish the connection, we
 now verify that they coincide with the operators similarly denoted in [B.B., ch.9]. 

First, some remarks on notation.  For reasons that will become clearer below (\S 1.4-5), we change
the notation for the transverse coordinates from $y$ to $u$, which should be thought of as the 
rescaled coordinates $u = L^{-\frac{1}{2}} y$.  Objects in the adapted model will henceforth always be
expressed in terms of $u-$coordinates. For instance, multiplication by $u_j$ will be denoted simply by $u_j$
and differentiation in $u_j$ by $D_{u_j}$. Also,   to quote easily some
basic facts about the metaplectic representation from Folland [F], we will   
conform to the following `transposed' notation  for the remainder of this section: symbols of operators
will be denoted $\sigma(p,q)$ rather than $\sigma(q,p)$, and the corresponding
  Weyl pseudodifferential operator  will be
 denoted by $\sigma(D,x)$ (later we will also use the more standard
notation $\sigma^w(x,D)$). An element of $Sp(n,\R)$ will be denoted by ${\cal A}$ and
 its transpose by ${\cal A}^*$.  
\medskip

\noindent(1.3.4) {\bf Proposition} 
$$\begin{array}{ll}
\Lambda_j =\sum_{k=1}^n (i y_{jk} D_{u_k} -\frac{dy_{jk}}{ds} u_k) &
\Lambda_j^*=\sum_{k=1}^{n} (-i \overline{y}_{jk}D_{u_k} -
\overline{\frac{dy_{jk}}{ds}}u_k) \end{array} $$
\medskip

\noindent{\bf Proof}:   This follows from the metaplectic covariance of the Weyl calculus, that
is from the identity [F, Theorem 2.15]
$$(\sigma \cdot {\cal A})(D,x) = \mu({\cal A}^*) \sigma (D,x) \mu({\cal A^*})^{-1}.$$
If we set $a = {\cal A}^*$ we get (in an obvious vector notation)
$$\mu(a)^*(u_j + iD_{u_j})\mu(a) = i(Y_j\cdot D_{u}) - \dot{Y}_j \cdot u$$
(compare [BB, ch.9]). \qed

\medskip

As for the tangential operator, we have: 
\medskip

\noindent(1.3.5) {\bf Proposition} {\it The image ${\cal L}$ of $D_s$ under $\mu$ is given by:
  $${\cal L}:= \mu(a)^{*} D_s \mu(a) = D_s - \frac{1}{2}(\sum_{j=1}^n D_{u_j}^2 + \sum_{ij=1}^n K_{ij}(s)
u_iu_j).$$}
\medskip

\noindent{\bf Proof}:  The left side is equal to
$(D_s + \mu(a_s)^{*} D_s\mu( a_s)).$ 
  To evaluate the second term, we use that both ReY(s)
and ImY(s) are Jacobi fields, and that Jacobi's equation is equivalent to the linear system
$\frac{D}{ds}(Y,P) = J H (Y, P).$  Here, $P=\frac{DY}{ds},$ J is the standard 
complex structure on $R^{2n},$ and
$$H =\left (\begin{array}{ll} K & 0 \\ 0 & I \end{array} \right ) $$
where $K$ is the curvature matrix and $I$ is the identity matrix [B.B, (9.2.9)]. 
 Hence, the second term is $\frac{1}{i} d\mu (JH)$ with
$d\mu$ the derived metaplectic representation.   But $\frac{1}{i}d\mu(JH) =
1/2(\sum_{i=1}^{n} \partial_{u_i}^2 - \sum_{ij=1}^{n} K_{ij}(s)u_i u_j)$ [F]. \qed
\medskip

We now consider the appropriate Hilbert space (quantized mapping cylinder) in the
adapted model. We first note:
\medskip

\noindent(1.3.6) {\bf Proposition}{ 

(i) $$\mu(a^{-1}) \gamma_o(s,u):=U_o(s,u)  = (det Y(s))^{-1/2} e^{i\half \langle
\Gamma(s) u,u\rangle}$$ where $\Gamma(s) := \frac{dY}{ds} Y^{-1}.$

(ii) $$\mu(a^{-1})\gamma_q:=U_q = \Lambda_1^{q_1}...\Lambda_n^{q_n} U_o. $$}
\medskip

\noindent{\bf Proof}:
Let us recall the
action of the metaplectic representation on Gaussians [F, Ch. 4.5].  For $Z$ in the Siegel upper
half space (n x n complex matrices $Z= X + i Y$ with $Y >> 0$),
 we have the Gaussian
 $$\gamma_Z(x) :=e^{i \half <Zx,x>}$$ on $\R^n.$  The action of an element
${\cal A} \in Mp(2n, \R)$ on the Gaussian is given by  
$$\mu ({\cal A}^{* -1}) \gamma_Z = m({\cal A}, Z) \gamma_{\alpha({\cal A})Z}$$
where 
$$m({\cal A},Z)=det^{-1/2}(CZ+D), \alpha({\cal A})Z =(AZ+B)(CZ+D)^{-1}$$
for 
$${\cal A}= \left( \begin{array}{ll} A & B \\ C & D \end{array} \right )$$
 (see [F, 4.65].)  Writing 
$$ Y(s) = Im Y(s) i + Re Y(s), \;\;\;\;\;\;\;\Gamma(s) = (Im \dot{Y}(s) i + Re \dot{Y})(Im Y(s) i + Re
Y(s))^{-1}$$
we see that $ m({\cal A}, iI) \gamma_{\alpha({\cal A})iI}(u) = (det Y(s))^{-1/2} exp (\frac{i}{2} <\Gamma(s)
u,u>$  if 
$${\cal A}= \left( \begin{array}{ll} Im \dot{Y}(s) & Re \dot{Y} \\Im Y(s) &  ReY(s)\end{array} \right ).$$
The formula (i)  follows from ${\cal A} = a^*$ and (ii) is an immediate
consequence of (i).\qed

Consider now the periodicity properties of the above data under $s\rightarrow
s+L$.  We observe that $a_s$ fails to be periodic in $s$ for two
reasons: first, due to the holonomy
of the frame $e(s)$, and  secondly due to the monodromy of the Jacobi fields 
$Y(s)$.  Indeed, we have:
$$a_{(s+L)} = T a_{s} r_{\alpha}^{-1}(L) \leqno(1.3.7a)$$
with $r_{\alpha}$ as in \S1.2 and where $T$, the holonomy matrix,
 is the 2n by 2n block diagonal matrix
with equal diagonal blocks $t=(t_{ij})$ satisfying
$$e_i(L) = \sum_{j=1}^n t_{ij}e_j(0).$$  It is of
 course the lift to $T^*\R^n$ of a rotation
on the base. The two properties can be summarized by writing
$$\tilde{a}_{s+L} = T \tilde{a}_s  \leqno(1.3.7b)$$
with $\tilde{a}_s:=a_s r_{\alpha}(s)$.

As in \S1.2 we will reformulate (1.3.7a-b) in terms of quantum mapping cylinders.
First, we put
$$C^{\infty}_T(\R\times \R^n):= \{ f \in C^{\infty}(\R\times \R^n):
f(s+L,u) = \mu(T) f(s,u)\}\leqno (1.3.8)$$
and let ${\cal H}_T$ denote its closure with respect to the obvious inner product
over $[0,L)\times \R^n$.  Note that the metaplectic operator $\mu(T)$ is simply
$$\mu(T) f(u) = f(t^{-1}u)$$
and hence that  
$$C^{\infty}_T(\R\times \R^n)\sim C^{\infty} (N_{\gamma})$$
where the isomorphism is simply the pull-back by the exponential map defined by
the frame $e(s)$.  Thus:
\medskip

\noindent(1.3.8){\bf Proposition}\;{\it 
\medskip

(i) Let $P(s,D_s,u,D_u)$ be a partial differential operator on $N_{\gamma}$,
 expressed
in the coordinates $(s,u).$  Then
$$P(s+L,D_s,u,D_u) =\mu(T) P(s,D_s,u,D_u)\mu(T)^*  ;$$

(ii) The functions 
$$\mu(\tilde{a}_s)(\phi^o_{kq}) =\phi_{kq}:=e^{i r_{kq} s}U_q(s,u) $$
define a smooth orthonormal basis of ${\cal H}_T$}.  
\medskip

\noindent{\bf Proof}
\medskip

(i) It suffices to prove this when $P$ is a vector field given in the local
normal coordinates by $a_o(s,u)D_s + \sum_{j=1}^n
a_j(s,u)D_{u_j}.$   
 Since the metaplectic operator $\mu(T)$ corresponding to $T$ is the operator
$f(u) \rightarrow f(t^{-1} u),$ we have
$$\mu(T)P\mu(T))^* = a_o(s,t^{-1}u)D_s + \sum_{j=1}^n a_j(s, t^{-1}u) t_{jk}D_{u_k}.$$
On the other hand, the vector field is well-defined on $N_{\gamma}$ if and only if
$$a_o(s+L, t^{-1}u)D_s +\sum a_i(s+L, t^{-1}u)t_{ij}D_{u_j} = 
a_o(s,u)D_s +\sum a_j(s,u)D_{u_j}.$$

 (ii)  Clear, since by (1.3.7b)
 $\mu(\tilde{a})$ intertwines the model and the quantum mapping cylinder of
$\mu(T).$\qed
\medskip

Remark: Statement (ii) is equivalent to 
$$U_q(s + L, t^{-1} u) = e^{-i \kappa_q} U_q(s,u)$$  
(correcting the formula stated in [B.B., (9.3.25)].) 
\medskip

It follows that
 we may write a smooth function $f$ in the adapted model in the form
$$f(s,u) = \sum_{k=o}^{\infty} \sum_{q \in {\bf N}^n}
 \hat{f}(k,q) e^{ir_{kq}s}U_q(s,u).
 \leqno(1.3.9)$$
\medskip

\noindent{\bf \S1.4 Metric scaling and weights}
\medskip

As mentioned above, $\Delta$ and hence the wave invariants have well-defined weights under the metric rescaling
$g \rightarrow \epsilon^2 g$. Since the wave invariants will be expressed in terms of QBNF coefficients, it
is natural to ask how the latter scale.  The question is not really well-posed since the QBNF coefficients
are coefficients with repsect to Harmonic oscillators $\partial_{y_j}^2 + y_j^2$ whose scaling behaviour
depends on the choice of coordinates.  To amplify this point, we record how various metric objects
scale under metric re-scaling.

In the following table, $(s,y)$ denote the Fermi normal coordinates relative to $g$, $(s,u)$ denote
the scaled coordinates $(s, L^{-\half}y)$ and $p_u$ denotes the symplectic coordinates dual to $u$.
\bigskip

\begin{tabular}{r||l} $g$ & $\epsilon^2 g$ \\ \hline
 $L, \;\;\;\;\;\;(s,y)$ & $\epsilon L, \;\;\;\;\;\;(\epsilon s, \epsilon y)$\\ \hline
$\partial_s,\;\;\; \partial_{y_j},\;\;\;e_j$ & $\epsilon^{-1}\partial_s,\;\;\; \epsilon^{-1} \partial_{y_j},
\;\;\;\epsilon^{-1} e_j$\\ \hline
$ y_{ij} := g(Y_i,e_j)$ & $\epsilon^{\frac{1}{2}} y_{ij}
 = \epsilon^2 g( \epsilon^{-\frac{1}{2}} Y_j,\epsilon^{-1}
e_j)$ \\ \hline
$K_{ij} = g(R(\partial_s,e_i)\partial_s,e_j)$ & $\epsilon^{-2} K_{ij}$ \\ \hline
$u = L^{-\half}y \;\;\;\;\; p_u = L^{\half}\eta$ & $\epsilon^{\half}u \;\;\;\;\;\;\; \epsilon^{-\half}p_u$ \\
\hline
$\Lambda_j = \sum_{k=1}^n (iy_{jk}D_{u_k} - \dot{y}_{jk}u_k)$ &
$\sum_{k=1}^n (iy_{jk}D_{u_k} - \dot{y}_{jk}u_k)$ \\ \hline
$ {\cal L}:= D_s - \frac{1}{2}(\sum_{j=1}^n D_{u_j}^2 + \sum_{ij=1}^n K_{ij}(s)
u_iu_j)$ &$ \epsilon^{-1} [D_s - \frac{1}{2}(\sum_{j=1}^n D_{u_j}^2 + \sum_{ij=1}^n K_{ij}(s)
u_iu_j)]$\\ \hline

\end{tabular}
\bigskip

The entry $y_{ij} \rightarrow \epsilon^{\half} y_{ij} $ 
 follows from the scale invariance of the Jacobi equation  together with the
normalization condition 
$$g(Re Y_i, Im \dot{Y}_i) - g(Re \dot{Y}_i, Im Y_i) = Constant.$$
which implies that $Y_i \rightarrow \epsilon^{-\half}Y_i$.

We observe that the creation/annihilation operators, hence the harmonic oscillators $\Lambda_j^*\Lambda_j$,
 of the adapted
model are scale-invariant,  and that the distinguished element ${\cal L}$ has weight -1. These are the
desired scaling properties and we would like the basic and twisted models to possess them as well.  As
they stand, these models do not scale properly if we interpret the $(s,y)$-coordinates as Fermi normal
coordinates.  However, they do scale properly if we interpret the $y$ coordinates as weightless. 
\medskip

\noindent{\bf \S1.5 Scaled adapted model and intertwining operators}
\medskip

 To
avoid confusion, we now introduce  the weightless coordinates
 $x = L^{-1}y = L^{-\half}u, \xi = L \eta = L^{\half}p_u$ and henceforth use them exclusively for the
scaled adapted, basic and twisted models.  In the following table we record how the various objects
appear in the weightless coordinates.  We also record the various intertwining operators, since they
get altered when we used weightless coordinates.
For instance, intertwining by $\mu(r_{\alpha})$ above is weightless but
that by $\mu(a)$ is of mixed weight. 
\bigskip

\begin{tabular}{r|l} Adapted model & Scaled adapted  model \\ \hline \hline
$\Lambda_j$ & $ \sum_{k=1}^n (i y_{jk}
 L^{-\frac{1}{2}} D_{x_k} - L^{-\frac{1}{2}} \frac{d y_{jk}}{ds} x_k)$ \\ \hline
${\cal L}$ &$D_s - \frac{1}{2}(\sum_{j=1}^n L^{-1} D_{u_j}^2 + \sum_{ij=1}^n L K_{ij}u_i u_j)$ \\ \hline
 $U_o(s,u)$ & $U_{L o}(s,x)  = (det Y(s))^{-1/2} e^{ \frac {i}{L} \langle \Gamma(s) x,x
\rangle}|dx|^{\half}$\\ \hline 
\end{tabular}
\bigskip

The following  intertwining operators will arise in the construction of the normal form:
\bigskip

\begin{tabular}{r|l|l}  From * to * & Classical & Quantum \\ \hline
Ad.model to Sc.ad.model & $(u,p_u) \rightarrow  (L^{-\half} u, L^{\half}p_u) = (x,\xi)$ & 
$\mu(D_{L^{\half}})$ \\ \hline
Sc.Ad.Model to Tw.model & $(x,\xi) \rightarrow (L^{- \half}Re Y x + L^{-\half}Im Y  \xi,
L^{\half}Re \dot{Y}  x +L^{\half} Im \dot{Y}\xi )$ & $\mu(  D_{L^{\half}} \cdot a )$\\ \hline
\end{tabular}
\bigskip

Above, the notation $\mu(D_r)$ refers to the metaplectic (dilation) operator 
corresponding to the symplectic matrix
$$D_r := \left( \begin{array}{ll} r & 0 \\ 0 &  r^{- 1} \end{array} \right),\;\;\;\;\;
\mu(D_r) f(s,y):= f(s, r^{-1} y).$$

\section{Semi-classical normal form of the Laplacian }

We return now to $\sqrt{\Delta}$, which as in the
introduction will be identified with its transfer to $N_{\gamma}$ under the
exponential map.  The finite jets of this transfer are globally well- defined  on
$N_{\gamma}$, so we will often treat $\Delta$ as if it too were globally well-defined.
  Our purpose  is to define the
{\it semi-classically rescaled Laplacian} $\Delta_h$ and to put $\Delta_h$ into a 
{\it semi-classical normal form}.  This
is the crucial  preliminary step in putting $\Delta$ itself into normal form.

In view of \S1.5, there are two rescalings at hand: the semi-classical rescaling $\Delta_h$ and the metric
$L$-rescaling above.  The two rescalings have quite distinct origins, so we have kept them separate.

To motivate the rescalings and the emergence of semi-classical asymptotics, let us recall
that the quasi-modes associated to $\gamma$ have the form
$$\Phi_{kq}(s,\sqrt{r_{kq}}y) =e^{ir_{kq}s} \sum_{j=0}^{\infty}
r_{kq}^{-\frac{j}{2}} U_q^{\frac{j}{2}}(s, \sqrt{r_{kq}}y,r_{kq}^{-1})\leqno(2.1)$$
with $U_q^o = U_q$ (see [B.B]). 
 The intertwining operator $W_{\gamma}$ to the normal form is then
the operator defined by the equations
$$W_{\gamma} \phi_{kq}(s,y) = \Phi_{kq}(s, \sqrt{r_{kq}}y).\leqno(2.2)$$
 The higher order terms $U^{\frac{j}{2}}$,  hence $W_{\gamma}$, are determined by the conditions
$$\Delta_y e^{ir_{kq}s} U_{kq}(s, \sqrt{r_{kq}}y,r_{kq}^{-1}) \sim
\lambda_{kq}e^{ir_{kq}s}
 U_{kq}(s, \sqrt{r_{kq}}y,r_{kq}^{-1})\leqno(2.3)$$
with $\lambda_{kq}$ given by (0.9).

Now write
$$r_{kq} = \frac{1}{h_{kq} L}, \;\;\;\;\;\;\;\;h_{kq}:= (2\pi k + \sum_{j=1}^n (q_j +
\frac{1}{2}\alpha_j))^{-1} $$
 so that the Planck constants $h_{kq}$  are metric-independent.
  In the scaled Fermi coordinates $u$ of \S1.3-5, the quasi-modes have then the form 
$$\Phi_{kq}(s, h^{-\frac{1}{2}}u) = e^{i\frac{ s}{h_{kq} L}} U_{kq}(s,
\sqrt{h_{kq}}^{-\half}u,h_{kq})\leqno(2.1')$$ and the eigenvalue problem (2.3) becomes
$$\Delta_u e^{\frac{i}{h_{kq}L}s}
U(s, h_{kq}^{-\half}u,h_{kq}) = \lambda(h_{kq})
e^{\frac{i}{h_{kq}L}s}U(s, h_{kq}^{-\half}u,h_{kq}). \leqno(2.3')$$
We are thus led to study the asymptotic eigenvalue problem
$$\Delta_u e^{\frac{i}{hL}s}U(s, h^{-\half} u,h) = \lambda(h) e^{\frac{i}{hL}s}U(s,
h^{-\half}u,h)\leqno(2.4)$$ on $C^{\infty}(\R^1\times \R^n)$ with $U(s, u, h)$ and
$\lambda(h)$ asympotic series in $h$.  Since $\Delta_u$ comes from an operator on
$N_{\gamma}$,
 the eigenvalue problem is taking place
on the quantized mapping cylinder ${\cal H}_T$ (\S1.3) of the adapted model.

We also note that the local expression for the Laplacian in the $(s,u)$ coordinates  is the same as in
the $(s,y)$ (Fermi) coordinates:
$$\Delta_u = \frac{1}{\sqrt{g}} \sum_{ij=o}^n \partial_{u_i} g^{ij}\sqrt{g} \partial_{u_j} = \Delta_y|_{
u_j:= L^{-\frac{1}{2}}y_j}$$
 since $g^{ij} \rightarrow L^{-1} g^{ij}$ and $\partial_{u_i} \rightarrow L^{\frac{1}{2}}\partial_{y_i}.$ 
Momentarily we are going to rescale the coordinates again to the weightless $x$-coordinates of \S1.5, and
again the Laplacian will be given by the usual expression.
Hence it will be a simple matter to pass back and forth between the $(s,y), (s,u)$ and $(s,x)$ expressions.

It is natural at this point to introduce the  unitary operators $T_h$ and $M_h$ on
${\cal H}_T$ or equivalently on the 1/2-density version 
$L^2_T(\R^1\times \R^n, \Omega_{1/2})$ given by
$$T_h (f(s,u)|ds|^{1/2}|du|^{1/2}):= h^{-n/2} f(s, h^{-\half}u) |ds|^{1/2}|du|^{1/2} \leqno (2.5a)$$
$$M_h(f(s,u)|ds|^{1/2}|du|^{1/2}) := e^{\frac{i}{hL}s} f(s,y)|ds|^{1/2}|du|^{1/2} \leqno (2.5b).$$
We easily see that:
$$  T_h^* D_{u_j} T_h = h^{-\half} D_{ u_j} \leqno (2.6)$$
$$ T_h^* u_i T_h = h^{\half} u_i$$
$$ M_h^*D_s M_h =((hL)^{-1} + D_s)$$
$$ [M_h, u_i]=[M_h, D_{ u_i}]= [M_h, T_h] =
 [T_h, D_s]=0.$$

\noindent{\bf Definition} {\it The rescaling of an operator $A_u = a(s,D_s,u,D_u)$ of
the adapted model  is the operator
$$A_h := T_h^*M_h^*AT_hM_h\leqno(2.7)$$}

We observe that the operation of rescaling is weightless.
In particular, the rescaled Laplacian $\Delta_h$  in the sense
of (2.7) is of weight -2.  To calculate it, we  first note that the (1/2-density) Laplacian in scaled
Fermi normal coordinates is given by the expression 
  $$\Delta_u =-( J^{-1} \partial_{ s} J g^{oo} \partial_{s} + \sum_{ij=1}^n
J^{-1}\partial_ {u_i} g^{ij}J \partial_{ u_j})\leqno(2.8)$$ 
where $\partial_{x} := \frac{\partial}{\partial x},$ and where
$ J=J(s,u)= \sqrt{g}$ is the volume density in these coordinates.
 To obtain a self-adjoint operator
with respect to the Lesbegue density  $|ds| |du|$,  we
 replace $\Delta$ by the unitarily equivalent  1/2-density Laplacian
 $$\Delta_{1/2} := J^{1/2} \Delta J^{-1/2},$$ 
which can be written  in the form:
$$-\Delta_{1/2} = J^{-1/2}\partial_s g^{oo}J \partial_s J^{-1/2}
+\sum_{ij =1}^{n} J^{-1/2}\partial_{u_i} g^{ij} J \partial_{ u_j}
J^{-1/2}\leqno(2.9)$$
$$\equiv g^{oo}\partial_s^2 + \Gamma^o \partial_s + 
 \sum_{ij=1}^n g^{ij} \partial_{u_i}\partial_{u_j} + \sum_{i=1}^{n} \Gamma^{i}
\partial_{u_i} + \sigma_o.$$
 From now on, we will only use $\Delta_{\half}$ and denote it simply by $\Delta.$

We then have:
$$-M_h^* \Delta M_h = -(hL)^{-2}g^{oo} + 2i(hL)^{-1}g^{oo} \partial_s + i(hL)^{-1}\Gamma^o +
\Delta \leqno (2.10)$$
Conjugation with $T_h$ then gives
$$-\Delta_{h} =  -(hL)^{-2} g^{oo}_{[h]} + 2i(hL)^{-1}g^{oo}_{[h]}\partial_s + i(hL)^{-1}\Gamma^o_{[h]}+
 h^{-1}( \sum_{ij=1}^n g^{ij}_{[h]}\partial_{u_i}\partial_{u_j}) + h^{-\half}(\sum_{i=1}^{n} \Gamma^{i}_{[h]}
\partial_{u_i}) + (\sigma)_{[h]},\leqno(2.11)$$
 the subscript $[h]$ indicating to dilate the coefficients of the operator in the form,
$f_h(s, u):=f(s, h^{\half} u).$  

Expanding the coefficients in Taylor series at $h=0$, we obtain the asymptotic expansion
$$\Delta_h \sim \sum_{m=0}^{\infty} h^{(-2 +m/2)}{\cal L}_{2-m/2} \leqno (2.12)$$
where ${\cal L}_2 = L^{-2},$ ${\cal L}_{3/2}=0$ and where
$${\cal  L}_1 = 2 L^{-1}[i  \frac{\partial}{\partial s} + \half \{\sum_{j=1}^{n}
\partial_{u_j}^2 -
\sum_{ij=1}^{n} K_{ij}(s) u_i u_j\}] \leqno (2.13).$$
We observe that ${\cal L}_1$ is $2 L^{-1}$ times the distinguished element ${\cal L}$ of the adapted model,
and that ${\cal L}$ has weight $-1$ when the $u$-variables are given their natural
weights $1/2$. 

As discussed in \S1. 4-5, it will be helpful to rescale the variables once again to make them weightless.
Hence we change variables to $x = L^{-\half} u$ and rewrite $\Delta_h$ and the ${\cal L}_{2-\frac{n}{2}}$'s
in terms of the $x$-variables. For instance we will henceforth write ${\cal L}$ in the form: 
$${\cal L} = i  \frac{\partial}{\partial s} + \half [\sum_{j=1}^{n}L^{-1} \partial_{x_j}^2 -
\sum_{ij=1}^{n} L K_{ij}(s) x_i x_j].$$

 We further note that
the operators ${\cal L}_{2-m/2}$ now satisfy the periodicity condition (1.3.8ii) : Indeed,
as noted above, $\Delta$ has this property when expressed in normal coordinates, and the various 
transversal rescalings and
the conjugations by
$T_h,M_h$  preserve it. To indicate that an operator has this periodicity
property and hence acts on ${\cal H}_T$,
 we will subscript the appropriate spaces of operators with a ``T." 

  From (2.11) we see that the terms in ${\cal L}_{2 - \frac{m}{2}}$
are of the form 
$$\begin{array}{lll}x^m & x^{m-2} &x^{m-4}\\x^{m-2}D^2_x &x^{m-3}D_x \\
x^{m-2}D_s &x^{m-4}D_s &x^{m-4}D_s^2\end{array}$$
hence 
 $${\cal L}_{2 -\frac{m}{2}} \in \Psi^2(S^1)_T\otimes{\cal E}^{m-4}_{\epsilon}
+\Psi^1(S^1)_T\otimes{\cal E}^{m-2}_{\epsilon}
+\Psi^o(S^1)_T\times{\cal E}^m_{\epsilon}\leqno(2.14a)$$
where the subscript $\epsilon$ indicates that the Weyl symbol is a polynomial with
the parity of $m$.  Moreover, in the $s$ variable, ${\cal L}_{2 -\frac{m}{2}}$ is a polynomial
in $D_s$ of degree at most two, so we can refine (2.14a) to the statement
$${\cal L}_{2 -\frac{m}{2}} \in C^{\infty}_T(S^1_L,{\cal E}^{m-4}_{\epsilon})D_s^2
+C^{\infty}_T(S^1_L,{\cal E}^{m-2}_{\epsilon})D_s
+C^{\infty}_T(S^1_L,{\cal E}^m_{\epsilon})\leqno(2.14b)$$

Comparing with (1.3.4) we see that
$${\cal L} = \mu({\cal A}_L^*) D_s \mu({\cal A}_L^*)^{-1}$$
where ${\cal A}_L$ is the weightless Wronskian matrix $D_{L^{\half}} a^*$, that is
$${\cal A}_{L} := \left ( \begin{array}{ll} L^{\half} Im \dot{Y} & L^{\half} Re \dot{Y} \\
L^{-\half} Im Y & L^{-\half} Re Y \end{array} \right ).$$
This  motivates the conjugation of (2.11) to the (untwisted) model. We therefore put
$${\cal D}_h =\mu({\cal A}_L^*)^{-1} \Delta_h \mu({\cal A}_L^*)$$
which has the asymptotic expansion
$${\cal D}_h \sim \sum_{m=o}^{\infty} h^{(-2 + \frac{m}{2})}
{\cal D}_{2 -\frac{m}{2}} \leqno(2.15)$$
  with  ${\cal D}_2 = I, {\cal D}_{\frac{3}{2}}=0,
{\cal D}_1=D_s$. 
Conjugation by $\mu({\cal A}_L^*)$  preserves weights, homogeneity and parity in the variables $(x,D_x)$. It
also transforms $D_s$ into $D_s$ plus a term quadratic in $(x,D_x)$. 
  Hence we find easily that ${\cal D}_{2 - \frac{m}{2}}$ has weight -1 for each m and that
$${\cal D}_{2 -\frac{m}{2}} \in \Psi^2(\R)\otimes{\cal E}^{m-4}_{\epsilon} 
 +\Psi^1(\R^)\otimes{\cal E}^{m-2}_{\epsilon}
+\Psi^o(\R)\otimes{\cal E}^m_{\epsilon}\leqno(2.16a)$$
or, analogously to (2.14b),
$${\cal D}_{2 -\frac{m}{2}} \in C^{\infty}(\R,{\cal E}^{m-4}_{\epsilon})D_s^2 
 +C^{\infty}(\R,{\cal E}^{m-2}_{\epsilon})D_s
+C^{\infty}(\R,{\cal E}^m_{\epsilon})\leqno(2.16b)$$

Of course, conjugation with $\mu({\cal A}_L^*)$ also alters the periodicity
 property of the terms
${\cal D}_{2 -\frac{m}{2}}$.  From (1.3.7a) and (1.3.8) we see in fact
that they transform like operators in the twisted model, i.e. on the quantum
mapping cylinder of $r_{\alpha}(L).$  More precisely, we have
$${\cal D}_{2 -\frac{m}{2}}|_{s+L} = \mu({\cal A}_L^*)^{-1}({\cal L}_{2
-\frac{m}{2}}\mu({\cal A}_L^*)|_{(s+L)}) =\mu(r_{\alpha}(L))\mu({\cal A}_L^*)^{-1}\mu(T)^* {\cal L}_{2
-\frac{m}{2}}|_{(s+L)}\mu(T)\mu({\cal A}_L^*)\mu(r_{\alpha}(L))^*\leqno(2.17a)$$
$$=\mu(r_{\alpha}(L))\mu({\cal A}_L^*)^*{\cal L}_{2 - \frac{m}{2}}\mu({\cal A}_L^*)\mu(r_{\alpha}(L)^*
= \mu(r_{\alpha}(L)){\cal D}_{2 -\frac{m}{2}}\mu(r_{\alpha}(L))^*.$$
Equivalently, in terms of the matrix elements in the basis of Hermite functions,
we have
$$\langle{\cal D}_j \gamma_q, \gamma_r\rangle (s+L)=e^{-i(\kappa_q - \kappa_r)s}
\langle{\cal D}_j \gamma_q, \gamma_r\rangle (s) \leqno(2.17b).$$  To indicate 
that these operators act on ${\cal H}_{\alpha}$
 we henceforth subscript the appropriate spaces of operators  with
an $\alpha$.

To render these terms  periodic in $s$, we have to conjugate further
to the (untwisted) model under
$\mu(r_{\alpha})$.  We record the resulting expressions, since they will be
 used later on. In the notation $\tilde{{\cal A}}_L^*(s) := {\cal A}^*_L(s) \cdot r_{\alpha}(s)$
 the principal operator becomes,  by Proposition (1.3.8),
$$\mu(\tilde{{\cal A}}_L^*) {\cal L} \mu(\tilde{{\cal A}})^{-1} = {\cal R} \leqno(2.18)$$
 
 Since conjugation
by $\mu(r_{\alpha})$ also preserves weight, homogeneity and parity,
we further have:
$${\cal R}_h := \mu(\tilde{a})\Delta_h\mu(\tilde{a})^*\sim \sum_{m=o}^{\infty} h^{(-2 + \frac{m}{2})}
{\cal R}_{2 -\frac{m}{2}} \leqno(2.19)$$
with ${\cal R}_2 = I, {\cal R}_{\frac{3}{2}}=0, {\cal R}_1 = {\cal R}$ 
and with all coefficients of weight -2 and periodic, that is,
$${\cal R}_{2 -\frac{m}{2}} \in \Psi^2(S^1)\otimes{\cal E}^{m-4}_{\epsilon} 
 +\Psi^1(S^1)\otimes{\cal E}^{m-2}_{\epsilon}
+\Psi^o(S^1)\otimes{\cal E}^m_{\epsilon}\leqno(2.20a)$$
or, analogously to (2.14b),
$${\cal R}_{2 -\frac{m}{2}} \in C^{\infty}(S^1_L,{\cal E}^{m-4}_{\epsilon}){\cal R}^2 
 +C^{\infty}(S^1_L,{\cal E}^{m-2}_{\epsilon}){\cal R}
+C^{\infty}(S^1_L,{\cal E}^m_{\epsilon})\leqno(2.20b).$$

 Our aim is now to put $\Delta$, or $\Delta_h$ for certain $h$,
into a semi-classical normal form.  
This normal form will be, at first, only a formal normal form for $\Delta_h$ as
a formal $h-$pseudo-differential operator.  Here,
a formal $h$-pseudodifferential operator of order m on $\R^N$ is the Weyl
quantization
$$a^w(x, hD; h) u(x) =(2 \pi h)^{-n} \int\int e^{\frac{i}{h}\langle x-y, \xi \rangle}
a(\frac{1}{2}(x+y), \xi; h) u(y) dy d\xi  \leqno(2.21a)$$
of an amplitude $a$ belonging to the space $S_{cl}^m(T^*\R^N)$ of asymptotic
sums
$$a(x, \xi, h) \sim h^{-m}\sum_{j+o}^{\infty} a_j(x, \xi) h^j   \leqno(2.21b)$$
with $a_j \in C^{\infty}(R^{2N}).$  Such operators form an algebra 
$\Psi^*_h(\R^N)$ under composition, with $\Psi^m_h$ the subspace of mth order elements.
(See [Sj] for further background and references).  We will also be concerned with
the slightly different situation of $h$- pseudodifferential operators on $S^1 \times \R^n$,
which are defined similarly using local coordinates.
Combining the $h-$ filtration with the previous filtrations of $\Psi^*(\R)\otimes{\cal W}^*
(\R^n)$ we get the triply filtered algebra
$$\Psi_h^{(*,*,*)}(\R\times \R^n)$$
$$\Psi_h^{(k,m,\frac{n}{2})}:= h^k\Psi^m(\R)\otimes{\cal W}^{\frac{n}{2}},$$
with $k+m+\frac{n}{2}$ the total order of an element (and similarly with $S^1$ replacing
$\R$.)

 The following lemma will prepare for the normal form. We state it in terms of the ${\cal R}$-operators
of the model since the periodicity properties are simplest there.  
The notation ``$A|_o$" will be used for the
restriction  of an operator $A \in \Psi^*(S^1 \times \R^n)$ of the model to elements
of ${\cal R}$-weight zero.  Equivalently,
 after conjugation by $\mu(r_{\alpha})$ to the twisted model,
 to elements of ${D_s}$-weight zero , that is, to functions
independent of $s$ in ${\cal H}_{\alpha}$. If we write the latter $A$ in the form
$A_2 D_s^2 + A_1 D_s + A_o$, then $A|_o$ is  $A_o|_o$ (restricted to weight 0
elements).
We have:
\medskip

\noindent(2.22){\bf Lemma} {\it There exists an $L$-dependent $h$-pseudodifferential operator
$W_h=W_h(s, x, D_x)$ on $L^2(S^1_L \times \R^n)$ such that,
for each $ s\in S^1_L$, $$W_h(s,x,D_x): L^2(\R^n)\rightarrow L^2(\R^n) $$
is unitary, and such that
$$W_h^* {\cal R}_h W_h \sim - h^{-2}L^{-2} + 2 h^{-1}L^{-1}{\cal R} + \sum_{j=0}^{\infty} h^{\frac{j}{2}}
 {\cal R}^{\infty}_{2-\frac{j}{2}}(s,D_s,y,D_y)$$
where
\medskip

(i) $ {\cal R}^{\infty}_{2-\frac{j}{2}}(s, D_s, x, D_x) =
  {\cal R}^{\infty,2}_{2-\frac{j}{2}} {\cal R}^2 +
   {\cal R}^{\infty,1}_{2-\frac{j}{2}} {\cal R} + {\cal R}^{\infty,o}_{2-\frac{j}{2}},$ with
${\cal R}^{\infty,k}_{2-\frac{j}{2}} \in C^{\infty}(S^1_L, {\cal E}_{\epsilon}^{j -2k});$
\medskip

(ii) $ {\cal R}^{\infty}_{2-j}(s, D_s, x, D_x)|_o = {\cal R}^{\infty,o}_{2-j}(s, x, D_x)|_o =
  f_j(I_1,...,I_n)|_o$ for certain polynomials $f_j$
 of degree j+2 on $\R^n,$ i.e.
 $f_j(I_1,...,I_n) \in {\cal P}^{j+2}_{\cal I}$
\medskip

(iii) $ {\cal R}^{\infty}_{2-\frac{2k+1}{2}}
(s, D_s, x, D_x)|_o = {\cal R}^{\infty, o}_{2-\frac{2k+1}{2}}(s,x, D_x)|_o = 0;$

(iv) The harmonic oscillators are weightless and all of the ${\cal R}$'s have weight -2. }
\medskip

\noindent{\bf Proof}:  The operator $W_h$ will be constructed as the asymptotic product
$$W_h:=\mu(r_{\alpha})^*\Pi_{k=1}^{\infty}W_{h \frac{k}{2}}\mu(r_{\alpha}) \leqno(2.23)$$
of weightless unitary $h$-pseudodifferential operators on $\R^n$, with
$$W_{h \frac{k}{2}}:= exp(ih^{\frac{k}{2}} Q_{\frac{k}{2}})\leqno(2.23a)$$ 
and with $h^{\frac{k}{2}} Q_{\frac{k}{2}} \in h^{\frac{k}{2}}\C^{\infty}(S^1_L)\otimes {\cal E}^{k+2}$
of total order 1.
The product will converge, for each s, to a unitary operator in $\Psi_h^o( \R^n)$
(see [Sj] for discussion of asymptotic products). 

To see what is involved, 
  we first construct a weightless $Q_{\frac{1}{2}}(s, x, D_x)
 \in C^{\infty}(S^1_L)\otimes{\cal E}^3_{\epsilon}$ such that
$$e^{-ih^{\frac{1}{2}}Q_{\frac{1}{2}}} {\cal R}_h e^{ih^{\frac{1}{2}}Q_{\frac{1}{2}}}|_o
=[- h^{-2}L^{-2} + 2 h^{-1}L^{-1}{\cal R} + {\cal R}^{\half}_o + \dots]|_o\leqno(2.24a)$$
where the dots $\dots$ indicate higher powers in $h$.  
The
 operator $Q_{\frac{1}{2}}$ then must satisfy
the commutation relation
$$\{ [L^{-1}{\cal R},Q_{\frac{1}{2}}]+ {\cal R}_{\frac{1}{2}} \}|_o= 0.\leqno(2.24b)$$ 
 To solve for $Q_{\half}$,
it is convenient to conjugate back to the ${\cal D}_{2 - \frac{m}{2}}$'s of the
twisted model
 by $\mu(r_{\alpha})$, since this transforms ${\cal R}$
into $D_s$.  The commutation relation thus becomes
$$\{ [L^{-1} D_s,\mu(r_{\alpha})^*Q_{\frac{1}{2}}\mu(r_{\alpha})]+ 
{\cal D}_{\frac{1}{2}} \}|_o= 0, \leqno(2.24c)$$
that is,
$$L^{-1}\partial_s\{\mu(r_{\alpha})^*Q_{\frac{1}{2}}\mu(r_{\alpha})\}|_o =
 - i\{{\cal D}_{\frac{1}{2}}\}|_o\leqno(2.24d)$$
where $\partial_s A$ is the Weyl operator whose complete symbol is the $s$-derivative
of that of $A$.  Since (2.24d) is simpler than (2.24b), we henceforth conjugate everything
by $\mu(r_{\alpha})$, and relabel the operators $\mu(r_{\alpha})^*Q\mu(r_{\alpha})$ by
$\tilde{Q}.$  The resulting  ${\cal D}$'s  then have the twisted model
periodicity properties (2.17a-b).  Our problem is then to
 solve (2.24d) with an operator $\tilde{Q}_{\half}$ satisfying (2.17a-b).  

 To do so, we first observe that under conjugation by $\mu(r_{\alpha})$,
 elements of ${\cal R}$-weight zero transform to
elements of $D_s$- weight zero, and hence it
suffices to solve for the matrix elements $\langle \tilde{Q}_{\frac{1}{2}} \gamma_q, \gamma_r
\rangle$.  
  It follows from (2.17b) that the solution is unique for $r \not=q$, while for
$r=q$ the matrix element $\langle \tilde{Q}_{\half}\gamma_q,\gamma_q
\rangle$ is periodic and hence a necessary and sufficient condition for solvability
by an operator with the correct periodicity  is
that, for all $q \in {\bf N}^n,$
$$\int_o^L \langle {\cal D}_{\half} \gamma_q, \gamma_q \rangle ds = 0 
  \;\; \leqno(2.25).$$
In view of (2.16a-b), we  have
$${\cal D}_{\frac{1}{2}} \in \Psi^1(\R)_{\alpha}\otimes{\cal E}^{1}_{\epsilon} +
\Psi^o(\R)_{\alpha}\otimes{\cal E}^{3}_{\epsilon},$$ in fact a simple 
calculation shows that ${\cal D}_{\frac{1}{2}} \in 
\Psi^o(\R)_{\alpha}\otimes{\cal E}^{3}_{\epsilon}$.  (Recall  that the
subscript $\alpha$ indicates the periodicity property (2.17a-b).)

We now observe that if $A \in \Psi^*(\R) \otimes {\cal E}_{\epsilon}^m$ then
$$A:\left \{ \begin{array}{ll}
 L^2_{\pm} \rightarrow L^2_{\pm} & 
m\;\; \mbox{even}\\ : L^2_{\pm} \rightarrow
 L^2_{\mp} & m \;\;\mbox{odd} \end{array} \right \},\leqno(2.26)$$
where $L^2_{+}$ (resp. $L^2_{-}$) denotes the subspace of
$L^2(\R)\otimes L^2(\R^n)$ spanned by even (resp. odd) functions
of $y \in \R^n$ with coefficients in $s$.

It follows that
$$ {\cal D}_{\frac{1}{2}} : L^2_{+} \rightarrow L^2_{-},$$
This parity reversing property implies the vanishing of the diagonal matrix elements in the basis
$\{\gamma_q\},$ hence (2.25) does hold.  A unique solution of (2.24d) is specified by the
condition that 
$$\int_o^{2 \pi}\langle \tilde{Q}_{\half} \gamma_q, \gamma_q \rangle ds = 0.$$

To solve the equation (2.24d) we rewrite it in terms of  complete Weyl symbols.
We will use the notation $A(s,x,\xi)$ for the complete Weyl symbol of the
operator $A(s,x,D_x)$.  Then (2.24d) becomes
$$L^{-1}\partial_s \tilde{Q}_{\half}(s,x,\xi)= -i {\cal D}_{\half}|_o(s,x,\xi)
\leqno(2.27a)$$
with
$$\tilde{Q}_{\half}(s + L,x,\xi) = \tilde{Q}_{\half}(s,r_{\alpha}(L)(x,\xi)).$$
We solve (2.27a) with the Weyl symbol
$$\tilde{Q}_{\half}(s,x,\xi) = \tilde{Q}_{\half}(0,x,\xi) + L \int_0^s
-i {\cal D}_{\half}|_o(u,x,\xi)du$$
where $\tilde{Q}_{\half}(0,x,\xi)$ is determined by the consistency condition
$$\tilde{Q}_{\half}(L,x,\xi) - \tilde{Q}_{\half}(0,x,\xi) = 
L \int_0^L -i {\cal D}_{\half}|_o(u,x,\xi)du \leqno(2.27b)$$
or in view of the periodicity condition in (2.27a),
$$\tilde{Q}_{\half}(0,r_{\alpha}(x,\xi)) - \tilde{Q}_{\half}(0,x,\xi) = 
L \int_0^L -i {\cal D}_{\half}|_o(u,x,\xi) du.\leqno(2.27c)$$
To solve, we use that ${\cal D}_{\half}|_o(u,x,\xi)$ is a polynomial of degree 3 in
$(x,\xi)$. Also, as is customary in such calculations, we switch to complex
coordinates $z_j = x_j + i\xi_j$ and $\bar z_j = x_j - i \xi_j$ in which the
action of $r_{\alpha}(L)$ is diagonal.
 With a
little abuse of notation, we will continue to denote the Weyl symbols, qua functions
of the $z_j,\bar z_j$'s by their previous expressions. We also suppress the
subscripts by using vector notation $z, \bar z$ and $e^{i\alpha}$.
 Thus, (2.27c) becomes
$$\tilde{Q}_{\half}(0,e^{i\alpha}z, e^{-i\alpha}\bar z) -
\tilde{Q}_{\half}(0,z,\bar z) = 
L\int_0^L -i {\cal D}_{\half}|_o(u,z,\bar z ) du.\leqno(2.28)$$
We now use that ${\cal D}_{\half}(u,z,\bar z )$ is a polynomial of degree 3
to solve (2.27c).  If we put
$$\tilde{Q}_{\half}(s,z,\bar z) = \sum_{|m|+|n|\leq 3} q_{\half;mn}(s) z^m \bar z^n$$
and 
$$ {\cal D}_{\half}|_o(s,z,\bar z ) du = \sum_{|m|+|n|\leq 3} d_{\half;mn}(s)
z^m \bar z^n$$
then (2.28) becomes
$$\sum_{|m|+|n|\leq 3} (1 - e^{(m - n) \alpha}) q_{\half; mn}(0) z^m \bar z^n
= -i L^2 \sum_{|m|+|n|\leq 3} \bar d_{\half;mn}
z^m \bar z^n \leqno (2.29).$$
Since there are no terms with $m=n$ in this (odd-index) equation, and since
the $\alpha_j$'s are independent of $\pi$ over $\Z$, there is no obstruction to
the solution of (2.29).

  For simplicity of notation we will express the solution
in the form
$$\tilde{Q}_{\frac{1}{2}}|_o =-i L \int {\cal D}_{\frac{1}{2}}(s)ds|_o\leqno(2.30)$$
where $\int$  denotes the integration procedure just defined, that is, the indefinite integral
satisfying (2.27b). 
 Since the
integration only involves the $s$-coefficients, the solution is a pseudodifferential
operator on
$\R^n$ with the same order,  same order of vanishing, and same parity 
as the restriction of ${\cal D}_
{\frac{1}{2}}$ to elements of weight zero.  We then extend it 
all of to ${\cal H}_{\alpha}$
as a pseudodifferential operator with the same properties by stipulating that
$$\tilde{Q}_{\frac{1}{2}}M_h = M_h \tilde{Q}_{\frac{1}{2}}. \leqno(2.31a)$$   Then, as
desired, 
$$\tilde{Q}_{\half} \in \Psi^o(\R^1) \otimes {\cal E}^3_{\epsilon}.\leqno(2.31b)$$
The conjugate by $\mu(r_{\alpha})$ then defines a periodic operator satisfying 
(2.24b)
and hence a unitary element $W_{h \half} \in \Psi^o_h(S^1 \times \R^n)$
 satisfying (2.24a).
The twisted unitary operator with exponent $\tilde{Q}_{\half}$, i.e.the
image of $W_{h \half}$ under
conjugation by
$\mu(r_{\alpha})$, will be   decorated with a tilde, $\tilde{W}_{h \half}.$ 

 Since $h^{\half}\tilde{ Q}_{\half}$ is of total order 1,
     $h^{\half} ad(\tilde{ Q}_{\half})$ (with $ad(A)B:=[B,A]$)
 preserves the
total order in $\Psi_h^{(*,*,*)}$, and hence $\tilde{W}_{h \half}$ 
is an order-preserving automorphism.
 It is moreover independent of $D_s$ and has an odd polynomial Weyl symbol, so that
$$h^{\half} ad(\tilde{Q}_{\half}): h^{\frac{k}{2}}\Psi^l(\R) \otimes {\cal E}^m_{\epsilon}
\rightarrow h^{\frac{k+1}{2}} [\Psi^{l-1}(\R) \otimes {\cal E}^{m+3}_{\epsilon} + 
\Psi^{l}(\R)
\otimes {\cal E}^{m+1}_{\epsilon}]. \leqno(2.32)$$
Finally, since ${\cal D}_{\half}$ has weight -2, $L$ has degree 1 and since
 $\int_o^s (\cdot)ds$ has degree 1 in the metric scaling, we see that $\tilde{Q}_{\half}$ is weightless.
Alternatively, the $d_{\half;m,n}$'s have weight -2,the variables $z$ have weight 0 and hence the
$q_{\half;m,n}$'s have weight 0.

 Consider now the element
$${\cal D}^{\half}_h:=\tilde{W}_{h \frac{1}{2}}^* {\cal D}_h \tilde{W}_{h \frac{1}{2}} \in \Psi^2_h(\R^1\times
\R^n).$$
Using only the $h$-filtration, we  expand
$${\cal D}^{\half}_h\sim \sum_{n=o}^{\infty} h^{-2 + \frac{n}{2}} \sum_{j+m=n} \frac{i^j}{j!}
(ad\tilde{Q}_{\half})^j {\cal D}_{2 - \frac{m}{2}}\leqno(2.33)$$
$$:= h^{-2}L^{-2} + h^{-1}L^{-1}D_s +\sum_{n=3}^{\infty}
 h^{-2 + \frac{n}{2}}{\cal D}^{\half}_{2 - \frac{n}{2}}.$$ 
An obvious induction using (2.32) and (2.16a-b) gives that
$$ad(\tilde{Q}_{\half})^j {\cal D}_{2 - \frac{m}{2}} \in C^{\infty}(\R, {\cal E}_{\epsilon}^{m+j-4})D_s^2
+ C^{\infty}(\R, {\cal E}_{\epsilon}^{m+j-2})D_s + C^{\infty}(\R, {\cal E}_{\epsilon}^{m+j}).$$
It follows that
${\cal D}^{\half}_{2 - \frac{n}{2}}$ has the same filtered structure (2.16b) as 
${\cal D}_{2 - \frac{n}{2}}.$

We carry this procedure out one more step before arguing inductively, since the even
steps behave differently from the odd ones.  
We thus seek an element $\tilde{Q}_1(s,x,D_x)\in \Psi^*(S^1 \times \R^n)$ and an element
$f_o(I_1,...,I_n) \in {\cal A}$ so that
$${\cal D}^1_h:=\tilde{W}_{h 1}^*{\cal D}^{\half} \tilde{W}_{h 1} = h^{-2}L^{-2} + h^{-1}L^{-1}D_s +
h^{-\half}{\cal D}^{\half}_{\half} + {\cal D}^1_o(s,D_s, x,D_x) + \dots$$
with 
$${\cal D}_o^1(s,D_s,x, D_x)|_o = f_o(I_1,...,I_n)\leqno(2.34a)$$
with $ \tilde{W}_{h 1}= e^{i h \tilde{Q}_1},$ and where the dots signify
terms of higher order in $h$.  Note that ${\cal D}^1_{\half} = {\cal D}^{\half}_{\half}$, so that
(2.33a) would imply
$$\{h^{\half}{\cal D}^1_{\half} +{\cal D}_o^1\}|_o = f_o(I_1,...,I_n).
\leqno(2.34b)$$
 The condition on $\tilde{Q}_1$ is then
$$\{[D_s,\tilde{Q}_1] + {\cal D}_o^{\half}\}|_o = f_o(I_1,...,I_n) \leqno(2.35a)$$
 or equivalently
$$\partial_s\tilde{Q}_1|_o =\{-{\cal D}_o^{\half} + f_o(I_1,...,I_n)\}|_o
\leqno(2.35b).$$ 

As above, we first consider the solvability of (2.35b) one matrix element at a time. 
The off-diagonal matrix element equations again have a unique solution, but because
${\cal D}_o^{\half}$ is even there is now a condition on the solvability, with a
periodic solution,  of the diagonal ones: 
$$\frac{1}{L}\int_o^L \langle {\cal D}_o^{\half} \gamma_q, \gamma_q \rangle ds =
  f_o(q_1 + \half,...,q_n + \half). \leqno(2.36)$$
These conditions (2.36) determine $f_o$ and hence the operator $f_o(I_1,...,I_n).$  To
analyse the properties of $f_o$ and of $\tilde{Q}_1$, we use (1.1b), (2.16b) and
(2.31c) to conclude  that
$${\cal D}_o^{\half} = {\cal D}_o + ad\tilde{Q}_{\half} (ad\tilde{Q}
_{\half}D_s + {\cal D}_{\half})$$
$$\in C^{\infty}(\R,{\cal E}^o_{\epsilon})D_s^2
 + C^{\infty}(\R, {\cal E}^{2}_{\epsilon})D_s + C^{\infty}(\R, {\cal E}^4_{\epsilon}).$$
 We also observe that (2.32) is equivalent to
$$f_o(I_1,...,I_n) = \frac{1}{L} \int_o^L \int_{T^n} V_t^* {\cal D}^{\half}_o|_o V_t 
dt ds \leqno(2.37)$$
where $T^n = \R^n / {\Z}^n$ is the n-torus, where $t \in T^n$ and where
$$V(t_1,...,t_n):=expit_1 I_1 \cdot  \dots \cdot expit_nI_n.$$
Since $V_t$ belongs to the metaplectic representation, we have by metaplectic 
covariance of the Weyl calculus that
$$f_o(I_1,...,I_n) \in {\cal P}_{\cal I}^2 \subset{\cal E}^4 \leqno (2.38).$$
  We then can
solve for $\tilde{Q}_1$ in the form
$$\tilde{Q}_1(s,y,D_y):= -i L [\int \{{\cal D}^{\half}_o - f_o(I_1,...,I_n)\}|_ods] \in
 \Psi^o(\R^1) \otimes {\cal E}^4_{\epsilon} 
\leqno(2.39)$$
where $\int$ is the indefinite integration in $s$ with zero mean value and where the
operator on the right is 
 extended to $L^2(\R^1 \times \R^n)$ using (2.35b).

To make the integration more concrete and to analyze the solution, we again
express everything in terms of complete Weyl symbols.  Thus we rewrite (2.35b)
in the form
$$L^{-1} \partial_s \tilde{Q}_1(s,z, \bar z) = -i \{{\cal D}^{\half}_o|_o (s,z,\bar z) -
f_o(|z_1|^2,\dots, |z_n|^2)\} \leqno (2.40a)$$
or equivalently
$$\tilde{Q}_1(s,z,\bar z) = \tilde{Q}_1(0,z,\bar z) -i L\int_0^s
[{\cal D}^{\half}_o|_o (u,z,\bar z) -
f_o(|z_1|^2,\dots, |z_n|^2)] du\leqno (2.40b)$$  
and solve simeltaneously for $\tilde{Q}_1$ and $f_o$. The consistency condition
determining a unique solution is that
$$\tilde{Q}_1(L,z,\bar z) = \tilde{Q}_1(0,z,\bar z) -i  L\int_0^L
[{\cal D}^{\half}_o|_o (u,z,\bar z) -
f_o(|z_1|^2,\dots, |z_n|^2)] du. \leqno(2.41a)$$
or in view of the twisted periodicity condition
$$\tilde{Q}_1(0,e^{i\alpha}z,e^{-i\alpha}\bar z) - \tilde{Q}_1(0,z,\bar z) = -i L \{\int_0^L
{\cal D}^{\half}_o|_o (u,z,\bar z)du  -
L f_o(|z_1|^2,\dots, |z_n|^2) \}.\leqno (2.41b)$$
As before we use that ${\cal D}^{\half}_o|_o (u,z,\bar z)$ is a polynomial of degree
4 to solve the equation.  We put
$$ \tilde{Q}_1(s,z,\bar z) = \sum_{|m|+|n|\leq 4} q_{1;mn}(s) z^m \bar z^n, \;\;\;\;\;\;
 f_o(|z_1|^2,\dots, |z_n|^2) = \sum_{|k|\leq 2} c_{o k} |z|^{2k} \leqno(2.42)$$
and  
$${\cal D}^{\half}_o|_o (s,z,\bar z)du :=\sum_{|m|+|n|\leq 4} d_{o; mn}^{\half}(s) z^m \bar
z^n,\;\;\;\;\;\;\;\;\;\bar d^{\half}_{o;mn}:= \frac{1}{L}\int_o^Ld^{\half}_{o;mn}(s)ds $$
As above, we can solve for the off-diagonal coefficients,
$$q_{1; mn}(0) = -i L^2 (1 - e^{i (m-n)\alpha})^{-1} \bar d_{o; mn}^{\half} \leqno(2.43a)$$
and must set the diagonal coefficients equal to zero.  The coefficients $c_{o k}$
are then determined by
$$c_{o k} =   \bar d_{1; kk}^{\half}. \leqno(2.43b)$$

It is evident that $\tilde{Q}_1$ and $f_o(I_1,\dots,I_n)$ are
even polynomial pseudodifferential  operators of degree 4 in the variables $(x,D_x)$, that $\tilde{Q}_1$
is weightless under metric rescalings and that the coefficients $c_{o k}$ are of weight -2. They are
essentially the same as the residual QBNF invariants.

We now proceed inductively to define self-adjoint polynomial pseudodifferential operators
$\tilde{Q}_{\frac{j}{2}}, f_j(I_1,...,I_n)$,
 unitary $h$-pseudodifferential operators
$\tilde{W}_{h \frac{j}{2}}:=exp(i h^{\frac{j}{2}} \tilde{Q}_{\frac{j}{2}})$ 
and approximate semi-classical
normal forms ${\cal D}_h^{\frac{k}{2}}$
such that:

\noindent(2.44)$$\begin{array}{l}(i)\;\;{\cal D}_h^{\frac{k}{2}}:=
\tilde{W}_{h \frac{k}{2}}^*\tilde{W}_{h \frac{k-1}{2}}^*
 \dots \tilde{W}_{h \frac{1}{2}}^* {\cal D}_h
\tilde{W}_{h \half} \dots \tilde{W}_{h \frac{k-1}{2}}\tilde{W}_{h \frac{k}{2}};\\

(ii)\;\;{\cal D}_h^{\frac{k}{2}} \sim h^{-2} + h^{-1}D_s +
\sum_{n=3}^{\infty} h^{-2 + \frac{n}{2}}{\cal D}^{\frac{k}{2}}_{2 - \frac{n}{2}};\\

(iii)\mbox{ for}\; k \geq n-2,
{\cal D}^{\frac{k}{2}}_{2 - \frac{n}{2}} = {\cal D}^{\frac{n-2}{2}}_{2 -\frac{n}{2}}
=[D_s, \tilde{Q}_{\frac{n}{2}-1}] + {\cal D}^{\frac{n-3}{2}}_{-\frac{n}{2}}\\

(iv)\mbox{for}\;\;k\geq n-2,
{\cal D}^{\frac{k}{2}}_{2-\frac{n}{2}}|_o \left \{
\begin{array}{ll}= f_j(I_1,...,I_n) & n=2j\\ =0 & n=2j+1 \end{array} \right\} \\

 (v)\;\;{\cal D}^{\frac{k}{2}}_{2-\frac{m}{2}} \in C^{\infty}(\R,
{\cal E}^{m-4}_{\epsilon})D_s^2
 +C^{\infty}(\R, {\cal E}^{m-2}_{\epsilon})D_s
+C^{\infty}(\R, {\cal E}^{m}_{\epsilon});\\

(vi)\;\;\tilde{Q}_{\frac{k}{2}} \in C^{\infty}(\R,{\cal E}^{k+2}_{\epsilon});\\

(vii) \;\;\mbox{the conjugates under} \;\;\mu(r_{\alpha})\;\;\mbox{of the above operators
are periodic of period L in s};\\

(viii)\;\;f_j(I_1,...,I_n) \in {\cal P}^{j+2}_{\cal I};

(ix) \;\; wgt(f_j) = -2. \end{array}$$

To check the induction, let us assume these properties hold for $k \leq N$. Then

$${\cal D}_h^{\frac{N+1}{2}}=\tilde{W}_{h \frac{N+1}{2}}^*{\cal D}^{\frac{N}{2}}_h
\tilde{W}_{h \frac{N+1}{2}}= \leqno(2.45)$$
$$\begin{array}{l}=h^{-2}L^{-2} + h^{-1}L^{-1}D_s + \sum_{n=3}^{N+2}h^{-2+\frac{n}{2}}{\cal
D}^{\frac{N}{2}}_{2-\frac{n}{2}}
 +h^{\frac{N-1}{2}}\{ [L^{-1} D_s, \tilde{Q}_{\frac{N+1}{2}}]+ {\cal D}^{\frac{N}{2}}_{-\frac{N-1}{2}}\}
+ \sum_{n=N+4}^{\infty} h^{-2 + \frac{n}{2}}{\cal D}^{\frac{N}{2}}_{2-\frac{n}{2}}\\
+\sum_{m=2}^{\infty}h^{m\frac{N+1}{2}}\frac{i^m}{m!}(ad \tilde{Q}_{\frac{N+1}{2}})^m (D_s)
+\sum_{m=1}^{\infty}\sum_{p=3}^{\infty} h^{m\frac{N+1}{2}+\frac{p}{2}-2}\frac{i^m}{m!}
(ad \tilde{Q}_{\frac{N+1}{2}})^m 
({\cal D}^{\frac{N}{2}}_{2 -\frac{p}{2}}).\end{array}$$
We  note that
$${\cal D}^{\frac{N+1}{2}}_{2-\frac{N+3}{2}} =
 [L^{-1} D_s, \tilde{Q}_{\frac{N+1}{2}}] + {\cal D}^{\frac{N}{2}}
_{- \frac{N-1}{2}} $$
and that further conjugations will not alter
this (or lower order) term(s), proving (iii).

When $N+1$ is odd, $\tilde{Q}_{\frac{N+1}{2}}$ must solve:
$$\{ [L^{-1} D_s, \tilde{Q}_{\frac{N+1}{2}}]+ {\cal D}^{\frac{N}{2}}_{-\frac{N-1}{2}}\}|_o = 0
 \leqno(2.46a)$$
while if $N+1$ is even $\tilde{Q}_{\frac{N+1}{2}}$ and $f_{\frac{N-1}{2}}$ must solve
$$\{ [L^{-1} D_s, \tilde{Q}_{\frac{N+1}{2}}]+ {\cal D}^{\frac{N}{2}}_{-\frac{N-1}{2}}\}|_o = f_{\frac{N-1}{2}}
(I_1,...,I_n). \leqno(2.46b)$$
 
If $N+1$ is odd, (2.45v) implies that the diagonal matrix elements of
${\cal D}^{\frac{N}{2}}_{-\frac{N-1}{2}}|_o$ vanish.  Hence, as in the case of $\tilde{Q}_{\frac{1}{2}}$,
the solution of (2.46a) is given by
$$\tilde{Q}_{\frac{N+1}{2}}|_o:=
 - i L \int{\cal D}^{\frac{N}{2}}_{-\frac{N-1}{2}}|_o ds \leqno(2.47a)$$
where the constant of integration $\tilde{Q}_{\frac{N+1}{2}}|_o(0,z,\bar z)$ is defined so that the solution
satisfies the periodicity condition analogous to (2.28).
In the even case, as in the case of $\tilde{Q}_1$ we set
$$f_{\frac{N-1}{2}}(q_1 + \half,...,q_n + \half):=\frac{1}{L}\int_o^L 
\langle  {\cal D}^{\frac{N}{2}}_{-\frac{N-1}{2}} \gamma_q, \gamma_q \rangle ds,
 \leqno(2.47b)$$
or equivalently
$$f_{\frac{N-1}{2}}(I_1,...,I_n):=\frac{1}{L}\int_o^L\int_{T^n}
V_t^*{\cal D}^{\frac{N}{2}}_{-\frac{N-1}{2}}|_o V_tdtds, \leqno(2.47c)$$
and
$$\tilde{Q}_{\frac{N+1}{2}}|_o:=
 -L\int\{{\cal D}^{\frac{N}{2}}_{-\frac{N-1}{2}} - f_{\frac{N-1}{2}}
(I_1,...,I_n)\}|_ods. \leqno(2.47d)$$
As above, $\tilde{Q}_{\frac{N+1}{2}}|_o$ is 
then extended to all of ${\cal H}_{\alpha}$ thru
(2.31b). The indefinite integrations can be precisely defined, as above, by
expressing everything as a polynomial in $z,\bar z$ and solving the resulting
algebraic equations for the constant of integration.
 By (2.47 c-d) the solution has the parity of ${\cal
D}^{\frac{N}{2}}_{-\frac{N-1}{2}}$, i.e. the parity of $N-1,$ and by   (2.45)  the
parity property  propagates to the case $N+1$.  The formula (2.47c) also shows that
for even $N-1$,
$f_{\frac{N-1}{2}}(I_1,...,I_n)$ has the same order and weight  properties as 
${\cal D}^{\frac{N}{2}}_{-\frac{N-1}{2}}$.  Since the latter may be written in
the form $aD_s^2 + bD_s + c,$ with $c \in C^{\infty}(\R, {\cal E}^{N+3}_{\epsilon})$,  the former
lies in ${\cal P}_{\cal I}^* \cap {\cal E}^{N+3}_{\epsilon}$, implying $(viii)$.  The formula
(2.47d) then implies $(vi)$.  The periodicity property (vii) is maintained
throughout. 

It follows immediately that
$$\begin{array}{ll} {\cal D}^{\infty}_{2 - \frac{n}{2}} &=
{\cal D}^{\frac{n-2}{2}}_{2 -\frac{n}{2}}\\
&= [D_s,\tilde{ Q}_{\frac{n}{2}-1}] + 
 {\cal D}^{\frac{n-3}{2}}_{2-\frac{n}{2}}.  \end{array}\leqno(2.48)$$
Conjugating back under $\mu(r_{\alpha})$, we see that
statement $(i)$ of the Lemma then follows from (2.48) and from (2.44iii-iv) and 
statement (ii) then follows from (2.44 viii). The weight property (ix) is visible from (2.47c). \qed
\medskip

As a transitional step to the quantum Birkhoff normal form, let us rephrase the above Lemma in terms  of the
actual  Fermi normal coordinates and inverse Planck constants
 $\{ r_{kq}\}$ in place of  $(Lh)^{-1}$.  We thus
 introduce the space
${\cal O}^*(N_{\gamma}, \Gamma, \{r_{kq}\})$ of semi-classical
 Hermite distributions associated
to the non-homogenenous isotropic manifold $\Gamma:=\R^+ \gamma$, i.e.
$$\Gamma:=\{(s, \sigma, y,\eta) \in T^*(N_{\gamma}): \sigma = 1, (y, \eta) = (0,0) \}.$$
It is, by definition, the union of the spaces ${\cal O}^{m/2}$ 
of elements of order $\frac{m}{2}, m \in{\bf N}$,
 given in  Fermi normal coordinates by asymptotic sums of the form
 $$u(s,y, r_{kq}^{-1}) \sim (r_{kq})^{\frac{m}{2}} e^{i r_{kq}s} \sum_{n=o}^
{\infty} (r_{kq})^{-\frac{n}{2}} f_n(s, \sqrt{r_{kq}}
y)$$
with $f_n\in C^{\infty}_q(\R,{\cal S}( \R^n))$.
Here, ${\cal S}( \R^n)$ is the Schwarz space, and the subscript q denotes the space
of functions of this form satisfying
$$f(s+L, y) = e^{i\kappa_q s} f(s, ty).$$  Aside 
from  the restriction to Schwartz functions and the half-integral orders, it is the image
under the rescaling operator $T_{(r_{kq})^{-1}}$ of the space
${\cal O}^*(N_{\gamma}, \Lambda, \{r_{kq}\})$  of semi-classical
Lagrangean distributions associated to the non-homogenenous Lagrangean
$$\Lambda:= graph(ds)=\{(s, \sigma, y, \eta): \sigma=1, \eta=0.\}.$$ Similar spaces
of oscillatory functions could be, and implicitly are, defined in the model $S^1 \times \R^n$
but the scaling aspect comes out most naturally on $N_{\gamma}$.  (For background on
semi-classical Lagrangean distributions, see [CV2]).

Let us denote by ${\cal O}^*_o(N_{\gamma}, \Lambda,\{r_{kq}\})$
 the subspace of elements of ${\cal O}^*(S^1 \times \R^n, \Lambda)$  
satisfying ${\cal L} u(s,y, r_{kq}^{-1} ) \sim 0.$  Let us also denote by $W_{kq}$ the
transfer of $W_{r_{kq}^{-1}}$ above to the normal bundle, i.e.
$$W_{kq}:= \mu(\tilde{a}_s) W_{r_{kq}^{-1}} \mu (\tilde{a}_s)^*.$$ For the sake of simplicity we will
be a little negligent here of the rescalings.
Under this operator, the kernel of ${\cal L}$
  goes over to the space $\tilde{{\cal O}}
^*_o(N_{\gamma}, \Lambda, \{r_{kq}\})$ of elements annihilated by $W_{kq}
^* {\cal L} W_{kq}.$
  It is clear from the above lemma
that this space is stable under $W_{kq}* \Delta_{r_{kq}^{-1}} W_{kq}.$  We may interpret
what was proved in the lemma as giving a semi-classical normal form for $\Delta$ or
$\Delta_{r_{kq}^{-1}}$ in the following sense:

\medskip

\noindent(2.49) {\bf Theorem} \; {\it The $ r_{kq}^{-1}$-pseudodifferential operator
$$W_{kq}: {\cal O}^*_o(N_{\gamma},\Lambda, \{r_{kq}\}) \rightarrow \tilde
{{\cal O}}^*_o(N_{\gamma},\Lambda, \{r_{kq}\})$$
has the properties:

$$W_{kq}^* W_{kq}^{-1} - I\sim W_{kq}W_{kq}^*-I \sim 0 $$ 
$$W_{kq}^* T_{r_{kq}^{-1}}^{-1} \Delta T_{r_{kq}^{-1}} W_{kq} \sim
 {\cal L}^2 + f_o(I_{\gamma1},...,I_{\gamma n}) 
+ \frac{f_1(I_{\gamma1},...,I_{\gamma n})}{{\cal L}} + \dots.$$
$$ W_{kq}  e^{ir_{kq}s}U_q(s, y) = e^{ir_{kq}s}W_{kq}U_q(s, y)
\Phi_{kq}(s,y)$$
(see (2.1)).}

\section{Normal form: Proof of Theorem B}

The intertwining operators $W_{kq}$ will now be assembled into the
Fourier-Hermite -series integral operator
$$W_{\gamma} : L^2(S^1_L \times \R^n, dsdy) \rightarrow L^2(S^1_L \times \R^n, dsdy)\leqno(3.1)$$
$$W_{\gamma} \sum_{(k,q) \in {\bf N}^{n+1}} \hat{f}(k,q) e^{ir_{kq}s}U_q(s,y) =
\sum_{(k,q) \in {\bf N}^{n+1}} \hat{f}(k,q) e^{ir_{kq}s} W_{kq}U_q(s,y).$$
 Also, the
dilation operators will be assembled into the  operator
$$T : L^2(S^1_L \times \R^n, dsdy) \rightarrow L^2(S^1_L \times \R^n, dsdy)\leqno(3.2)$$ 
$$T \sum_{(k,q) \in {\bf N}^{n+1}} \hat{f}(k,q) e^{ir_{kq}s}U_q(s,y) =
\sum_{(k,q) \in {\bf N}^{n+1}} \hat{f}(k,q) e^{ir_{kq}s}U_q(s,\sqrt{r_{kq}}y).$$ 
By theorem (2.49) we then have, at least formally,
$$W_{\gamma}^{-1} T^{-1} \Delta T W_{\gamma} \sim
 {\cal L}^2 + f_o(I_{\gamma1},...,I_{\gamma n}) 
+ \frac{f_1(I_{\gamma1},...,I_{\gamma n})}{{\cal L}} + \dots.\leqno(3.3)$$
The purpose of this section is to make  this equivalence precise.  Since it is independent of the
model we will carry out the proof in the basic model and use the weightless Fermi normal coordinates
$(s,x)$.  For the sake of notational simplicity the $\tilde{W}'s$ and $\tilde{Q}$'s, transferred back
to the model in this way, will be written as $W$'s and $Q$'s. Also in place of the $r_{kq}^{-1}$'s we will
use the weightless $h_{kq}$'s with $h_{kq}^{-1} = (2\pi k + \sum_{j=1}^n(q_j + \half)\alpha_j) = L r_{kq}.$ 
\medskip

\noindent(3.4) \; {\bf Propostion} {\it $\; T W_{\gamma}T^{-1}$ is a
 (standard) Fourier integral operator, well-defined and invertible 
on the microlocal neighborhood (0.1) in $T^*(S^1_L \times \R^n)$.}
\medskip

\noindent{\bf Sketch of Proof}:  

To simplify, we first consider the unitarily equivalent operator
$\tilde{T}\tilde{W}\tilde{T}^{-1}$ in the microlocal neighborhood (0.1) in the
twisted model, with
$$\tilde{W}: {\cal H}_{\alpha} \rightarrow {\cal H}_{\alpha} \leqno(3.5)$$
$$\tilde{W}(e^{i\frac{s}{Lh_{kq}}}\gamma_q):=e^{i\frac{s}{Lh_{kq}}}W_{h_{kq}}\gamma_q,$$ and
with $\tilde{T}$ the dilation operator like (3.2) relative to the basis
$e^{i\frac{s}{Lh_{kq}}}\gamma_q$.  
We then factor $\tilde{T}\tilde{W}\tilde{T}^{-1}$ as the product 
 $\tilde{T}\tilde{W}\tilde{T}^{-1}
= j^*\tilde{T} V \tilde{T}^{-1} $ where:
 $$ V: {\cal H}_{\alpha} \rightarrow L^2_{loc}(\R \times \R \times
\R^n)\leqno(3.6)$$
$$V=\Pi_{j=o}^{\infty} exp[i D_s^{-\frac{j}{2}} Q_{\frac{j}{2}}(s',y,D_y)]$$ that is,
$$Ve^{i\frac{s}{Lh_{kq}}}\gamma_q(x):=e^{i\frac{s}{Lh_{kq}}}W_{h_{kq}}(s',x,D_x)\gamma_q(x),$$
and where
 $$j^*: C^{\infty}(\R \times \R \times \R^n) \rightarrow C^{\infty}(\R \times \R^n)
\leqno(3.7)$$
$$j^*f(s,x)=f(s,s,x)$$
is the pullback under the partial diagonal embedding.  We note that the variable $s'$ in
$V$ occurs essentially as a parameter, so we can (and often will) regard V as a one parameter
family of operators $V_{s'}$ on ${\cal H}_{\alpha}.$

 It is not clear that the infinite
product in (3.6) is well-defined, nor what kind of operators the factors are.  On the 
second point, we note that $D_s$ is elliptic in the set (0.1) in $T^*(S^1_L \times
\R^n)$ and that, as an operator-valued function of $s'$,
$$\tilde{T}D_s^{-\frac{j}{2}} Q_{\frac{j}{2}}(s',x,D_x)
 \tilde{T}^{-1} \in C^{\infty}(S_L^1,\Psi^1(S_L^1
 \times \R^n)). \leqno(3.8)$$
Here of course $C^{\infty}(S^1_L, {\cal A})$ denotes the smooth functions
 of $s'$ with values
in the algebra ${\cal A}$; to simplify the notation we do not indicate the
possible twisting.
Indeed, we first observe that conjugation by $\tilde{T}$ converts the mixed polyhomogeneous-
isotropic algebra $\Psi^*(S^1_L)\otimes {\cal W}^*$ into the pure polyhomogeneous
algebra $\Psi^*(S^1_L)\otimes\Psi^*(\R^n)$.  This essentially follows from the fact that 
$$\tilde{T} a^w(x,D_x) \tilde{T}^{-1} = a^w(D_s^{\half} x, D_s^{-\half} D_x)$$  
 (see [G.2] or [CV] for further details).  Obviously, 
$\tilde{T}D_s\tilde{T}^{-1} \in \Psi^1(S^1_L \times \R^n)$,
and by (2.44(vi)) we  also have $$\tilde{T} Q_{\frac{j}{2}} \tilde{T}^* \in C^{\infty}
(S^1_L,\Psi^{j+1}(S^1\times \R^n)).$$   
  A simple calculation shows further that  $\tilde{T}D_s^{-\frac{j}{2}} 
Q_{\frac{j}{2}}(s',x,D_x)\tilde{T}^{-1}$ is  a
one-parameter family of pseudodifferential
operators of real principal type over $ S^1_L \times \R^n$
 whose principal symbols have nowhere radial Hamilton
vector field.  It follows in a well-known way that these operators
are microlocally conjugate to $D_1$ (see [Ho IV, Proposition 26.1.2]) and hence that
their exponentials give a smooth one-parameter family  of Fourier Integral operators. 
 Consequently, after conjugation with $\tilde{T}$, the factors in (3.6) are for each $s'\in S^1$
 microlocally well-defined Fourier Integral operators in the
neighborhood (0.1), with smooth dependence on $s'$.
As for the infinite product, we note that the principal symbol of   $\tilde{T}D_s^{-\frac{j}{2}}
 Q_{\frac{j}{2}}\tilde{T}^{-1}$ has the form $\sigma^{-\frac{j}{2}} q_{\frac{j}{2}}(s',
\sqrt{\sigma} y,\frac{1}{\sqrt{\sigma}}\eta)$, where $q_{\frac{j}{2}}$ is a homogenous
polynomial in $(y,\eta)$ of degree j+2 (2.36(iv)).  Its exponential can therefore be
constructed microlocally in (0.1) as 
a one parameter family of Fourier integral operators over $S^1_L\times \R^n$
 with phase functions of the form
$$\Psi_{\frac{j}{2}}(s',s,s''\sigma, x,x'', \xi) =  (s-s'')\sigma + 
\langle x-x'', \xi \rangle + \psi_{\frac{j}{2}}(s',s,s'',\sigma, x,x'', \xi)$$
with
$$\psi_{\frac{j}{2}} \in O_{j+2}S^1_{cl}
 \leqno(3.9b)$$
Here, $S^1_{cl}$ denotes the space of classical symbols of first order in (0.1) and $O_k$ denotes the
elements which vanish to order k along $(x,\xi)=(0,0).$  (Henceforth we will use the
symbol $O_k$ for objects of any kind which vanish to order k on this set.)
 It follows that for fixed $s'$ and sufficiently small $\epsilon$, 
the phase parametrizes the graph
of a homogeneous canonical transformation of $T^*(S^1_L \times \R^n)$
 which is $C^{j+1}$- close to the identity.
Any finite product $\Pi_{j=o}^N$ in (3.7) is therefore a clean composition and the phase 
$\Psi_N$ of the product
parametrizes the corresponding compostion $\chi_N$ of canonical transformations.  By
(3.9b), we have $\chi_{N+1} - \chi_N \in O_{N+1}$ and so there exists a smooth
one parameter family homogeneous
canonical transformations $\chi_{s'}$, equal to the identity along $(x,\xi)=(0,0)$ such
that $\chi_{s'} - \chi_{s',N} \in O_{N+1}$    and with a generating function $\psi_{\infty}$ 
satisfying $\psi_{\infty} - \psi_N \in O_N$ (cf. [Sj, Proposition 1.3]).  Regarding the amplitudes,
the situation is of a similar kind.  Denoting the amplitude of $\tilde{T}exp(D_s^{-\frac{j}{2}}
 Q_{\frac{j}{2}})\tilde{T}^{-1}$ by $a_{\frac{j}{2}}(s',s,\sigma, y \eta)$ we have
$$a_{\frac{j}{2}} \equiv 1\;\;\; mod O_{j+2}S^o. \leqno(3.10)$$
It follows that the amplitudes $a_N$ of the products $\Pi_{j=o}^N$ satisfy
$a_{N+1} - a_{N} \in O_{N+1}S^o$ and hence there exists an element $a_{\infty}$ satisfying
$a_{\infty} - a_N \in O_{N+1}S^o$.  The pair $(a_{\infty}, \psi_{\infty})$ then determine
a one-parameter family of Fourier integral operators of order zero which agrees with $\Pi_{j=o}^N$ up to an
error which also vanishes to order $N$.  Finally, composing with $j^*$ just sets
$s'=s$ in the kernel, which visibly remains  Fourier integral, with phase function parametrizing
the graph of a canonical transformation $\chi.$
 Invertibility  then follows from
the fact that $\tilde{W}^*{W}$ is a pseudodifferential operator with principal symbol the square
of $j^*(\sigma_V)$, which is easily seen to equal 1.      

The proposition  follows by expressing
$$TW_{\gamma}T^{-1} = T \mu(a)^*\tilde{T}^{-1} \tilde{T} W \tilde{T}^{-1}\tilde{T}\mu(a)T^{-1}$$
and noting that $\tilde{T}\mu(a)T^{-1}$ is also a standard Fourier Integral operator. \qed

We now give the proof of the quantum normal form Theorem B for $\sqrt{\Delta}$, stated
in an equivalent form in terms of $W_{\gamma}.$  
 As in the introduction,
the notation $A\equiv B$ means that the complete (Weyl) symbol of $A-B$ vanishes to
infinite order at $\gamma.$  We also use the notation $O_j\Psi^m$ for pseudodifferential
operators of order m whose Weyl symbols vanish to order j at $(y,\eta)=(0,0)$.  Here,
pseudodifferential operator could mean of the standard polyhomogeneous kind, or of
the mixed polyhomogeneous-isotropic kind as in $\Psi^k(S_L^1)\otimes {\cal W}^l$, in which
case the total order is defined to be $m=k+l$. To simplify notation, we will denote the
space of mixed operators of order m by $\Psi_{mx}^m(S^1_L\times \R^n)$.
\medskip
 
\noindent(3.11) \;{\bf Theorem B} {\it $\;\;$  Let $TW_{\gamma}T^{-1}$ be the Fourier Integral
operator of Proposition (3.4), 
defined over a conic neighborhood of $R^+\gamma$ in $T^*(S^1_L \times \R^n)$.  Then:
   $$W_{\gamma}^{-1} T^{-1} \sqrt{\Delta} T W_{\gamma} \equiv
 P_{1}({\cal L},I_{\gamma1},...,I_{\gamma n}) 
+ P_{o}({\cal L},I_{\gamma1},...,I_{\gamma n}) + \dots,$$
where
$$P_1({\cal L}, I_{\gamma1},...,I_{\gamma n})\equiv {\cal L} +
\frac{ p_{1}^{[2]}(I_{\gamma1},...,I_{\gamma n})}{L {\cal L}} 
+ \frac{p_{2}^{[3]}(I_{\gamma1},...,I_{\gamma n})}{(L{\cal L})^2} + \dots \leqno(3.14b)$$
$$P_{-m}({\cal L}, I_{\gamma1},...,I_{\gamma n})\equiv \sum_{k=m}^{\infty} 
 \frac{p_{k}^{[k-m]}(I_{\gamma1},...,I_{\gamma n})}{(L{\cal L})^j}$$
with $p_{k}^{[k-m]}$, for m=-1,0,1,..., homogenous of degree l-m in the variables
$(I_{\gamma1},.
..,I_{\gamma n}) $ and of weight -1.  }
\medskip

\noindent{\bf Proof}:  

 As a semi-classical expansion in the ``parameter" $h = \frac{1}{L {\cal L}}$,
 (3.3) may be rewritten in terms of $\sqrt{\Delta}$ :
$$W_{\gamma}^{-1} T^{-1} \sqrt{\Delta} T W_{\gamma} \sim
 {\cal L} +\frac{ p_{1}(I_{\gamma1},...,I_{\gamma n})}{L{\cal L}} 
+ \frac{p_{2}(I_{\gamma1},...,I_{\gamma n})}{(L{\cal L})^2} + \dots.\leqno(3.12)$$ 
From the fact that the numerators
$f_j(I_{\gamma 1},...,I_{\gamma n})$ in (3.3) are polynomials of degree j+2 (2.36viii)
and of weight -2, the
numerators $ p_{k}(I_{\gamma1},...,I_{\gamma n})$ are easily seen to be
 polynomials of degree $k+1$ and of weight -1. Hence they may be expanded in homogeneous terms 
$$p_k = p_k^{[k+1]} + p_k^{[k]} + \dots p_k^{[o]}, \leqno(3.13)$$
 with  $p_k^{[j]}$ the term of degree j and still of weight -1.  The right side of (3.12) can
then be expressed as a sum of homogeneous operators:
$$P_1({\cal L}, I_{\gamma1},...,I_{\gamma n}) + P_o({\cal L}, I_{\gamma1},...,I_{\gamma n})
+ \dots\leqno(3.14a)$$
with
$$P_1({\cal L}, I_{\gamma1},...,I_{\gamma n})\equiv {\cal L} +
\frac{ p_{1}^{[2]}(I_{\gamma1},...,I_{\gamma n})}{L{\cal L}} 
+ \frac{p_{2}^{[3]}(I_{\gamma1},...,I_{\gamma n})}{(L{\cal L})^2} + \dots \leqno(3.14b)$$
$$P_{-m}({\cal L}, I_{\gamma1},...,I_{\gamma n})\equiv \sum_{k=m}^{\infty} 
 \frac{p_{k}^{[k-m]}(I_{\gamma1},...,I_{\gamma n})}{(L{\cal L})^k}.$$

 The remainder term in (3.12)  can be described as follows:
$$W_{\gamma}^{-1} T^{-1} \sqrt{\Delta} T W_{\gamma} -
[ {\cal L} +\frac{ p_{1}(I_{\gamma1},...,I_{\gamma n})}{L{\cal L}} 
+ \frac{p_{2}(I_{\gamma1},...,I_{\gamma n})}{(L{\cal L})^2} +\dots
+\frac{p_{m}(I_{\gamma1},...,I_{\gamma n})}{(L{\cal L})^m}]\leqno(3.15)$$
$$ \in \oplus_{k=o}^{m+1}
O_{2(m+1-k)}\Psi^{1-k}_{mx}(S^1_L\times \R^n).$$

The remainder terms in (3.14b) are given by:
$$P_1({\cal L}, I_{\gamma1},...,I_{\gamma n}) - [{\cal L} +
\frac{ p_{1}^{[2]}(I_{\gamma1},...,I_{\gamma n})}{L{\cal L}} 
+ \dots + \frac{p_{N}^{[N+1]}(I_{\gamma1},...,I_{\gamma n})}{(L{\cal L})^k}]\leqno(3.16a)$$
$$ \in O_{2(N+2)}\Psi^1_{mx}(S^1_L\times \R^n)$$

$$P_{-m}({\cal L}, I_{\gamma1},...,I_{\gamma n}) - \sum_{k=m}^{N} 
 \frac{p_{k}^{[k-m]}(I_{\gamma1},...,I_{\gamma n})}{(L{\cal L})^k}\in O_{2(N+1-m)}
\Psi^{-m}_{mx}(S^1_L\times \R^n).\leqno(3.16b)$$

Hence the expansion (3.12) is also asymptotic in the sense of $\equiv.$  For the statement of Theorem B
in the introduction, we only need to conjugate under $\mu(\tilde{a}).$\qed

\noindent{\bf Remark} \;\;The remainder in (3.15) could be equivalently described by
saying that its complete symbol lies in the symbol class $S^1_{cl} \cap S^{1, 2(m+1)}$
of  Boutet-de-Monvel [BM].

\section{Residues and Wave Invariants: Proof of Theorem C}

The aim now is to use the normal form of $\sqrt{\Delta}$ near $\gamma$ to calculate
the wave invariants $a_{\gamma k}$ associated to $\gamma$.  In terms of the model
(see the statement of Theorem B in the introduction) we may write the normal form as
$${\cal D} = \overline{\psi}(I_1,...,I_n)[{\cal D}_N + {\cal B}_N] \leqno(4.1)$$
with 
$${\cal D}_N:= {\cal R} + \frac{p_1(I_1,...,I_n)}{L{\cal R}} + \dots +
\frac{p_N( I_1,...,I_n)}{(L{\cal R})^N}$$
and with
$${\cal B}_N  \in \oplus_{k=o}^{N+1}
O_{2(N+1-k)}\Psi^{1-k}_{mx}(S^1_L\times \R^n).$$
Our first observation is that we can drop a sufficiently high remainder term ${\cal B}_N$
in the calculation of a given wave invariant.   We prove this by working out, roughly,
 which parts of
the complete symbol of a general first order pseudodifferential operator $P$ of real principal
type contribute to the wave invariant of a given order associated to a non-degenerate
closed bicharacteristic $\gamma.$  

 In the following proposition, $(s, y)$
will denote any coordinates in a tubular neighborhood $U $
of the projection of $\gamma$ to M  with the property that the equation for $\R^+ \gamma$
in the associated symplectic coordinates $(s, \sigma, y ,\eta)$ is $y=\eta=0.$  Also,
 $Op$ we will denote a fixed quantization of symbols in a conic neighborhood (0.1) of
 $\R^+ \gamma$ to
pseudodifferential operators, and  $p(s, \sigma, y, \eta) \sim p_1 + p_o + \dots $ 
 will denote the
complete symbol of $P$.  The Taylor expansion of $p_j$ at $\R^+\gamma$ will
be written in the form, 
$$p_j(s, \sigma, y, \eta) = \sigma^j
 p_j(s, 1, y, \frac{\eta}{\sigma})= \sigma^j (p_j^{[o]} + p_j^{[1]} + \dots)$$
with 
$p_j^{[m]}(s, 1, y, \frac{\eta}{\sigma})$ homogeneous of
 degree $m$ in $(y, \frac{\eta}{\sigma}).$  We  set $P_j:=Op(p_j), P_j^{[m]} := Op(p_j^{[m]})$
and $P^{\leq N}_j = \sum_{m\leq N} P_j^{[m]}$.  Finally,
    we will denote by
 $\tau_{\gamma k}(P)$
the coefficient of $(t - L + i0)^k log (t - L + i0)$ (or of $(t - L + i0)^{-1}$ in the case of $k=-1$)
in the singularity expansion of the microlocal unitary group $e^{itP}$ near 
 $\gamma$.  

\medskip

\noindent(4.2)\;\;{\bf Proposition}  
$$\tau_{\gamma k}(P) = \tau_{\gamma k}(P_1^{\leq 2(k+2)}  + P_o^{\leq 2(k+1)}
 + \dots P_{-k -1}^{o}).$$
\medskip

\noindent{\bf Proof}: \;\;As mentioned several times, the wave invariant $\tau_{\gamma k}(P)$ is
given by the the non-commutative residue
$$\tau_{\gamma k}(P) = res D_t^k e^{itP} |_{t=L} = res P^k e^{itP}|_{t=L} \leqno(4.3)$$
(see [Z.1] or [G.2]).   We first describe how (4.3) leads to a local formula for
 $\tau_{\gamma k}(P)$ in terms of the jets  of the amplitudes and phases of  
a microlocal parametrix
$$F_{\gamma}(t,x,x') = \int_{\R^n} e^{i\phi(t,x,x',\theta)} a(t,x,x',\theta)
\psi(x,\theta) d\theta \leqno(4.4)$$ for $e^{itP}$
near $\{L\} \times R^+\gamma \times R^+\gamma$.  The remainder of the proof will
connect this data to the complete symbol of $P$.

The amplitude $a$ in (4.4) is a classical symbol, hence has the
expansion $ a \sim \sum_{j=0}^{\infty} a_{-j}$ with $a_{-j}$ homogeneous of degree $-j$ for
$|\theta|\geq 1$.  Since the residue or wave invariant depends only on the wave kernel in
a microlocal neighborhood of $\{L\}\times \gamma \times \gamma$ it is unchanged  by
the insertion homogeneous cut-off 
 $\psi$, supported in (0.1) and equal to one on a smaller conic neighborhood of
$\gamma$.    The phase $\phi(t,x,x',\theta)$
 parametrizes a piece ofthe graph of the Hamilton flow $exptH_{p_1}$ of
 $p_1$ near $\{L\}\times \gamma,$ and may be assumed to consist of only $n$ phase
variables.   
  Writing $x=(s,y)$ and $\theta = (\sigma,
\eta)$ as above and choosing $\psi$ of the form $\psi(\frac{\eta}{\sigma})$, we get
$$\tau_{\gamma k}(P) = Res_{s=0} \int_{\R^+}\int_{\R^{n-1}}\int_U a^{(k)}
(L,s,y,s,y, \sigma, \eta) \psi(\frac{\eta}{\sigma})
 e^{i \phi(L,s,y, s,y, \sigma, \eta)}
|\sigma|^{-s} d\sigma d\eta dv(s,y) \leqno(4.5)$$
where for certain universal constants $C_{\alpha \beta}$,
$$a^{(k)} := \sum_{|\alpha| + |\beta| = k} C_{\alpha \beta} (D_t^{\alpha_1}\phi)^{\beta_1} 
\dots (D_t^{\alpha_n}\phi)
^{\beta_n} D_t^{\alpha_{n+1}} a.\leqno(4.5.1)$$  We note here that
 the residue can be calculated using
any gauging of the trace, in particular by powers $|\sigma|^{-s}$ of the elliptic symbol
$\sigma$ in (0.1).
Changing variables to $\tilde{\eta}:= \frac{\eta}{\sigma}$ and denoting by $a_{k-j}^{(k)}$
the term of order $k-j$ in the polyhomogenous expansion of the kth order amplitude
 $a^{(k)}$, and by $B_{\epsilon}$ the ball of radius $\epsilon$ in $\R^n_{\tilde{\eta}}$, we get
$$\tau_{\gamma k}(P) = \leqno(4.6)$$
$$\sum_{j=o}^{\infty} Res_{s=0} \int_U \int_{R^+}\int_{B_{\epsilon}}
 e^{i \sigma \phi(L,s,y,s,y,1,\tilde{\eta})} a_{k-j}^{(k)}(L,s,y,s,y,1,\tilde{\eta})
\psi(s,y, \tilde{\eta})
 \sigma^{n+k-1-j-s} d\sigma d\tilde{\eta} dv(s,y) $$
$$=\sum_{j=0}^{\infty} \int_U\int_{B_{\epsilon}}\psi a_j^{(k)} 
 (\phi(L,s,y,s,y,1,\tilde{\eta}) + i0)^{s + j- n - k} d\tilde{\eta} dv(s,y).$$ 
The non-degeneracy of $\gamma$ as a closed geodesic implies its non-degeneracy as the
critical manifold of $\phi(s,y,s,y,1, \tilde{\eta})$.  Hence we have
 $$\phi(L,s,y,s,y,1, \tilde{\eta}) = \phi(L,s,0,s,0,1,0) + Q_s(y,\tilde{\eta}) + g(s,y,\tilde{\eta})$$
with $Q$ a non-dengenerate symmetric bilinear form and with $g = 0_s(|(y,\tilde
{\eta})|^3). $ Expanding 
$$(Q_s(y,\tilde{\eta}) +i0 + g)^{\lambda} = \sum_{p=o}^{\infty} C_{(\lambda, p)}
 (Q_s (y ,\tilde{\eta}) + i0)^{\lambda - p} g^p$$
we get that
$$\tau_{\gamma k}(P) = \sum_{j,p = o}^{\infty} C_{(s+j-n-k,p)} Res_{s=0}
\int_{\gamma}\int_{|y| < \epsilon}\int_{B_{\epsilon}}
 (Q_s(y, \tilde{\eta}) + i0)^{s + j - n - k - p} (g^p a^{(k)}_{k-j}\psi J(s,y)) ds dy d\tilde{\eta} 
\leqno(4.7)$$
with $J(s,y)$ the volume density and with binomial coefficients $C_{(\lambda,p)}$.
The family of distributions$ (Q_s(y, \tilde{\eta}) + i0)^{\lambda}$ is meromorphic
with simple poles at the points $\lambda = -(n-1) - r, r=0,1,2,\dots$ and with
residue $C_{r} (Q_s^{-1}(D_y,D_{\eta}))^r \delta(y,\eta)$ for certain constants
$C_r$ (see, [G.S., p.276]).  Here we have dropped the tilde in the notation for $\eta.$ 
For $s=0$ there are possible poles when $k + p + 1 - j \geq 0$
with residues
$$C_r C_{(j -n -k,p)} \int_{\gamma} Q_s^{-1}(D_y, D_{\eta})^{k + p + 1 - j}
(g^p a^{(k)}_{k-j} J)|_{(y= \eta = 0)}ds. \leqno(4.8)$$
However, the residue vanishes unless $2(k + p + 1 -j) \geq 3p$ since $g^p$ vanishes to
order $3p,$ constraining the sum in (4.6) to $2j+p \leq 2(k+1).$  Hence we have, with
new constants $C_{jkp}$,
$$\tau_{\gamma k}(P) = \sum_{j,p: 2j + p\leq 2(k+1)} C_{jkp} \int_{\gamma}
Q_s^{-1}(D_y,D_{\eta})^{k+p+1-j}(g^p a^{(k)}_{k-j} J)|_{(y,\eta)=(0,0)} ds. \leqno(4.9)$$

It follows first that $j \leq k+1$, hence that only the terms $a^{(k)}_k,\dots, a^{(k)}_{-1}$ in
the amplitude  $a^{(k)}$ contribute to (4.9).  To determine the
corresponding range of terms $a_j$ in the amplitude $a$, we observe from (4.5.1) that
$$a^{(k)}_{k-j} := \sum_{|\alpha| + |\beta| = k}
 C_{\alpha \beta} (D_t^{\alpha_1}\phi)^{\beta_1} \dots (D_t^{\alpha_n}\phi)
^{\beta_n} D_t^{\alpha_{n+1}} a_{k-j-|\beta|},\leqno(4.10)$$  
hence that only the terms $a_o,a_{-1},\dots,a_{-k-1}$ contribute to (4.9).  

The same kind of calculation determines $res A e^{itP}|_{t=L}$, for
any $k$th order pseudodifferential operator $A$ supported microlocally in (0.1) in
terms of the complete symbol of $A$ and the phases and amplitudes in (4.3).  In particular,
we see that $res D_t^k  A e^{itP}|_{t=L}=0$ if $A$ has order $-(k+2)$.  This implies:
$$res D_t^k e^{itP} = res D_t^k e^{it (P_1 + \dots P_{-(k+1)})} \leqno (4.11)$$
since $e^{itP} = A_k(t) e^{it(P_1 + \dots + P_{-(k+1)})}$ with $A_k(t)-I \in \Psi^{-k-2}.$
Indeed, with $H_k= P_1+ \dots P_{-(k+1)}, V_k= P-H_k, \tilde{B}(t):= e^{-itH_k}B(t)e^{itH_k}$,
we have
$$D_t \tilde{A}_k = \tilde{V}_k(t) \tilde{A}_k(t), \;\;\;\;\;\tilde{A}_k(0)=I.$$
Thus $\tilde{A}_k(t)$ is the time-ordered exponential of a pseudodifferential operator
of order $-k-2$,  hence equal to $I$ modulo $\Psi^{-k-2}$;  consequently so is $A_k(t).$

We next consider which jets $j^m_{\gamma}a_{-i}$ of $a_o, \dots a_{-k-1}$ 
and $j^m_{\gamma}\phi$ of $\phi$ contribute to (4.9).  
Here by $j^m_{\gamma}a_{-i}$ we mean the Taylor polynomial of degree $m$ of $a_{-i}(t, s, y,
s,y, \sigma, \eta)$ in the $(y,\eta)$ variables for $t$ near $L$.  From (4.9-10) it is evident
that the maximum number of $(y,\eta)$ derivatives on $a_{-i}$ occurs when $p=0$ in
(4.9) and when $|\beta| = k, k-j-|\beta| = -i$ in (4.10).  It follows that at most $2(k+1-i)$
derivatives fall on $a_{-i}$ in (4.9).
Also, the maximum number of $(y,\eta)$ derivatives on the phase occurs in terms where a
factor of $g$ is differentiated the maximum number of times; namely in terms with $p=1$
and with $j=0,$ in which the phase is differentiated $2(k+2)$ times.  Hence
$$\tau_{\gamma k}(P) \;\;\;\mbox{depends only on}\;\;\;\;\; j^{2(k+1)}_{\gamma}a_o, 
j_{\gamma}^{2k}a_{-1}, \dots j_{\gamma}^oa_{-(k+1)};\;\;j^{2(k+2)}_{\gamma}\phi. \leqno(4.12)$$

To give bounds on
the jets of the terms $p_j (j=1,0,\dots -(k+1))$ in the symbol of $P$  which contribute to 
the jets $j^{2(k+1-i)}_{\gamma}a_{-i}$, we now have to consider  some details of the construction of the
parametrix (4.3).    
 We first recall that the
amplitudes are obtained by solving transport equations of the form
$$D_t a_{-j} = \sum_{m=0}^j [{\cal P}_{1, m+1} a_{-j+m} + 
\sum_{\nu =o}^{j-m}{\cal P}_{-\nu,m}a_{-j+m+\nu}] \leqno(4.13)$$
where ${\cal P}_{-\nu,m}={\cal P}_{-\nu,m}(\phi, s,y, D_s, D_y)$
 is a differential operator of order $\leq m$ obtained from $Op(p_{-\nu})$ as follows:
$$e^{-i\rho \phi}Op(p_{-\nu}) e^{i \rho \phi} \sim \sum_{m=0}^{\infty}\rho^{-\nu - m}
{\cal P}_{-\nu,m}  \leqno(4.14.1)$$
 $\phi$ being the phase in the microlocal parametrix (see [Tr, Ch.VI,(5.28)]).  Consistently
with (4.11), only the terms
$p_1, p_o, \dots, p_{-(k+1)}$ in the symbol of $P$ contribute to
the transport equation for $a_{-i} (i=o,\dots,-(k+1).$ To 
determine which jet of $p_{-\nu}$ is involved in ${\cal P}_{-\nu,m}$,
 we further recall [loc.cit.,Theorem 3.1] that the expansion in
(4.14.1) is obtained by re-arranging terms in the expansion  
$$e^{-i\rho \phi}Op(p_{-\nu}) e^{i \rho \phi}u \sim \sum_{ \alpha} \frac{1}{\alpha!}
\partial_{\xi}^{\alpha} p_{-\nu}(s,y, \rho d\phi) {\cal N}_{\alpha}(\phi; \rho, D_s,D_y)u
\leqno(4.14.2)$$
with
$${\cal N}_{\alpha}(\phi; \rho, D_x) u = D_{x'}^{\alpha} e^{i\phi_2(x,x')}u(x') |_{x=x'}$$
$$\phi_2(x,x')= \phi(x)-\phi(x') - \langle x-x', \nabla \phi(x) \rangle.$$
In view of the fact that ${\cal N}_{\alpha}$ is a polynomial in $\rho$ of degree $\frac{|\alpha|}
{2}$, we see that $|\alpha| \leq 2m$ in any term of (14.4.2) contributing to ${\cal P}_{-\nu,m}$.
Hence ${\cal P}_{-\nu,m}$ involves at most $2m$ derivatives of $p_{-\nu}$, and so the
transport equation for $a_{-i}$ involves at most $2(i-\nu)$ derivatives of $p_{-\nu}$ 
($\nu=-1, 0, \dots, i.$)  

To draw these kinds of conclusions about the amplitudes $a_{-i}$ themselves, as opposed
to the coefficients in their transport equations, we have to take into account the initial
conditions in the transport equations.  Unfortunately, the defining initial conditions occur
at $t=0$ while we are interested in the long time behaviour at $t=L$.  Were there no 
conjugate pairs along geodesics near $\gamma$, the transport equations could be solved
up to $t=L$ and since the initial conditions can be taken in the form  $a_o|_{t=0}\equiv 1,
a_{-i}|_{t=0} \equiv 0$ we could conclude that the $a_{-i}$'s  are integrals of polynomials
in at most $2(i-\nu)$ derivatives of $p_{-\nu}$ and in at most $2(i+1)$ derivatives of the
phase (coming from the ${\cal N}_{\alpha}$'s).  However, in the case of elliptic closed
geodesics, there will always be conjugate pairs for $L$ sufficiently large and we cannot
construct the parametrix so simply.  Rather we will use the group property $U(L) = U(L/N)^N$
with N sufficiently large that a parametrix for $U(L/N)$ can be constructed by the geometric
optics method.  We then have to eliminate all but $n$ phase variables in the power, which
will lead to crude but serviceable bounds on the order of the jets.

We therefore begin with the construction of a short time parametrix $F_{\delta k}$,
valid for $|t| < \delta,$ of the form (4.4) with
  phase  $\phi(t,x,x',\theta)= S(t, x, \theta) - \langle x', \theta \rangle,$
satisfying $\partial_t S + p_1(x, d_x S) = 0, S|_{t=o} = \langle x, \theta \rangle$, and with
amplitude $a\sim \sum_{j=0}^{k+1}a_{-j}$ satisfying (4.13) on  $|t| \leq \delta$ 
 and with initial conditions
 $a_o|_{t=0}=1, a_{-j}|_{t=0} = 0 (j >0)$.  
By the observations above, we  have
that  $a_{\delta, -i}(t,x,x',\theta)$ 
involves at most $2(i-\nu)$ derivatives of $p_{-\nu} (\nu \leq i)$ for $t$ in this interval.
 As is well-known, we then have
 $e^{itP} - F_{\delta k}(t) = R_{ok}(t)$ with
the remainder $R_{\delta k}(t)$ a Fourier integral operator of order $-(k+2)$ associated
 to the graph of
$exptH_{p_1}$ for $|t|< \delta$. 
Since $F_{\delta k}$ is of order zero, it follows that for any $N$ and for $|t| < \delta,$
$e^{itNP} \equiv  F_{\delta k}^N(t)$ modulo Fourier integral operators
 of order $-k-2$ in this class.  Hence, the desired parametrix (4.4) can be constructed by choosing $N$ so that
$L/N < \delta$ and by eliminating all 
 but $n$ phase variables in $F_{\delta k}^N$ 
by the stationary phase method. 
Just as in (4.11) the remainder terms with a factor of $R_{\delta k}$ will not contribute to 
$\tau_{\gamma k}(P)$, nor will terms $a_{-i}$ in the final amplitude with $i> k+1.$

To complete the proof of the proposition, we have to count the number of $(y,\eta)$-derivatives
of $p_{-\nu}$ which  appear in the terms
$$D_t^{k-|\beta|}Q_s^{-1}(D_y,D_{\eta})^{k+1-j}
\hat{a}_{k-j-|\beta|}|_{(y,\eta)=(0,0)},\leqno(4.15.1)$$
and the number of derivatives  of the final phase (hence of $p_1$) that occur in the terms
$$Q_s^{-1}(D_y,D_{\eta})^{k+p+1-j}g^p(D_t^{\alpha_1} \phi)^{\beta_1} \dots
(D_t^{\alpha_n} \phi)^{\beta_n}D_t^{\alpha_{n+1}}
 \hat{a}_{k-j-|\beta|}|_{(y,\eta)=(0,0)}\leqno(4.15.2)$$
with $\hat{a}_{-i}$  the terms in the final amplitude.
  We start from the fact that, in an obvious notation,
$$F_{k \delta}^N=\int_{\R^{nN}}\int_{M^{N-1}} e^{i (\phi_1 + \dots \phi_N)} 
\sum_{(i_1,\dots,i_N) \in {\bf N}^N}
 a_{-i_1}\dots a_{-i_N} d\theta_1 \dots d\theta_Ndv(x_1) \dots dv(x_{N-1})\leqno(4.16)$$
and that the elimination of all but $n$ phase variables by the stationary phase method gives
essentially the same formulae for the final amplitudes $\hat{a}_{-i}$ as in the case of
non-homogeneous Lagrangean distributions of the form $\int u(x,y) e^{i\omega f(x,y)} dx$
with $f(x,y)$ a real-valued function defined near $(0,0)$, with $f_x'(0,0)=0$, and with
$f''_{xx}(0,0)$ non-singular.  In this situation, we have (see [HoI, Theorems 7.7.5-6])
$$\int u(x,y)e^{i\omega f(x,y)}dx = C\frac{ e^{i \omega f(x(y),y)} }
{|det(\omega f''_{xx}(x(y),y))|^{-\half}} \sum_{j+ \frac{n}{2} < M} L_{f,j,y} u(x(y),y) \omega^{-j}$$
modulo terms of order $0(\omega^{-(M)})$, with $L_{f,j,y}$ a differential operator of
order 2j whose coefficients involve at 2j derivatives of $f''$.  It follows that $\hat{a}_{-i}$
has the form (with $I=(i_1,\dots, i_N), \Phi = \phi_1 + \dots \phi_N)$)
$$\sum_{I,j: |I| + j = i} L_{\Phi,j} a_{-i_1}\dots a_{-i_N}.$$
Since also $a_{-i}$ has the form
$$F_i(p_1,Dp_1,\dots,D^{2(i+1)}p_1, p_o,\dots, D^{2i}p_o,\dots, p_{-i}; \phi,\dots,
D^{2i}\phi)$$
with $F_i$  a multiple integral of polynomials, 
and with $D^k$ some differential operators of degree k,
 we see that also $\hat{a}_{-i}$ has the form
$${\cal F}_i(p_1,\dots,D^{2(i+1)},p_o,\dots,D^{2i}p_o,\dots,D^2p_{-i+1}, p_{-i};
\phi, \dots, D^{2(i+1)}\phi).$$
Hence as regards  number of derivatives of $p_1,\dots,p_{-i},\phi,$ the $\hat{a}_{-i}$'s
behave exactly as the $a_i$'s, that is, involve at most $2(i-\nu)$ derivatives of $p_{-\nu}.$
Hence (4.9)-(4.10) apply to the final amplitudes, and we conclude that 
$$\tau_{\gamma k}(P)\;\;\; \mbox{depends only on}\;\;\;\;\; j^{2(k+2)}_{\gamma}p_1, 
j_{\gamma}^{2(k+1)}p_o, \dots j_{\gamma}^op_{-(k+1)};\;\;j^{2(k+2)}_{\gamma}\phi. $$
Since $j^{2(k+2)}\phi$ depends only on $j^{2(k+1)}p_1$, the proof of the proposition is
complete..\qed
\medskip

\noindent(4.16) {\bf Corollary} 
 $$\tau_{\gamma k} ({\cal D}) = \tau_{\gamma k}({\cal D}_{k+1}).$$
\medskip

The following Lemma combines the previous results in a form applicable to the
calculation of the wave invariants.  The notation $TrAP^{-s}$ (with $P\in \Psi^{1}$
elliptic) is short for the zeta function obtained by
meromorphic  continuation of the trace from $Res >>0$ to $\C.$
\medskip

\noindent(4.17) {\bf Lemma}
$$\tau_{\gamma k}(\sqrt{\Delta}) = Res_{s=0} D_t^k
Tr\overline{\psi}({\cal R}, I_1,\dots,I_n) e^{it {\cal D}_{k+1}} {\cal R}^{-s}|_{t=L}$$
\medskip

\noindent{\bf Proof}\;\; Since ${\cal R}$ is elliptic in the essential support
 of $\overline{\psi}$, the trace on the right side is well defined and has a meromorphic
continuation to $\C$ (cf. [G.2], [Z.1,5]).  The residue is of course $\tau_{\gamma k}({\cal D})$
by the previous proposition.  Hence it suffices to show that 
$$\tau_{\gamma k}(\sqrt{\Delta}) = \tau_{\gamma k}({\cal D}). \leqno (4.18)$$
This however follows from the previous proposition combined with Theorem B:  Indeed,
the proposition shows that $\tau_{\gamma k}(A) = \tau_{\gamma k}(B)$ if $A\equiv B$
in the sense of Theorem B.  Also, $\tau_{\gamma k}(W A W_1) = \tau_{\gamma k}(A)$
if $W_1$ is a parametrix for $W$ on the essential support of $A$, as follows from
the tracial property $res WAW_1 = res W_1 W A$  of the residue (see [Z.1,5] for instance).\qed
\medskip

We come now to the calculation of the wave invariants as residue traces of the normal form wave group.
But we will simplify (4.17) further before evaluating the residue trace.  First, we rewrite
everything in terms of $D_s$, and $H_{\alpha}$ using (0.3-4).  Since
$$\frac{p_{\nu}(I_1,\dots,I_n)}{(L{\cal R})^{\nu}} = \frac{p_{\nu}(I_1,\dots,I_n)}
{(LD_s)^{\nu}} (I -\nu \frac{H_{\alpha}}{L D_s} +\half \nu (\nu-1) 
(\frac{H_{\alpha}}{L D_s})^2 + \dots)$$
and since we can drop the  $D_s^{-(k+1) + \nu}H_{\alpha}^{k+1 -\nu}$ and higher terms
by Proposition (4.2), ${\cal D}_{k+1}$ can be written in the form 
$${\cal D}_{k+1} \equiv D_s + H_{\alpha} + \frac{\tilde{p}_1(I_1,\dots,I_n,L)}{L D_s} +
 \frac{\tilde{p}_2(I_1, \dots,I_n,L)}{(LD_s)^2}
+\dots+\frac{\tilde{p}_{k+1}(I_1,\dots,I_n,L)}{(LD_s)^{k+1}} \leqno(4.19)$$
modulo terms which make no contribution to $\tau_{k\gamma}$.  Secondly, we can use
$L D_s$ rather than ${\cal R}$ as the gauging elliptic operator in (4.17).  To simplify the
notation we will denote all but the first two terms on the right side of (4.19) by 
$${\cal P}_{k+1}(D_s, I_1,\dots,I_n,L):=\frac{\tilde{p}_1(I_1,\dots,I_n,L)}{L D_s} +
 \frac{\tilde{p}_2(I_1, \dots,I_n,L)}{(LD_s)^2}
+\dots+\frac{\tilde{p}_{k+1}(I_1,\dots,I_n,L)}{(LD_s)^{k+1}}. \leqno(4.20)$$
Then we have:
$$\tau_{\gamma k}(\sqrt{\Delta}) = Res_{z=0} Tr D_t^k\overline{\psi}(D_s,I_1,\dots,I_n)
e^{it[\frac{1}{L} (2 \pi L D_s + H_{\alpha} )+ {\cal P}_{k+1}]} (LD_s)^{-z} |_{t=L}. \leqno(4.21)$$
We can now give:
\medskip

\noindent{\bf Proof of Theorem C:}  By (4.21) and the fact that $e^{2\pi i L D_s} \equiv I$ on 
$L^2(S^1_L)$ we get
$$a_{k \gamma} = \tau_{\gamma k}(\sqrt{\Delta}) =$$
$$Res_{z=0} Tr [\frac{1}{L}2 (\pi L  D_s + H_{\alpha}) + {\cal P}_{k+1}]^k e^{i
H_{\alpha}}
 e^{iL{\cal P}_{k+1}}\overline{\psi}(\frac{I}{\epsilon D_s}) (LD_s)^{-z} \leqno(4.22)$$
with 
$$\overline{\psi}(\frac{I}{\epsilon D_s}):=\Pi_{j=1}^n\overline{\psi}(\frac{I_j}{\epsilon D_s}).$$
From the well-known spectra of $D_s, I_1,\dots,I_n$ we can rewrite (4.22) as
$$Res_{z=0}  \sum_{m=1} ^{\infty}  m^{-z}
\{ \sum_{(q_1,\dots,q_n) \in {\bf N}^n} 
[\frac{1}{L} (2\pi m + \sum_{j=1}^n (q_j + \half)\alpha_j) + {\cal P}_{k+1}
(m, q_1 + \half,\dots, q_n + \half)]^k \cdot \leqno(4.23) $$
$$ \cdot e^{iL{\cal P}_{k+1}(m, q_1 + \half,\dots, q_n + \half,L)}\overline{\psi}(\frac{q}{\epsilon m})
e^{i \sum_{j=1}^n (q_j + \half) \alpha_j}\} .$$
Regarding
$${\sum}_{(q_1,\dots,q_n) \in
{\bf N}^n}  \overline{\psi}(\frac{q}{\epsilon m})e^{i \sum_{j=1}^n (q_j + \half) \alpha_j}
\leqno(4.24)$$
as a smooth function of $\alpha \in \R^n$, we can further rewrite (4.23) as
$$Res_{z=0}  \sum_{m=1} ^{\infty}  m^{-z} \{
[\frac{1}{L} (2\pi m + \sum_{j=1}^n\alpha_j D_{\alpha_j}) + {\cal P}_{k+1}
(m, D_{\alpha_1},\dots, D_{\alpha_n},L)]^k \cdot$$
$$ e^{iL{\cal P}_{k+1}(m, D_{\alpha_1},\dots, D_{\alpha_n},L)}\sum_{(q_1,\dots,q_n) \in
{\bf N}^n} \overline{\psi}(\frac{q}{\epsilon m})
e^{i \sum_{j=1}^n(q_j + \half) \alpha_j}\} .\leqno(4.25)$$
 Since $ \overline{\psi}(\frac{q}{\epsilon m})$
is for each $m$ a finitely supported function of $q$, we can also rewrite (4.24) as
$$ \lim_{\delta \rightarrow 0} \Pi_{j=1}^n \overline{\psi}(\frac{D_{\alpha_j}}{\epsilon m})
\sum_{q_j = 0}^{\infty} e^{i \sum_{j=1}^n\alpha_j (q_j + \half) -\delta}$$
$$= \lim_{\delta \rightarrow 0} \Pi_{j=1}^n \overline{\psi}(\frac{D_{\alpha_j}}{\epsilon m})
\Pi_{j=1}^{n}\frac{e^{i\half \alpha_j}}{1- e^{i(\alpha_j +i \delta})}.\leqno(4.26)$$
Recalling the definition of $T(\alpha)$ (0.12a), combining (4.23)-(4.26), and
taking the limit as $\delta \rightarrow 0$, we get:
$$a_{\gamma k} = Res_{z=0} \sum_{m=1}^{\infty} m^{-z} 
\{[\frac{1}{L} (2\pi m + \sum_{j=1}^n \alpha_j D_{\alpha_j}) + {\cal P}_{k+1}
(m, D_{\alpha_1},\dots, D_{\alpha_n},L)]^k \cdot$$
$$ e^{iL{\cal P}_{k+1}(m, D_{\alpha_1},\dots, D_{\alpha_n},L)} \Pi_{j=1}^n 
\overline{\psi}(\frac{D_{\alpha_j}}{\epsilon m}) T(\alpha)\}. \leqno(4.27)$$
Here, we use the fact that $T(\alpha)$ is a
a temperate distribution on $\R^n_{\alpha}$ with singular support on $\cup_{j=1}^n
  \R \times\R \dots \times 2 \pi \Z \dots \times \R^n$ (the factor of $\Z$ occurring in the jth
position) and that the limit $\delta \rightarrow 0$ in (4.26) can be taken in ${\cal S}'$.  Since
the cut-off smooths out the singularity, the right side of (4.27) is a smooth function of
$\alpha$ and can be evaluated at the special values of $\alpha$ determined by $\gamma.$

Since ${\cal P}_{k+1}(m, D_{\alpha_1},\dots, D_{\alpha_n},L)$ is a symbol of order $-1$
in $m$ with coefficients given by polynomials in the operators $D_{\alpha_j}$, 
we can expand the kth power in (4.27)   as an
operator-valued polyhomogeneous function of $m$.  At least formally, we can also expand the
exponential $e^{iL{\cal P}_{k+1}(m, D_{\alpha_1},\dots, D_{\alpha_n},L)}$ in a power
series and then expand each term in the power series as a polynomial in $m^{-1}$.  Collecting
powers of m, the right side of (4.27) can be  put, at least formally, in the form 
$$Res_{z=0} \sum_{m=1}^{\infty} \sum_{j= o}^{\infty} m^{-z+ k -j} {\cal F}_{k,k-j}(D_{\alpha},L)
\overline{\psi}(\frac{D_{\alpha}}{m\epsilon}) T(\alpha),  \leqno(4.28)$$
with ${\cal F}_{k,k-j}(D_{\alpha},L)$ the coefficient of $m^{k-j}$ in (4.27).  To justify
this manipulation, we have to deal with the remainder term in the Taylor expansion of
the exponential.  We have
$$e^{ix} = e_N(ix) + r_N(ix)$$
 $$e_N(ix) = 1 + ix + \dots+ \frac{(ix)^N}{N!}$$
$$r_N(ix)=(ix)^{N+1}b_N(ix)$$
$$b_N(ix) = \int_o^1\dots\int_o^1 t_N^N t_{N-1}^{N-1}
 \cdots t_o^o e^{t_N t_{N-1} \cdots t_o ix}
dt_o \cdots dt_N.$$
Hence we have
$$e^{iL{\cal P}_{k+1}}:=  e_N(i L {\cal P}_{k+1}) +
(i L {\cal P}_{k+1})^{N+1} b_N(i L {\cal P}_{k+1})$$
with ${\cal P}_{k+1}$ short for 
${\cal P}_{k+1} (m, D_{\alpha_1}, \dots, D_{\alpha_n})$.
The $e_N$ term of course contributes a finite number of terms of the desired form (4.28).
For the remainder,
we expand  $(i L {\cal P}_{k+1})^{N+1}$ as a polynomial in $m^{-1}$ 
with coefficients given by operators $Q_{N p}(D_{\alpha_1},\dots,D_{\alpha_n})$ 
and observe that
each term has a factor of $m^{-N-1}$.  For each such term, we remove the coefficient operator
$Q_{N p}$
from the sum $\sum_m$, leaving only the factor of $b_N$.  We then rewrite the resulting sum
as a double sum $\sum_{m q}$ as in (4.23), replacing all operators in $D_{\alpha_j}$ by
their eigenvalues.
  Since $b_N(ix)$ is a bounded function and since
each term of the resulting sum has at least the factor  $m^{-z -N -1 + k}$ (possibly
multiplied by a negative power of $m$), we see that the remainder is a sum of terms
of the form
$$Res_{z=0} Q_{N p}( D_{\alpha_1},\dots,D_{\alpha_n}) \sum_{kq} m^{-z +k -N -1 -l}
b_N(i{\cal P}_{k+1}(m, q+ \half)) \overline{\psi}(\frac{q}{\epsilon m}) e^{i\langle
q+\half,
\alpha \rangle }.\leqno(4.29)$$
We then observe that the sum is bounded by $\sum_{m=1}^{\infty} m^{-Re z - N - 1 + k +n}$,
hence converges absolutely and uniformly for $Re z > -N + k + n$.  It follows that
for $N > (n+k)$ the sum in (4.29) defines a holomorphic function of $z$ in a half-plane
containing $z=0$ and since the operations of taking the residue in $z$ and  derivatives
in $\alpha$ commute, each term (4.29) is zero.  This justifies (4.28) and shows that 
 it is actually a finite sum in j, say $j < M$ (in fact M=(k+1)(n+k+1)).

The residue in (4.28) is therefore well-defined  and independent of $\epsilon$. Since
 $\overline{\psi}(\frac{D_{\alpha}}{m\epsilon}) T(\alpha) \rightarrow T(\alpha)$ as $\epsilon \rightarrow
\infty$  we must have
$$a_{\gamma k} =Res_{z=0} \sum_{m=1}^{\infty} \sum_{j= o}^{M} m^{-z + k -j} 
{\cal F}_{k,k-j}(D_{\alpha},L)
T(\alpha)  \leqno(4.30)$$  
$$=Res_{z=0} \sum_{j=0}^{M} \zeta(z + j - k){\cal F}_{k,k-j}(D_{\alpha},L) T(\alpha).$$
Here, $\zeta$ is the Riemann zeta-function, which has only a simple pole at $s=1$ with
reside equal to one.  It follows that the only term contributing to (4.29) is that with
$j = k + 1$ and hence we have
$$a_{\gamma k} =  {\cal F}_{k, -1}(D_{\alpha},L) T(\alpha) \leqno(4.30),$$
 completing the proof of Theorem C.\qed

\section{Local formulae for the residues: Proof of Theorem A}

The characterization of the wave invariants in Theorem A is reminiscent of that
of the heat invariants in [ABP] or [Gi], but involves non-local metric invariants 
near $\gamma.$  We begin by determining the metric data which contribute
to $a_{\gamma k}$.

As always, we assume that $\gamma$ denotes a primitive
closed geodesic, and denote the $\ell$th iterate of $\gamma$ by $\gamma^{\ell}$.  As
above, we use the exponential map along $N_{\gamma}\sim S^1_L \times \R^n$ to pull back the metric to a
metric on a neighborhood of $S^1_L\times \{0\}$, or more simply $S^1$, in $S^1_L\times \R^n$
with the same wave invariants as $g$ along $\gamma$.  This reduces the theorem to the case 
$M = S^1_L \times \R^n$, with  $S^1$  a closed geodesic of the metric.

We let $J^m_{S^1}$ denote the manifold of $m$-jets along $S^1$ of Riemannian metrics
 on $S^1\times \R^n$ with the property that $S^1$ is a closed geodesic.  We also let
$J^m_{S^1 ell}$ denote the open submanifold of m-jets of metrics 
for which $S^1$ is non-degenerate elliptic.  We also write $J^m_{\gamma ell}$ when
we wish to identify $\gamma$ and  $S^1$.
The density $I_{\gamma k}(s)ds$ of the $k$th wave invariant is then a map
$$I_{\gamma k} : J^{m_k}_{\gamma ell} \rightarrow \Omega^1(S^1)$$
where $\Omega^1(S^1)$ is the space of smooth densities along $S^1$ and where the
jet order $m_k$ will be shown below to be $m_k = 2k+4$.  Since $I_{\gamma k}
(s,g)ds$
is independent of the choice of coordinates on $S^1 \times \R^n$ we may express it
in terms of Fermi normal coordinates $(s,y)$ with respect to a fixed normal frame
$e(s)$ for g.  As usual, the
metric coefficients $g_{ij} (i,j=o,1,\dots,n)$ will be understood relative to these
coordinates. 
\medskip

\noindent{\bf \S5.1 The metric data in $I_{k \gamma}$}
\medskip

We now claim that $I_{k \gamma}$ is an invariant polynomial in the following data:

(i) the curvature tensor $R$ and its covariant derivatives $\nabla^m R$ with $m \leq 2k+2$,
contracted with respect to the  Fermi normal vector fields $\frac{\partial}{\partial s}$ and
$\sum_{j=1}^n c_j \frac{\partial}{\partial y_j}$ with $c_j \in \R$;

(ii) the components $Rey_{ij}(s), Im y_{ij}(s)$ of the normalized eigenvectors $Y_j \in
{\cal J}^{\bot}_{\gamma}\otimes \C$ of $P_{\gamma}$ relative to $e_1,\dots,e_n$,
and their first derivatives;

(iii) at most $2k+1$ indefinite integrals over $S_L^1$ of (i)-(ii).

(iv) the length $L$ and inverse length $L^{-1}$ of $\gamma^{\ell}$;
 
(v) the Floquet invariants $\beta_j = (1 - e^{i\alpha_j})^{-1}$;

Indeed,  by Theorem C, $I_{\gamma,k}(s,g)ds$ is a  density depending
 only  on the data contained in $T(\alpha), L$ and in
the coefficients of ${\cal F}_{k,-1}$.  The latter is identical to the data in the coefficients
in $\tilde{p}_1, \dots, \tilde{p}_{k+1}$, hence to that in $p_1,\dots,p_{k+1}$.  By
Proposition (4.2), these depends only on $j^{2k+4}_{\gamma} g$, and on the
coefficients of the symplectic matrices $r_{\alpha}(s), a_s$ which arise in the intertwining
of $\Delta_h$ to its normal form.  The coefficients of $r_{\alpha}$ depend only on L and  the
$\alpha_j$'s and those of $a_s$ depend only on the Jacobi field components $y_{ij},$ and
their first derivatives $\dot{y}_{ij}$.  The second and higher derivatives of $y_{ij}$
 can of course be eliminated by means of the Jacobi equation.   Regarding (ii), we 
 recall (see [Gr, Theorem 9.21]) that the Taylor series
expansions at $y=0$ of the  $g_{ij}(s,y)$'s, of the co-metric coefficients $g^{ij}(s,y)$, and of the
volume densities and their powers, involve only the curvature tensor $R$ and its covariant
derivatives $\nabla R,\dots, \nabla^m R$ contracted with respect to the normal Fermi
vector fields.  Since the coefficients ${\cal L}_{2-\frac{m}{2}}$ of the semi-classical
expansion of $\Delta_h$ depend only on this data, and since the intertwinings only introduce
coefficients in  data  (ii)-(iv), the $p_k(I)$'s can only involve the $g_{ij}$'s thru
 the curvature and its  covariant derivatives.

For the fact that $I_{\gamma k}(s)$ is a polynomial in (i)-(iv), it suffices to reconsider
the construction of the normal form.  Obviously the coefficients of the ${\cal L}$'s and
their metaplectic conjugates the ${\cal D}_{\frac{j}{2}}$'s  are polynomials in (i)-(iv).  It then
 follows from (2.39c)
 that the $f_j(I)'s$, hence  the $p_j(I)$'s, are also polynomials in this data.  Indeed, we argue
inductively from the
construction of the normal form that 
  the coefficients change in the step from ${\cal D}^{\frac{m}{2}}_{\frac{j}{2}}$
to ${\cal D}^{\frac{m+1}{2}}_{\frac{j}{2}}$ as a result of taking commutators with
operators whose coefficients are polynomials in (i)-(iv).  It follows that the coefficients
in ${\cal D}^{\frac{m+1}{2}}_{\frac{j}{2}}$ are also polynomials in this data, 
with possibly two more $D_s$-derivatives
due to the commutators with $D_s$.  Since the Jacobi data $y_{ij}, \dot{y}_{ij}$ enters
in thru a linear change of variables, the degree of the polynomial in this data will be
closely related to the degree in the $(y,\eta)$ variables, which is twice the degree in the
$I$ variables.  The degree of the polynomial in the Jacobi data is however not necessarily
the same as that (namely, k+2) in the $I$-variables: from the proof of Lemma (2.22),
we see that the algorithm for computing the polynomials $f_j(I_1,\dots,I_n)$  involves
taking operator
commutators (or Poisson brackets of symbols); this  lowers the order in the $I$-variables
but not in  the Jacobi data which are coefficients of the polynomials in the $I$'s.  We
will show below (see `Jacobi degrees') that the order in the Jacobi data is no more
than $6k + 6$.
\medskip

\noindent{\bf \S5.2 Weights of $I_{k\gamma} $ and of the data}
\medskip

To determine the weights of the polynomials in the data (i)-(iv), we now
extend (and in part recall) the table in \S 1.4 describing how the various 
 data transform under the metric rescaling $g\rightarrow g_{\epsilon}:=\epsilon^2 g.$ As above, $(s,y)$
always refer to Fermi normal coordinates relative to $g$, and for notational simplicity we put $s=y_o.$
The notations $\nabla, R,\Delta$ refer respectively to the Riemannian connection, curvature and Laplacian.
We distinguish $\Delta$ from the local expression
 $\frac{1}{\sqrt g} \sum_{i,j=o}^n \partial_{y_i} g^{ij} \sqrt{g}
\partial_{y_j}$ for $\Delta$ relative to the Fermi normal coordinate frame. 
\bigskip

\begin{tabular}{r||l} $g$ & $\epsilon^2 g$ \\ \hline
$g_{ij} = g(\partial_{y_i},\partial_{y_j})$ & $g_{ij} = \epsilon^2 g(\epsilon^{-1}\partial_{y_i},
\epsilon^{-1}\partial_{y_j})$ \\ \hline
$ D^{\alpha}_{s,y} g_{ij}$ & $\epsilon^{-|\alpha|}D^{\alpha}
_{s,y} g_{ij}$   \\ \hline
$\nabla,\;\;\;\;\;R,\;\;\;\;\Delta$ &$\nabla,\;\;\;\;\;R,\;\;\;\;\epsilon^{-2}\Delta$ \\ \hline
$\frac{1}{\sqrt g} \sum_{i,j=o}^n \partial_{y_i} g^{ij} \sqrt{g}
\partial_{y_j}$ & $\epsilon^{-2} \frac{1}{\sqrt g} \sum_{i,j=o}^n \partial_{y_i} g^{ij} \sqrt{g}
\partial_{y_j}$ \\ \hline
 $L, \;\;\;\;\;\;(s,y)$ & $\epsilon L, \;\;\;\;\;\;(\epsilon s, \epsilon y)$\\ \hline
$ y_{ij} := g(Y_i,e_j)$ & $\epsilon^{\frac{1}{2}} y_{ij}
 = \epsilon^2 g( \epsilon^{-\frac{1}{2}} Y_j,\epsilon^{-1}
e_j)$ \\ \hline
$K_{ij} = g(R(\partial_s,e_i)\partial_s,e_j)$ & $\epsilon^{-2} K_{ij}$ \\ \hline
\end{tabular}
\bigskip

The trace of the wave group thus scales as
$$Tr e^{it\sqrt{\Delta}} \rightarrow Tr e^{ i\frac{t}{\epsilon} \sqrt{\Delta}} \leqno(5.2.1)$$
from which it follows that
$$a_{\gamma k} = res \sqrt{\Delta}^{k} e^{i L \sqrt{\Delta}} \rightarrow \epsilon^{-k}a_{
\gamma k}. \leqno (5.2.2a)$$
This can also be seen from the fact that $a_{\gamma k}$ is the coefficient of
$(t - L + i0)^k log (t - L + i0)$ which is homogeneous of degree k modulo smooth functions
of $t$.  Since $a_{\gamma k} = \int_{\gamma} I_{\gamma k}(s) ds$, and the integral 
scales like $ds$, we also have
$$I_{\gamma k} (s, \epsilon^2 g) = \epsilon^{-k+1} I_{\gamma k}(s). \leqno(5.2.2b)$$

As a check on the normal form, let us verify that $wgt(a_{\gamma k})= -k$ directly from a weight
analysis of the normal form.  It is obvious that
${\cal R} = \frac{1}{L}(LD_s +H_{\alpha})$
scales like $\epsilon^{-1}$ and so does  each term in the expansion
$$\sqrt{\Delta} \sim {\cal R} + \frac{p_1(I_1,\dots,I_n)}{L{\cal R}}+\dots.$$
Moreover the $p_k(I_1,\dots,I_n)$'s and ${\cal P}_k$'s are also of weight -1, 
as determined in Lemma (2.22) and Theorem B. Since the transition
from  $p_{\nu}$ to the $\tilde{p}_{\nu+r}$'s only involves multiplication of  $p_{\nu}$
with the $r$th power of $\frac{H_{\alpha}}{LD_s}$ it is obvious that $\tilde{p}_{\nu}$
also has weight -1. Expanding the exponent in the residue calculation, we get that
$$a_{\gamma k} =\sum_{N=0}^{k+1} \frac{i^N}{N!}
res [ \frac{1}{L} (LD_s + H_{\alpha}) + {\cal P}_k]^k (L {\cal
P}_k)^Ne^{iH_{\alpha}} \leqno(5.2.3)$$ 
with $[ \frac{1}{L} (LD_s + H_{\alpha}) + {\cal P}_k]$ of weight -1 and with 
$(L {\cal P}_k)$ of weight 0.  Hence $wgt(a_{\gamma k}) = -k.$
\medskip

\noindent{\bf \S5.3 Wave invariants and QBNF coefficients}
\medskip

  They are related as follows: The 
terms in (5.2.3)  contributing nontrivially  to  $a_{\gamma k}$ have the form
$$L^{-k + j_1 + \dots j_{k+1} }
res (LD_s)^{\ell - (j_1 + 2j_2 + \dots (k+1)j_{k+1})} H_{\alpha}^{j_o} \tilde{p}_1^{j_1} 
\tilde{p}_2^{j_2}\dots \tilde{p}_{k+1}^{j_{k+1}}e^{iH_{\alpha}}\leqno(5.3.1a)$$
with 
$$j_o + j_1 + \dots j_{k+1} + \ell = N + k, \;\;\;\;\; j_1 + 2j_2 + \dots (k+1)j_{k+1} = \ell +
1.\leqno(5.3.1b)$$
 By the argument in the proof of Theorem C ( \S 4), taking the residue removes the factors of
$L D_s$ and replaces an expression $res (LD)^{-1} F(I) e^{i H_{\alpha}}$ by the value of
$F(D_{\alpha} + \half) T(\alpha)$ at a regular point.  It follows that
$$a_{\gamma k} = \sum C_{k;j,\ell} L^{-k + j_1 + \dots j_{k+1} }[ H_{\alpha}^{j_o} \tilde{p}_1^{j_1} 
\tilde{p}_2^{j_2}\dots \tilde{p}_{k+1}^{j_{k+1}}](D_{\alpha} + \half) T(\alpha) \leqno(5.3.2)$$
where the sum is taken over the indices specificed in (5.3.1b).

\medskip

\noindent{\bf \S 5.4 Jacobi degrees}
\medskip

We now show that the degree of $I_{\gamma k}$ as a polynomial in the Jacobi data
$y_{ij}, \dot{y}_{ij}$ is $\leq 6k+6.$  
The proof  a detailed review of the
construction of the polynomials $f_j(I_1,\dots,I_n)$ in Lemma (2.22) as well as the
relation between the Jacobi degrees of the polynomials $f_j$'s and of the operators
${\cal F}_{k, -1}$ in Theorem C.  For the sake of brevity, and since it is quite routine,
we will be a little sketchy in  a few of the details. We will use the notation 
$J-deg$ for the degree of a polynomial in the metric data above with respect to the
Jacobi field components.

\noindent{\bf (5.4.1) Lemma} {\it $J-deg \tilde{p}_{j+1} = J-deg\; p_{j+1} = J-deg
\;f_j= 6j+6.$}
\medskip

\noindent{\bf Proof} From the formulae (4.19) relating the $\tilde{p}_j (I_1,\dots,I_n)$'s and
the $p_j(I_1,\dots,I_n)$, and from the fact that $\frac{H_{\alpha}}{D_s}$ has Jacobi
degree 0, we see that 
$$J-deg \tilde{p}_j = \mbox{max}\;\;\{J-deg p_1,\dots,J-deg p_j\} \leqno(5.4.2a)$$
On the other hand, the relation between the $p_j$'s and $f_j$'s is given by comparing
(3.3) and (3.12):
$$f_j = \frac{2}{L}p_{j+1} + \sum _{i+r+2 =j+1} p_{i+1} p_{r+1}  \leqno(5.4.2b)$$
Let us assume $J-deg f_j = 6j+6$: an easy induction using (5.4.2b) then shows that $J-deg 
p_{j+1}= 6j+6$, and another using (5.4.2a) shows that $J-deg \tilde{p}_{j+1} = 6j+6.$

Hence we must prove (5.4.2b). By (2.39c) we have that 
$$J-deg f_j = J-deg {\cal D}_{-j}^{j+ \frac{1}{2}}. \leqno(5.4.3)$$
To determine the latter, we need to recall that  ${\cal D}_{-j}^{j+ \frac{1}{2}}$ is
the $h^j$-term in the expression
$$Ad(e^{i h^{j + \frac{1}{2}} Q_{j+ \frac{1}{2}}})Ad( e^{i h^j Q_j})\dots Ad( e^{i
h^{\frac{1}{2}} Q_{\frac{1}{2}}}){\cal D}_h $$
with $Ad(V)A := V^{-1}AV$.    Using that $Ad(e^{iB}) = e^{i ad(B)},$ and expanding everything
in a formal $h$-series, we find that ${\cal D}_{-j}^{j+ \frac{1}{2}}$ is the $h^j$-term 
in the series
$$\sum_{n, N_{\frac{1}{2}},\dots, N_{2j -1}} h^{\frac{n}{2} -2} 
(h^{j+ \frac{1}{2}})^{N_{j+ \frac{1}{2}}}(h^j)^{N_{j}}\dots (h^{ \frac{1}{2}})^{N_{ \frac{1}{2}}}
ad(Q_{j + \frac{1}{2}})^{N_{j+ \frac{1}{2}}} \dots ad(Q_{\frac{1}{2}})^{N_{\frac{1}{2}}}
{\cal D}_{2 - \frac{n}{2}}. \leqno (5.4.4)$$

We claim
\smallskip

\noindent{\bf (5.4.5) Claim}: {\it Suppose $J-deg Q_{\frac{m}{2}} = 3m$ for $m\leq 2j+1$. Then:
$J-deg {\cal D}_{-j}^{j+ \frac{1}{2}}= 6j + 6$.}
\smallskip

\noindent{\bf Proof}: The $h^{j}$ term in (5.4.4) is the sum of terms with indices satisfying  
$\frac{n}{2} - 2 + N_{j+ \frac{1}{2}}(j+ \frac{1}{2}) + \dots + N_{\frac{1}{2}}( \frac{1}{2})
=  j.$  Multiplying by 6 and using the hypothesis we find that 
$$J-deg {\cal D}_{-j}^{j+ \frac{1}{2}} = n + 6 \sum_{m=1}^{2j+1} \frac{m}{2} N_{\frac{m}{2}}
= n +6[ j+ 2 -\frac{n}{2}].$$
The claim now follows from the fact that  $n\geq 3$ in any term in the sum.\qed
\smallskip

\noindent{\bf (5.4.6) Claim}: {\it Suppose $J-deg {\cal D}^{\frac{j}{2}}_{-\frac{j}{2} +
\frac{1}{2}} = 3j+3.$  Then:  $J-deg Q_{\frac{j+1}{2}} = 3j+3$.}
\smallskip

\noindent{\bf Proof}:  Follows from (2.39c,d) which implies that
$$J-deg Q_{\frac{j}{2}+ \frac{1}{2}} = J-deg 
{\cal D}^{ \frac{j}{2}}_{-\frac{j}{2}+ \frac{1}{2}}. \leqno(5.4.7)$$\qed
\smallskip

We then show:
\smallskip

\noindent{\bf (5.4.8) Claim}: {\it $J-deg Q_{\frac{j}{2}} = 3j$.}
\smallskip

\noindent{\bf Proof}:    We  prove this by  induction on $j$.  For $j=1$ it holds
by explicit calculation:  from (2.27a), $deg Q_{\frac{1}{2}} = deg {\cal D}_{\frac{1}{2}}
= deg {\cal L}_{\frac{1}{2}} = 3.$  It follows then by Claim (5.4.5) (with j=0)
that $J-deg {\cal D}_o^{\frac{1}{2}} = 6$, and then by Claim (5.4.6) (with j=1)
that $J-deg Q_1 = 6.$  The rest of the induction proceeds similarly and the details
are left to the reader. \qed
\smallskip

The proof of Lemma (5.4.1) is completed by combining (5.4.3) and Claims (5..4.5) and (5.4.6).\qed
\medskip

Finally we have
\medskip

\noindent{\bf(5.4.9)  Lemma:} {\it $J-deg a_{k\gamma} \leq 6k + 6$.}
\medskip

\noindent{\bf Proof}:  As in \S 5.3, the kth residual 
terms in $[(D_s + H_{\alpha}) + {\cal P}_k]^k e^{i H_{\alpha}}e^{ iL {\cal P}_k}$
involve only the products
$$ \tilde{p}_1^{j_1} \tilde{p}_2^{j_2}\dots \tilde{p}_{k+1}^{j_{k+1}}$$
of the $\tilde{p}_j$'s with 
$$j_o + j_1 + \dots j_{k+1} + \ell = N + k, \;\;\;\;\; j_1 + 2j_2 + \dots (k+1)j_{k+1} = \ell + 1,
\;\;\;\;\;\;\;\;\ell \leq k, N \leq k+1.\leqno(5.4.10)$$
Using that $J-deg \tilde{p}_{m} \leq 6m$,  the Jacobi degree of any such term cannot
exceed $6(j_1 + \dots + j_{k+1}).$  We claim that the maximum possible value, subject to
the constraints (5.4.10), is 6k+6.  Indeed, since $j_1 + \dots + j_{k+1} = \ell + 1 -
(j_2 + 2j_3 + \dots kj_{k+1}) \leq k + 1 - (j_2 + 2j_3 + \dots kj_{k+1})$ the maximum
value is achieved when $j_1 = k+1$ and all other $j_m$'s are zero.\qed  
\medskip

This completes the proof of Theorem A. \qed

\section {Quantum Birkhoff normal form coefficients}

In this section we give a brief summary of the algorithm for calculating the quantum Birkhoff normal
form coefficients $B_{\gamma;k,j}$.  From these coefficients one could determine the wave invariants
$a_{\gamma;k}$ as described in \S 5.3, but since the $B_{\gamma;k,j}$'s are  simpler spectral invariants we choose
to concentrate on them.  In the next section, we apply the algorithm to the calculation of the coefficients
$B_{\gamma;o,j}$ in dimension 2.
\medskip

\noindent{\bf Preliminaries}
\medskip

\noindent{\bf (6.1) Definition }{\it
The quantum Birkhoff normal form (QBNF) coefficients are the coefficients of the
monomials $I^{j} := I_1^{j_1}\dots I_n^{j_n} = |z_1|^{2j_1}\dots |z_n|^{2j_n}$ in the complete Weyl
symbol $\tilde{p}_k(|z|^1_1,\dots,|z|^2_n)$ of  the coefficient operator
$\tilde{p}_k(I_1,\dots,I_n)$ of the model normal form of Theorem B:
$$\tilde{p}_k(|z|^1_1,\dots,|z|^2_n) = \sum_{j\in \Nb^n: |j| \leq k+1} B_{\gamma; k,j} |z|^{2j}.$$}
\medskip

There is of course something arbitrary in the emphasis on the $\tilde{p}_k(I_1,\dots,I_n)$'s here.
The coefficients of the monomials of the complete symbol of the $p_k(I_1,\dots,I_n)$'s
 are equally spectral invariants and we
only prefer the $\tilde{p}_k(I_1,\dots,I_n)$'s  to maintain consistency with the
terminology of [G.1-2]. Moreover, the coefficients of the operator monomials $I^{j}$
in either the $p_k$'s or $\tilde{p}$'s are also spectral invariants and in view of
the relation between wave invariants and the QBNF
 (\S 4), it is the operator coefficients which are most simply related to the wave invariants. 

The crucial point for the effective calculation of the QBNF coefficients is that they can be 
obtained from the coefficients of the complete symbol
of the semi-classical normal form (SCNF) 
$$W_h^*{\cal R}_h W_h|_o \sim h^{-2} 
+ \sum_{j=0}^{\infty} h^{k} f_k (I_1, \cdots,I_n)$$
of Lemma (2.22) (after restriction to weight 0.)
  Indeed, 
 after substituting $h = {\cal R}^{-1},$ the square of the QBNF is the SCNF.
 Hence, the key invariants are really the coefficients of the monomials in the complete symbols
$$f_k(|z_1|^2, \dots,|z_n|^2) = \sum_{j \in \Nb^n, |j|\leq k+2} c_{\gamma; kj} |z|^{2j}$$ of the
coefficient operators $f_k(I_1,\dots,I_n)$. In the remainder of the section, we drop the subscript
$\gamma$ from the notation.

On the operator level, the relation between the $\tilde{p}_k$'s, $p_k$'s and $f_k's$ is very simple:
the operators $p_k(I_1,\dots,I_n)$ of Theorem B are related to the
operators $f_k(I_1,\dots, I_n)$ of Lemma (2.22) by:
$$f_k = \frac{2}{L} p_{k+1} + \sum_{i + j = k, i,j\geq 1} p_i p_j.\leqno(6.2a)$$
while the operators $\tilde{p}_j$ are related to  the operators $p_k$ by:
$$\tilde{p}_{\ell} = \sum_{k\geq 1, (j_1,\dots,j_k): k + |j|=\ell} C_{k,(j_1,\dots,
j_k)}  H_{\alpha}^{|j|} p_k$$
for certain universal (multinomial) constants $C_{k,(j_1,\dots, j_k)}$.  On the symbol level, the
relation is a little more complicated since one has to compose the symbols in the Weyl calculus.  In the
following we summarise the steps involved in calculating the coefficients of the monomials $|z|^{2j}$
in the complete symbols of the $f_k$'s. To distinguish on operator from its complete symbol we use the
notation $\hat{f}(I_1,\dots,I_n)$ for the Weyl operator with complete symbol $f(|z_1|^2,\dots,|z_n|^2).$
We also recall that the notation $|_o$ refers to the weight
  0 part of an operator relative to $D_s$: $(A_2(s, x, D_x) D_s^2 + A_1(s, x, D_x) D_s + A_o(s,x,D_x))|_o
=A_o(s,y,D_y).$  We will assume that the ${\cal L}$'s and ${\cal D}$'s have been expressed in terms of
the weightless Fermi normal coordinats $(s,x)$ of \S 1.5. The complex  coordinates $z$ are given by
$z=x + i\xi.$
\medskip

\noindent{\bf Summary of the algorithm}
\bigskip

\begin{tabular}{r|l} Step & Relevant data \\ \hline \hline
Step 1 & Express $  \hat{f}_{\frac{N-1}{2}}$ in terms of ${\cal D}_{2 - \frac{n}{2}}$'s \\ \hline 
Formula for $  \hat{f}_{\frac{N-1}{2}}$  & $  \hat{f}_{\frac{N-1}{2}}(I_1,...,I_n):=\frac{1}{L}\int_o^L\int_{T^n}
V_t^*{\cal D}^{\frac{N}{2}}_{-\frac{N-1}{2}}|_o V_tdtds$ (diagonal part; N odd)\\ \hline
Equivalently & 
$\hat{f}_{\frac{N-1}{2}}(q_1 + \half,...,q_n + \half):=\frac{1}{L}\int_o^L 
\langle  {\cal D}^{\frac{N}{2}}_{-\frac{N-1}{2}} \gamma_q, \gamma_q \rangle ds $\\ \hline
Weyl symbol & $ f_{\frac{N-1}{2}}(|z_1|^2,\dots, |z_n|^2) = \frac{1}{L} \sum_{|k|\leq
\frac{N+3}{2}} c_{\frac{N-1}{2};j} |z|^{2j}$\\ \hline 
Recursion for ${\cal D}^{\frac{N+1}{2}}_{-\frac{N-1}{2}}$  &
 ${\cal D}^{\frac{N+1}{2}}_{-\frac{N-1}{2}} =$ term of
order $h^{\frac{N-1}{2}}$ in \\
$\cdot$ & $ \sum_{n=3}^{N+2}h^{-2+\frac{n}{2}}{\cal D}^{\frac{N}{2}}_{2-\frac{n}{2}}
h^{\frac{N-1}{2}}\{ [L^{-1} D_s, \tilde{Q}_{\frac{N+1}{2}}]+ {\cal D}^{\frac{N}{2}}_{-\frac{N-1}{2}}\}$\\
$\cdot$ & $+ \sum_{n=N+4}^{\infty} h^{-2 + \frac{n}{2}}{\cal D}^{\frac{N}{2}}_{2-\frac{n}{2}}
+\sum_{m=2}^{\infty}h^{m\frac{N+1}{2}}\frac{i^m}{m!}(ad \tilde{Q}_{\frac{N+1}{2}})^m (L^{-1} D_s)$\\
$\cdot$ & $+\sum_{m=1}^{\infty}\sum_{p=3}^{\infty} h^{m\frac{N+1}{2}+\frac{p}{2}-2}\frac{i^m}{m!}
(ad \tilde{Q}_{\frac{N+1}{2}})^m 
({\cal D}^{\frac{N}{2}}_{2 -\frac{p}{2}})$\\ \hline
Operator Transport  : &
 $N+1$  odd: $\{ [L^{-1} D_s, \tilde{Q}_{\frac{N+1}{2}}]+ 
{\cal D}^{\frac{N}{2}}_{-\frac{N-1}{2}}\}|_o = 0$\\
equations for $\tilde{Q}$'s &
 $N+1$  even: 
$\{ [L^{-1} D_s, \tilde{Q}_{\frac{N+1}{2}}]+ {\cal D}^{\frac{N}{2}}_{-\frac{N-1}{2}}\}|_o =
f_{\frac{N-1}{2}} (I_1,...,I_n). $ \\ \hline
Symbol  & 
 $N+1$  odd: $ L^{-1} \partial_s \tilde{Q}_{\frac{N+1}{2}}(s,z,\bar z) = -i 
{\cal D}^{\frac{N}{2}}_{-\frac{N-1}{2}}|_o (s,z \bar z)$\\
Transport equations &
 $N+1$  even: 
$ L^{-1} \partial_s \tilde{Q}_{\frac{N+1}{2}}(s, z, \bar z)=-i ({\cal
D}^{\frac{N}{2}}_{-\frac{N-1}{2}}\}|_o(s,z,\bar z)
 + i f_{\frac{N-1}{2}} (|z_1|^2, \dots, |z_1|^2).$\\ \hline 
Homological  &
 $N+1$  odd: $  \tilde{Q}_{\frac{N+1}{2}}(0,e^{i\alpha}z, e^{-i\alpha}\bar  z) - 
\tilde{Q}_{\frac{N+1}{2}}(0,z,\bar z) = -i L \int_0^L {\cal
D}^{\frac{N}{2}}_{-\frac{N-1}{2}}|_o (s,z \bar z)ds$\\
equations &
 $N+1$  even: 
$ \tilde{Q}_{\frac{N+1}{2}}(0,e^{i\alpha}z, e^{-i\alpha}\bar  z)
-\tilde{Q}_{\frac{N+1}{2}}(0,z,\bar z) = -i L \int_0^L ({\cal
D}^{\frac{N}{2}}_{-\frac{N-1}{2}}\}|_o(s,z,\bar z) ds$ \\
$\cdot$ &
 $- f_{\frac{N-1}{2}} (|z_1|^2, \dots, |z_1|^2).$\\ \hline
 Solutions:& N+1 odd:
$\tilde{Q}_{\frac{N+1}{2}}(s, z, \bar z) = \tilde{Q}_{\frac{N+1}{2}}(0, z, \bar z) 
 -iL \int_o^s{\cal D}^{\frac{N}{2}}_{-\frac{N-1}{2}}|_o(t,z,\bar z) ds$\\ $\cdot$ &
N+1 even:
$\tilde{Q}_{\frac{N+1}{2}}(s,z,\bar z)= \tilde{Q}_{\frac{N+1}{2}}(0,z,\bar z)$\\
$\cdot$ &
 $-i L \int_o^s\{{\cal D}^{\frac{N}{2}}_{-\frac{N-1}{2}}|_o(t,z,\bar z) - f_{\frac{N-1}{2}}
(|z_1|^2,...,|z_n|^2)\}|_ods$\\ \hline  \hline
Step 2 & Express ${\cal D}_{2 - \frac{n}{2}}$'s in terms of metric data \\ \hline
Conjugate to ${\cal L}_{2 - \frac{n}{2}}$'s & ${\cal D}_{2 - \frac{n}{2}} = 
\mu ({\cal A_L}^*){\cal L}_{2 - \frac{n}{2}}\mu({\cal A_L}^*)^{-1}$ \\ \hline
${\cal L}_{2 - \frac{n}{2}}$ = &
$-(hL)^{-2} g^{oo}_{[h]} + 2i(hL)^{-1}g^{oo}_{[h]}\partial_s + i(hL)^{-1}\Gamma^o_{[h]}$ \\
$h^{\frac{n}{2}- 2}$-term of & $+
 h^{-1}( \sum_{ij=1}^n g^{ij}_{[h]}\partial_{x_i}\partial_{x_j}) + h^{-\half}(\sum_{i=1}^{n} \Gamma^{i}_{[h]}
\partial_{x_i}) + (\sigma)_{[h]}$ \\ \hline \hline
Step 3 & Solve for coefficients \\ \hline
Coeff's of ${\cal D}^*_{2 - \frac{k}{2}}|_o$ &   ${\cal D}^*_{2 - \frac{n}{2}}|_o(s,z,\bar z) =
=\sum_{|m|+|n|\leq 2k} d^*_{\frac{k}{2}-2; m,n} (s) z^m \bar z^n$ \\ \hline
Coeff's of ${\cal D}^{\frac{N}{2}}_{-\frac{N-1}{2}}|_o$ & $d^{\frac{N}{2}}_{\frac{N-1}{2}; m,n}(s) =
F^{\frac{N}{2}}_{\frac{N-1}{2}; m,n} (d_{\frac{k}{2}-2; m,n};k\leq \frac{N-1}{2})$ \\ \hline
Coeff's of $\tilde{Q}_{\frac{N+1}{2}}$'s  & $ \tilde{Q}_{\frac{N+1}{2}}(s,z,\bar z) = 
\sum_{|m|+|n|\leq N+3} q_{\frac{N+1}{2}; mn}(s) z^m \bar z^n$ \\ \hline 
Off-diagonal Coefficients & $q_{\frac{N+1}{2}; m,n}(0) = -i (1 - e^{i (m-n)\alpha})^{-1}
-i L \int_o^L d_{\frac{N-1}{2};m,n}(s) ds$\\ \hline
Diagonal Coefficients & $c_{\frac{N-1}{2}; k} =  \bar d^{\frac{N}{2}}_{\frac{N-1}{2}; k,k} =
\frac{1}{L}\int_o^L d^{\frac{N}{2}}_{\frac{N-1}{2}; k,k}(s) ds$ \\ \hline
\end{tabular}
\bigskip

\noindent{\bf Remark}  Since $\sqrt{\Delta}$ commutes with complex conjugation, the odd order terms in
its Weyl symbol are even functions under the involution $(s,\sigma, x,\xi) \rightarrow (s, -\sigma,x,
-\xi).$  Although we have restricted to the positive cone where $\sigma > 0$, one gets similar normal
forms for $\sigma < 0$ as long as $D_s$ is interpreted as $|D_s|$.  Since the Harmonic oscillators and
their powers have even symbols, it must be the case that  the QBNF coefficients of the even order terms  vanish. 
For instance the coefficient $B_{o;2}$ of $\frac{I}{D_s}$ must vanish automatically.

\section{Explicit formulae  in dimension 2}

To illustrate the algorithm, we carry out the
details of the calculation of $a_{\gamma o}$, or more importantly the normal form coefficients $B_{k;j}$ for k=0,
in dimension 2. The result may be summarized as follows:
\medskip

\noindent{\bf QBNF coefficients for k=0, dim =2}{\it The complete symbol of $f_o(I)$ in complex
coordinates $z = y + i\eta$ has the form $B_{o;4}|z|^4 + B_{o;0}$ where $B_{o; j}$
are given for both $j=0,4$ by weight -2  Fermi-Jacobi polynomials of the form: 
$$B_{o;j} =   \frac{1}{L}\int_o^L [a |\dot{Y}|^4 +   b_1 \tau |\dot{Y}Y|^2 +
b_2 \tau  Re (\bar{Y} \dot{Y})^2
+c \tau^2 |Y|^4 + d \tau_{\nu \nu} |Y|^4 +e \delta_{jo} \tau ] ds + $$
$$ +\frac{1}{L}\sum_{0\leq m,n \leq 3; m+n=3} C_{1;mn}\frac{sin((n-m)\alpha)}
{|(1 - e^{i(m-n)\alpha})|^2}
|\int_o^L  \tau_{\nu}(s)\bar{Y}^mY^n](s)ds|^2 $$
$$+\frac{1}{L}\sum_{0\leq m,n \leq 3; m+n=3} C_{2;mn}Im \{\int_o^L \tau_{\nu}(s)\bar
{Y}^mY^n(s)[\int_o^s\tau_{\nu}(t)\bar {Y}^nY^m](t)dt] ds \}.$$ }
\medskip

This corroborates the previous remark that the term $B_{o k}$ vanishes.  As will be seen, the
corresponding density is a total $\partial_s$-derivative. Also, the coefficient $\delta_{jo}$ of $\tau$ is the
Kronecker symbol, i.e. equals one if $j=0$ and vanishes if $j=4$.

The expressions for the normal form coefficients in higher dimensions are very similar, and it is only
for the sake of simplicity that we have  stated the result in  dimension 2. 

The result for the wave coefficient $a_{\gamma o}$ is then very similar, modulo the Floquet
factors.  Indeed,
 by (4.30) or by \S 5.3 the wave invariant is related to the normal form coefficients by
$$a_{\gamma o} = {\cal F}_{o, -1}(D_{\alpha}, L) T(\alpha)$$
with
$${\cal F}_{o, -1}(D_{\alpha}, L) = L \tilde{p}_1(D_{\alpha_1}+ \half,\dots,D_{\alpha_n}+\half,L)
= L p_1(D_{\alpha_1} + \half ,\dots,D_{\alpha_n}+ \half).$$
Hence 
$$\int_{\gamma} I_{\gamma o 2}ds = L p_1(D_{\alpha_1}+ \half,\dots, D_{\alpha_n}+\half)
 T(\alpha) = L p_1(D_{\alpha_1}+\half,\dots, D_{\alpha_n}+\half) \Pi_{j=1}^n 
(1- e^{i\alpha_j})^{-1}$$  
In
view of (2.33) and the fact that $p_1 = \half L f_o$ we have
$$a_{\gamma o} = \int_{\gamma} I_{\gamma o}ds =
\half L^2 f_o(D_{\alpha_1}+\half,\dots, D_{\alpha_n}+\half)
\Pi_{j=1}^n  (1- e^{i\alpha_j})^{-1}. \leqno(7.1)$$
From (7.1) and the formula for the complete symbol of $f_o$ one gets the explicit formula for
$a_{\gamma o}$ as stated in the Introduction.

To prove that the $B_{o;j}$ coefficients have the form claimed above, we begin with the expression for $f_o$
from (2.33) or from the table in \S 6:
$$f_o(I_1,\dots, I_n) = \frac{1}{L}\int_o^L\int_{T^n} V_t^*{\cal D}_o^{\half} V_t|_odt.\leqno(7.2) $$
  We also have, from (\S 2, (2.36) and below), that 
$${\cal D}^{\half}_o = {\cal D}_o + \half i [{\cal D}_{\half}, \tilde{Q}_{\half}]. \leqno (7.3)$$
We wish to evalute the coefficients of the complete symbol of $f_o$  in terms of integrals over
$\gamma$ of Fermi-Jacobi data.  At first, we will allow the dimension to be arbitrary; when it is time
to substitute in  metric expressions we will restrict to dimension 2.

Let us consider first the
second term on the right side, which simplifies the expression
$$[{\cal D}_{\half},\tilde{Q}_{\half}] + \half i^2 [[L^{-1} D_s,\tilde{Q}_{\half}], \tilde{Q}_{\half}]$$
by using $[L^{-1} D_s,\tilde{Q}_{\half}] = i {\cal D}_{\half}.$ We recall from \S 2 or from the table of \S6 that 
the complete symbol of $\tilde{Q}_{\half}$ is given (in complex vector notation) by
$$\tilde{Q}_{\half}(s,z,\bar{z}) = \tilde{Q}_{\half}(0,z,\bar{z}) + L \int_o^s {\cal
D}_{\half}|_o(t,z,\bar{z})dt,\;\;\;\;\;\;
\tilde{Q}_{\half}(0,e^{i\alpha }z,e^{-i\alpha}\bar{z}) - \tilde{Q}_{\half}(0,z,\bar{z}) = L \int_o^L {\cal
D}_{\half}|_o(t,z,\bar{z})dt.$$
We note that ${\cal D}_{\half}$ is independent of $D_s$ so that ${\cal D}_{\half}|_o = {\cal D}_{\half}.$
It follows that 
$$ [{\cal D}_{\half}(s), \tilde{Q}_{\half}(s)] =  [{\cal D}_{\half}(s), \tilde{Q}_{\half}(0)] +  [{\cal
D}_{\half}(s),\int_o^s {\cal D}_{\half}(t)dt]$$
so that the second term of (7.3) contributes to $f_o(I_1,\dots,I_n)$ the {\em diagonal part} of
$$\frac{i}{2L^2}\{[\tilde{Q}_{\half}(0,e^{i\alpha}z,e^{-i\alpha}\bar{z}), \tilde{Q}_{\half}(0,z,\bar{z}] +
\frac{1}{2}\{[\int_o^L {\cal D}_{\half}(s)ds, \int_o^s{\cal D}_{\half}(t)dt].\leqno(7.4 diag)$$ 
Here, the bracket $[,]$ denotes the  commutator of complete symbols in the sense of operator
(or complete symbol) composition.

To determine the diagonal part, we write (as usual)
$$\tilde{Q}_{\half}(0,z,\bar{z}) = \sum_{mn: |m|+|n|=3} q_{\half; mn}(0) z^m\bar{z}^n,\;\;\;\;\;\;\;
{\cal D}_{\half}(s,z,\bar{z}) = \sum_{mn: |m|+|n|=3} d_{\half; mn}(s) z^m\bar{z}^n.$$

Let us denote the operator composition of two complete Weyl symbols $a,b$ by $a \# b$, that is
$$Op^w(a) \circ Op^w (b) = Op^w (a \# b)$$
and recall (cf. [Ho III, Theorem 18.5.4]) that this composition has an asymptotic expansion 
$$a \# b \sim ab + i P_1(a,b) + \frac{i^2}{2!} P_2(a,b) + \dots$$
where $P_k(a,b)$ is the higher order Poisson bracket (or transvectant), a bidifferential operator
 given in complex notation by
$$P_k(a,b) = (\partial_z\partial_{\bar{w}} - \partial_{\bar{z}}\partial_w)^k a(z,\bar{z})b(w,\bar{w})|_{z=w}.$$
It is well known [loc.cit.] that the commutator is given symbolically by the odd expansion
$$a\#b - b\#a \sim \frac{1}{i} P_1(a,b) + \frac{1}{i^3 3!}P_3(a,b) +\dots$$
while anticommutators involve only the even transvectants.  Since the complete symbols in (7.4 diag) are
homogeneous polynomials of degree 3, the commutators involve only $P_1$ and $P_3$. One easily computes that
$$P_1(z^m\bar{z}^n, z^{\mu}\bar{z}^{\nu}) = C_{1;mn\mu\nu} z^{m+\mu -1}\bar{z}^{n + \nu -1}$$
where $C_{1;mn\mu\nu} = \half\sigma((m,n),(\mu,\nu))$ with $\sigma$ the standard symplectic inner product, and
that
$$P_3(z^m\bar{z}^n, z^{\mu}\bar{z}^{\nu}) = C_{(m,n),(\mu,\nu)} z^{m+\mu - 3}\bar{z}^{n + \nu -3}$$
for certain other coefficients $C_{(m,n),(\mu,\nu)}$.  A straightforward computation then shows that the
diagonal part in (7.4 diag) is the sum of the following terms:
$$\frac{i}{2L^2}\{[\tilde{Q}_{\half}(0,e^{i\alpha}z,e^{-i\alpha}\bar{z}), 
\tilde{Q}_{\half}(0,z,\bar{z}]\leqno(7.5.1)$$
$$ =
L^{-2}[C_{1;3030} q_{\half; 30}(0)q_{\half; 03}(0) e^{i3\alpha} + C_{1;2112}q_{\half;21}(0)q_{\half;12}(0)e^{i
2\alpha} +C_{1;1221} q_{\half; 12}(0)q_{\half; 21}(0)e^{i \alpha} + C_{1;0330}q_{\half; 03}(0)q_{\half;
30}(0)]|z|^4 +$$
$$+[C_{3;3030} q_{\half; 30}(0)q_{\half; 03}(0) e^{i3\alpha} + C_{3;2112}q_{\half; 21}(0)q_{\half; 12}(0)e^{i
2\alpha} +C_{3;1221} q_{\half; 12}(0)q_{\half; 21}(0)e^{i \alpha} + C_{3;0330}q_{\half;03}(0)q_{\half;30}(0)] $$
plus
$$\frac{1}{2}\int_o^L\int_o^s[\tilde{D}_{\half}(s,z,\bar{z}),
\tilde{D}_{\half}(t,z,\bar{z}]dsdt \leqno(7.5.2)$$
$$ =\int_o^L\int_o^s [C_{1;3030}
d_{\half;30}(s)d_{\half; 03}(t) + C_{1;2112}d_{\half;21}(s)d_{\half; 12}(t)
 +C_{1;1221}
 d_{\half;12}(s)d_{\half;21}(t) +
C_{1;0330}d_{\half;30}(s)d_{\half; 03}(t) dsdt] |z|^4 +$$
$$\int_o^L\int_o^s [C_{3;3030} d_{\half; 30}(s)d_{\half; 03}(t)  +
C_{3;2112}d_{\half; 21}(s)d_{\half; 12}(t) +C_{3;1221}
d_{\half;12}(s)d_{\half;21}(t) + C_{3;0330}d_{03}(s)d_{\half;30}(t) ]dsdt. $$ We
observe that there is no term of order $|z|^2$.  Using   
 that
$$q_{\half; mn} = 
-i L (e^{i\alpha \cdot(m-n)} -1)^{-1}  \int_o^L d_{\half; mn}(s)ds $$
 we reduce our problem to the calculation of the 
$d_{\half; mn}(t)$'s and the $d_{o;mn}(t)$'s.

To evaluate  the expressions ${\cal D}_{\half}$ and ${\cal D}_o$.
we conjugate back to the ${\cal L}$'s:
$${\cal D}_o = \mu({\cal A}_L^*) {\cal L}_o\mu(({\cal A}_L^*){-1},\;\;\;\;\;\;\;\;
{\cal D}_{\half} = \mu({\cal A}_L^*){\cal L}_{\half}\mu({\cal A}^{*}_L)^{-1} (7.6a)$$
with
$${\cal A}_L:=\left( \begin{array}{ll} L^{\half} Im\dot{Y} &  L^{\half} Re \dot{Y} \\  L^{-\half}Im Y & 
L^{-\half}Re Y
\end{array}\right).\leqno(7.6b)$$
 Using (2.11) we then compute the ${\cal L}$'s in (unscaled) Fermi normal coordinates as follows:
$${\cal L}_{\half}= -\sum_{\beta: |\beta|=3}
\frac{1}{3!} \partial_y^{\beta} g^{oo}(s,0) y^{\beta} \leqno(7.7.1/2)$$
$${\cal L}_o =  -\sum_{\beta: |\beta|=4}\frac{1}{4!} \partial_y^{\beta} g^{oo}(s,0) y^{\beta} 
-\sum_{\beta: |\beta|=2} \partial_x^{\beta}[g^{oo} D_s +
J^{-\half}D_s(g^{oo}J^{\half})]|_{y=0}y^{\beta} \leqno(7.7.1)$$
$$+\sum_{\beta: |\beta|=2}\sum_{ij=1}^n \partial_y^{\beta}g^{ij}(s,0) y^{\beta}\partial_{y_i}
\partial_{y_j} + \sum_{i,j=1}^n \partial_{y_j} \Gamma^i(s,0) y_j \partial_{y_i}$$
$$+ \partial_{s}^2 +C_n \tau(s,0)$$
where $\tau$ is the scalar curvature. The last term comes from $\Delta_{\half}1$.

We then change variables $y=L x$ and conjugate the symbols.  By metaplectic covariance of the Weyl
calculus,
the conjugations change the complete Weyl symbols of the ${\cal L}$'s (in the $x$ variables) by the linear 
symplectic transformation ${\cal A}_L$, i.e by the substitutions
$$\begin{array}{l} x \rightarrow L^{-\half} [(Re Y) x + (Im Y)\xi] = \half L^{-\half}[\bar {Y}\cdot z + Y \cdot
\bar{z}]
\\
\xi \rightarrow L^{\half}[Re \dot{Y} x + (Im \dot{Y})\xi] = \half L^{\half} [\bar{\dot{Y}}\cdot z +
\dot{Y}\bar{z}]
\end{array}.
\leqno (7.8)$$

It is evident that the normal form coefficients are going to be rather lengthy.
To give the flavor of the full calculations in the simplest setting, we now specialize to the case of
surfaces.  Later we will briefly extend the calculations to all dimensions.

\noindent{\bf Dimension 2}

  In dimension 2 we have (in scaled Fermi coordinates) 
$$g^{oo}(s,y) =1+ C_1 \tau(s) y^2 +C_2 \tau_{\nu}(s)y^3 + \dots\;\;\;\;\;\;\;g^{11}=1 \;\;\;\;
J(s,u) = \sqrt{g_{oo}}= 1 + C_1'\tau(s)y^2 + \dots $$, 
 for some constants $C_j,C_j'$ which will change from line to line.
Hence the 1/2-density Laplacian in Fermi normal coordinates equals
$$-\Delta = J^{-1/2}\partial_s g^{oo}J \partial_s J^{-1/2}
+ J^{-1/2}\partial_{y}  J \partial_{ y}
J^{-1/2}$$
$$\equiv g^{oo}\partial_s^2 + \Gamma^o \partial_s + 
  \partial_{y}^2+  \Gamma^{1}\partial_{y} + \sigma_o$$
and the rescaled Laplacian equals 
$$-\Delta_{h} =  -(Lh)^{-2} g^{oo}_{[h]} + 2i(Lh)^{-1} g^{oo}_{[h]}\partial_s +   i(Lh)^{-1} \Gamma^o_{[h]}-
g^{oo}_{[h]}\partial_s^2 + \Gamma^o_{[h]}\partial_s + \partial_y^2
+ h^{-\half} \Gamma^{1}_{[h]}\partial_{y} + (\sigma)_{[h]}$$
  Using the Taylor expansion of the metric coefficients one finds that the ${\cal L}$'s have the form
$${\cal L}_{\half} = C L^{-2}\tau_{\nu}(s) y^3,\;\;\;\;\;{\cal L}_o = C_1 L^{-2}y^4 \tau_{\nu \nu} 
+ C_2 L^{-1} y^2 \tau \partial_s + C_3 L^{-1} \tau_s y^2 - \partial_s^2+  C_4\tau y \partial_y  + C_5\tau.
\leqno(7.9)$$ All terms have weight -2.

We now switch to the weightless coordinates $y = L x$ and get:
$${\cal L}_{\half} = C L \tau_{\nu}(s) x^3,\;\;\;\;\;{\cal L}_o = C_1 L^{2}x^4 \tau_{\nu \nu} 
+ C_2 L^{1} x^2 \tau \partial_s + C_3 L \tau_s x^2 - \partial_s^2+  C_4\tau x \partial_x  + C_5\tau.
\leqno(7.9)$$

Making the linear symplectic substitutions above we first get
$${\cal D}_{\half}(s,z,\bar{z}) = C L^{-\half} \tau_{\nu}(s)
 ( [ \bar {Y}\cdot z + Y \cdot \bar{z}])^{3}$$
hence 
$$d_{\half; mn}(s) =  C_{mn;3}L^{-\half}\tau_{\nu} [\bar
{Y}^m\cdot Y^n](s)\leqno(7.10.1)$$
$$\bar{d}_{\half; mn}(s) =  C_{mn;3}L^{-\frac{3}{2}}\int_o^{L}\tau_{\nu} [\bar
{Y}^m\cdot Y^n](s)ds \leqno(7.10.2)$$
$$q_{\half; mn} =  -i C_{mn;3}(1 -
e^{i(m-n)\alpha})^{-1} L^{\half}\int_o^L\tau_{\nu} [\bar {Y}^m\cdot
Y^n](s) ds \leqno(7.10.3)$$
 with $m+n=3$ and
for certain coefficients $C_{mn;\beta}$. It is evident that $d_{mn}(s)$ is a Fermi-Jacobi polynomial
of weight -2 and of Jacobi degree 3, so that the terms (7.5.1-2)
 are of Jacobi degree 6  as stated in Theorem A. The
diagonal part has terms of degree $|z|^4$ and $|z|^o$ with coefficients of the form
  (with $m,n=0,...3; m+n=3$):
$$\frac{1}{L}\{[\int_o^L  \tau_{\nu}(s)\bar{Y}^mY^n](s)ds] [\int_o^L \tau_{\nu}(t)\bar
{Y}^nY^m(t)dt]\frac{e^{i(n-m)\alpha}}{|1-e^{i(n-m)\alpha}|^2}$$
$$ -
[\int_o^L
\tau_{\nu}(s)\bar {Y}^nY^m(s)ds][\int_o^L\tau_{\nu}(t)\bar {Y}^mY^n](t)dt]
\frac{e^{-i(m-n)}}{|1-e^{-i(m-n)}|^2}\}
\leqno(7.11.1)$$ and
$$\frac{1}{L}\{\int_o^L \tau_{\nu}(s)\bar {Y}^mY^n(s)[\int_o^s\tau_{\nu}(t)\bar
{Y}^nY^m](t)dt] ds - \int_o^L \tau_{\nu}(s)\bar {Y}^nY^m(s)[\int_o^s\tau_{\nu}(t)\bar
{Y}^mY^n(t)dt] ds\}. \leqno(7.11.2)$$

To calculate the complete symbol of $D_s$-weight 0,  ${\cal D}_o|_o$, of the second term ${\cal D}_o^{\half}$ we 
make the same linear substitution and eliminate any $D_s$ appearing all the way to the right.
 We also invert the relation 
$$\mu({\cal A}_L^*)^{-1} D_s \mu({\cal A}_L^*) = D_s - \half(L^{-1} \partial_x^2 + L \tau x^2)$$
to get
$$\mu({\cal A}^*) D_s \mu({\cal A}^*)^{-1} = D_s - \half\mu({\cal A}^*)(L^{-1}\partial_x^2 + L \tau x^2)
\mu({\cal A}^*)^{-1}$$
and transform the complete symbol of quadratic term by the symplectic substitution.  The result is
that ${\cal D}_o|_o (s,z,\bar{z})$ equals
$$ C_1 \tau_{\nu \nu}  [\bar {Y} z + Y  \bar{z}]^4\leqno(7.12.1)$$
$$+C_2 L^{-2} \tau [\bar {Y} z + Y  \bar{z}]^2 \# 
( [\bar{\dot{Y}} z + \dot{Y}\bar{z}]^2 + \tau L^{-2} [\bar {Y} z + Y \bar{z}]^2))\leqno(7.12.2)$$
$$ + C_3 \tau_s [\bar {Y} z + Y  \bar{z}]^2 \leqno(7.12.3)$$
 $$ -2 L^{-2}\partial_s( [\bar{\dot{Y}} z + \dot{Y}\bar{z}]^2 -L^{-2}\tau [\bar {Y} z + Y
\bar{z}]^2)\leqno(7.12.4.1)$$
$$ +L^{-2}\{ [\bar{\dot{Y}} z + \dot{Y}\bar{z}]^2 - L^{-2}\tau  [\bar {Y}\cdot z + Y \bar{z}]^2\}
\#\{ [\bar{\dot{Y}} z + \dot{Y}\bar{z}]^2 -L^{-2} \tau [\bar {Y} z + Y  \bar{z}]^2\} \leqno(7.12.4.2)$$
$$+ C_4 L^{-1}\tau (\bar {Y} z + Y  \bar{z})\#(\bar{\dot{Y}} z + \dot{Y}\bar{z}) +C_5 \tau.\leqno(7.12.5)$$

Our concern is with the diagonal part of the complete symbol, that is, with the terms involving
$|z|^4, |z|^2, |z|^o$, and more precisely with their integrals over $\gamma$.  To begin with, we observe that
the diagonal part of term (7.12.1) is purely of degree $|z|^4$ and its average over $\gamma$ equals 
$$ (Const.) |z|^4 \cdot \frac{1}{L} \int_o^L \tau_{\nu \nu} |Y|^4 ds. \leqno(7.13.1)$$
The diagonal part of term (17.12.2) contributes only the $P_o$ and $P_2$ terms, of degrees $|z|^4$ and $|z|^o$
respectively, whose averages over $\gamma$  have the form
$$(|z|^4 /or/ |z|^o) \cdot 
 \frac{1}{L}\int_o^L \tau [a\tau |Y|^4 + b Re (\bar{Y}\dot{Y})^2 + c|Y\dot{Y}|^2 ]ds
\leqno(7.13.2)$$ with constants $a,b,c,d$ which can differ between the two degrees.
The missing $P_1$-term vanishes: it is a multiple of the Poisson bracket
$$P_1 ([\bar {Y} z + Y  \bar{z}]^2, \tau [\bar {Y} z + Y \bar{z}]^2))$$
which simplifies to a term of the form
$$\tau [\bar{Y}^2\dot{Y}^2 - Y^2\dot{\bar {Y}}^2] = \tau (\bar{Y}\dot{Y} - Y \bar{\dot{Y}})(\bar{Y}\dot{Y} +
Y\bar{\dot{Y}}) = C \tau \frac{d}{ds} |Y|^2$$
by the symplectic normalization of the Jacobi eigenfield.    However the integral
over $\gamma$ of this term vanishes, that is
$$\;\;\;\;\;\;\;\;\int_o^L \tau_s |Y|^2=0. \leqno(7.13.3)$$
This can be seen from the  Jacobi equation, which implies: 
$$ [\bar{Y} (Y')'' + \tau_s|Y|^2 + \tau Y'\bar{Y}]=0;$$
integrating over $\gamma$ and integrating the first term by parts twice kills the outer terms and
hence the inner one. 
In a similar way,  the diagonal part of (7.12.3)
is   of pure degree $|z|^2$ with coefficient $ \tau_s |Y|^2$, so it makes no contribution either.  Nor does
 the  term  (7.12.4.1), which  is manifestly a total derivative and hence automatically has 
zero integral.  The term
 (7.12.4.2) is a composition square and hence contributes only a product $P_o$-term of degree $|z|^4$ and a
$P_2$-term of degree $|z|^o$, namely (for j=0,2) the diagonal part of
$$P_j(z^2\dot{\bar{Y}}^2 + 2|z|^2|\dot{Y}|^2 +\bar{z}^2\dot{Y}^2 + \tau(z^2\bar{Y}^2
 + 2 |z|^2 |Y|^2 + \bar{z}^2Y^2, z^2\dot{\bar{Y}}^2 + 2|z|^2|\dot{Y}|^2 +\bar{z}^2\dot{Y}^2 + \tau(z^2\bar{Y}^2
 + 2 |z|^2 |Y|^2 + \bar{z}^2Y^2)$$
whose average over $\gamma$ has the form
$$(|z|^4 /or/|z|^o)\cdot  \frac{1}{L}\int_o^L [a |\dot{Y}|^4 
+ b\tau Re (\dot{\bar{Y}}^2 Y^2) + c\tau |\dot{Y}Y|^2
+d \tau^2 |Y|^4]ds \leqno(7.13.4.2)$$
where again the coefficients may vary between the two degrees.
Finally, we  the first term of (7.12.5) obviously has no diagonal
part while obviously the second term contributes
$$C  \frac{1}{L}\int_o^L \tau ds. \leqno(7.13.5)$$
This completes the analysis of the QBNF coefficients $B_{o4}, B_{o2}, B_{o0}$.

\section{Appendix : The classical Birkhoff normal form}

The method of \S2 for putting $\sqrt{\Delta}$ into quantum Birkhoff normal form began,
essentially, by putting the linearization of $\sqrt{\Delta}$ at $\gamma$ into normal form
by a linear symplectic transformation, and then proceeded by induction on the jet filtration
at $\gamma$.  The purpose of this appendix is to describe, rather briefly and
sketchily, how to put the principal symbol of
$\sqrt{\Delta}$ into Birkhoff normal form by an analogous method.  (We have not
found this particular algorithm in the literature, but it is quite likely that it, or
a much better algorithm, is well-known).  We hope that the classical algorithm will
help clarify the procedure in the quantum case. 

In the usual Fermi symplectic normal coordinates $(s,\sigma,y,\eta)$, we may write
the principal symbol of $\sqrt{\Delta}$ in the form
 $$\tilde{H}(s, \sigma, y ,\eta):= (g^{oo}(s,y)\sigma^2 
+ \sum_{ij=1}^n g^{ij}(s,y)\eta_i \eta_j)^{\half}$$
$$= \sigma (g^{oo}(s,y) + \sum_{ij=1}^n\frac{\eta_i \eta_j}{\sigma^2})^{\half}.$$
Taking the Taylor expansion at $y=\eta=0$ we get
$$\tilde{H}(s,\sigma,y,\eta) = \sigma(1 + \half[\sum_{ij=1}^n K_{ij}(s)y_i y_j
 + \sum_{i=1}^n \frac{\eta_i^2}{\sigma^2}] + \dots)$$
from which we extract the {\it linearized symbol}
$$\tilde{h}(s, \sigma, y, \eta) = \sigma + \half \sum_{ij=1}^n K_{ij}(s)y_i y_j \sigma
+ \sum_{i=1}^n \frac{\eta_i^2}{\sigma}. \leqno (A.1)$$ 
We make the symplectic change of variables 
$$\phi: (s,\sigma, y, \eta) \rightarrow (s', \sigma', y', \eta'):=(s + \half \frac{y \dot \eta}
{\sigma}, \sqrt{\sigma}y, \frac{\eta}{\sqrt{\sigma}})\leqno(A.2)$$
which transforms $\tilde{h}$ into
$$h(s, \sigma, y, \eta) = \sigma + \half \sum_{ij=1}^n K_{ij}(s)
y_i y_j + \sum_{i=1}^n \eta_i^2\leqno(A.3)$$
and $\tilde{H}$ into
$$H(s, \sigma, y, \eta) = h + h^{[3]} + \dots$$
with $\dots$ denoting  terms of order 3 and higher in $(y, \eta).$  (Such terms are of
order 3/2 with respect to the order in the isotropic calculus, while $h$ is of order 1,
which is the rationale for calling $h$ the linearized symbol).

The first step in putting $H$ into Birkhoff normal form is to put $h$ into Birkhoff normal
form
$$\hat{h}= \sigma + \half \sum_{i=1}^n \alpha_i (y_i^2 + \eta_i^2)$$
 by means of a  symplectic map.  Equivalently, we wish to convert the Hamilton equations
$$\begin{array}{l}   \frac{d}{ds}s=1 \\  \frac{d}{ds}\sigma =0 \\ \frac{d}{ds} y=\eta \\
 \frac{d}{ds}\eta = -K(s) y\end{array}$$
 into the linear equations
$$\begin{array}{l}   \frac{d}{ds}s=1 \\  \frac{d}{ds}\sigma=0 \\  \frac{d}{ds}q = \alpha p \\
 \frac{d}{ds}p = - \alpha q.\end{array}$$
We first do this in just the $(y,\eta, q,p)$ variables, with a symplectic map of the form
$(y, \eta) \rightarrow (q,p) = B(s) (y, \eta)$.  The condition on $B$ is that
$$\dot{B}B^{-1} + B\tilde{K}B^{-1} = \alpha \dot J\leqno(A.4)$$
where $\tilde{K}$ is the block anti-diagonal matrix with blocks $I$ and $-K$,
where $J$ is the usual block anti-diagonal matrix with blocks $\pm I$ and
where $\alpha \dot J$ is the block anti-diagonal matrix with coefficients $\pm \alpha_j$ 
replacing $\pm 1$ in $J$.  The equation (A.4) is equivalent to
$$-\frac{d}{ds} B^{-1} + \tilde{K}(s)B^{-1} \alpha \dot B^{-1} J \leqno(A.5)$$
which has the solution
$$B(s) = r_{\alpha}(s) a_{s}^{-1}.\leqno(A.6)$$
To make the map symplectic in all the $(s, \sigma,y,\eta)$ variables, we observe that
 the map
$$\psi_1: (s, \sigma, y, \eta) \rightarrow (s, \sigma - \half\sum_{i=1}^n \alpha_i (y_i^2
+ \eta_i^2), r_{\alpha}(s)(y, \eta))$$
 is symplectic with respect to $ds \wedge d\sigma + dy \wedge d\eta$ and satisfies
$$\psi_1 ^* (\sigma +  \half\sum_{i=1}^n \alpha_i (y_i^2+ \eta_i^2)) = \sigma.$$  Also,
the map
$$\psi_2: (s, \sigma, y, \eta) \rightarrow (s, \sigma + f(s, y, \eta), a_{s}^{-1} (y, \eta))$$
with
$f(s, y, \eta):= -\half \sum_{ij=1}^n K_{ij}(s)y_iy_j + \sum_{i=1}^n \eta_i^2$ 
is symplectic with respect to $ds \wedge d\sigma + dy \wedge d\eta$ 
and satisfies
$$\psi_2^* (\sigma) = \sigma + f.$$
It follows that $\chi_1:=\psi_2 \cdot \psi_1$ pulls back $h$ to its Birkhoff normal form.

Let us now write
$$H_1:=\chi_1^*(H) = \hat{h} + h_1^{[3]} + h_1^{[4]} + \dots$$
with $h_1^{[k]}(s, \sigma, y,\eta)$ of order 1 and  vanishing to order k at $(y, \eta) = (0,0).$ 
Following the algorithm for putting a Hamiltonian with non-degenerate minimum at $0$
into Birkhoff normal form [AM, 5.6.8, p.500]
we seek a symplectic map of the form
$$\chi_2 = exp  ad \sigma^{\half} F_3 $$
with $F_3 = F_3(s, y, \eta)$ vanishing to order 3 and
with  $$\{ \hat{h}, [\chi_2^*H_1]^{\leq3}\} = 0.\leqno(A.7)$$
Here, $\{,\}$ denotes the
Poisson bracket with respect to $ds\wedge d\sigma + dy \wedge d\eta$, $[F]^{\leq k}$ denotes
the terms of vanishing order $\leq k$ and $adF \cdot G:=\{F, G\}.$   The equation (A.7)
is equivalent to
$$ h_1^{[3]} + \{\sigma^{\half}F_3, \hat{h} \} \in ker(ad \hat{h}) \cap {\cal P}^3\leqno(A.8)$$
where ${\cal P}^k$ denotes the space of homogeneous polynomials of degree k in $(y, \eta)$
with smooth coefficients in $s$ and with an overall factor of   $\sigma^{\frac{k}{2}-1}$. 
In terms of the complex coordinates $z_j = y_j + i \eta_j$ on $\R^{2n}$ we may write
$$F_3 = \sum_{|j|+|k| = 3} c_{jk}(s) z^j \overline{z}^k$$
$$h_1^{[3]} = \sum_{|j|+|k| = 3} a_{jk}(s) z^j \overline{z}^k$$
with $j=(j_1,\dots,j_n), k=(k_1,\dots,k_n)$, and  we may write 
the Lie derivative ${\cal L}$ with respect to the
Hamilton vector field of $\hat{h}$ as
$${\cal L}= \frac{\partial}{\partial s} + \sum_{i=1}^n \alpha_i(z_i \frac{\partial}{\partial z_i}
- \overline{z}_i \frac{\partial}{\partial \overline{z}_i})$$  The monomials in
$ker ad {\cal L} \cap {\cal P}^k$ are of the form $z^{j} \overline{z}^k$ with 
$\langle \alpha, j-k \rangle =0$ which implies $j=k$ with our assumptions on $\alpha.$
No such terms occur for odd k, so equation (A.8) thus becomes
$$\sum_{|j|+|k| = 3} \dot{c}_{jk}(s) + \alpha \cdot (j - k) c_{jk}(s) = - a_{jk}(s).\leqno(A.9)$$
Expanding the periodic (or more generally almost periodic) functions
$c_{jk}$ and $a_{jk}$ in  Fourier series 
$$c_{jk}(s) = \sum e^{ims} \hat{c}_{jk}(m) \;\;\;\;\;\;a_{jk}(s)=\sum e^{ims}\hat{a}_{jk}(m)$$
we can solve (A.9) with 
$$\hat{c}_{jk}(m) = \frac{\hat{a}_{jk}(m)}{im + \langle \alpha, j-k\rangle}.\leqno(A.10)$$

Writing $H_2 = \chi_2^* \cdot \chi_1^*H$, we arrive at the analogous problem for the
fourth order terms.  As in the quantum case,the even steps behave a little differently from the
odd ones since now there can be terms $H_2^{[4]o}$ in $H_2^{[4]}$
with $|j|=|k|$ and hence which lie in $ker ad{\cal L} \cap {\cal P}^4.$
Since these terms already commute with ${\cal L}$ it suffices to solve the analogue of
(A9) with only the coefficients $a_{jk}(s)$ coming from $H_2^{[4]}- H_2^{[4]o}$.  

This
puts the terms up to fourth order in normal form, and the process continues inductively
to define symplectic maps $\chi_{N} \cdot \chi_{N-1} \cdot \dots \cdot \chi_1$ which
pulls back H to a normal form up to degree N.  Since $\chi_N = I mod O_N$ the infinite
product defines a smooth symplectic map which pulls back $H$ to a normal form modulo
$O_{\infty}$, that is, to its Birkhoff normal form. 
\medskip

\section{Index of Notation}

In the following, $\tau_L$ will denote the translation operator
$\tau_Lf(s,y) = f(s+L,y)$ on functions on $\R\times\R^n$. Operators $A$ then transform
under $\tau_L$ by $\tau_L A \tau_L^*$.

\noindent{\bf NI.1: Model objects}

\begin{tabular}{r||l}  Model object & $\cdot$ \\ \hline \hline
Functions & $\tau_L f = f$\\ \hline
Operators & $\tau_L A \tau_L^* = A$ \\ \hline
Maximal abelian algebra & ${\cal A}= <{\cal R}, I_1,\dots,I_n>$ \\ \hline
Distinguished element &${\cal R}:=\frac{1}{L}( L D_s + H_{\alpha})$ \\ \hline
Harmonic Oscillator &
$H_{\alpha}:= \frac{1}{2}\sum_{k=1}^n \alpha_k I_k $\\ \hline
Gaussians/Hermites  & $\gamma_o(y) = \gamma_{iI}(q) =
e^{-\half |y|^2}\;\;\;\;\;\;\;\;\gamma_q = C_q A_1^{*q_1}\dots A_n^{*q_n}$ 
\\
\hline
${\cal A}$-Eigenfunctions& $\phi_{kq}^o(s,y):= e_k(s)\otimes\gamma_q(y)$\\ \hline
${\cal R}$-Eigenvalues &${\cal R}\phi_{kq}^o(s,y)
= r_{kq} \phi_{kq}^o(s,y),\;\;\;\;\;\;\;
r_{kq} = \frac{1}{L} (2 \pi k + \sum_{j=1}^n (q_j + \frac{1}{2})
\alpha_j)$ \\ \hline
Scaled Laplacian &  
${\cal R}_h := \mu(\tilde{a})\Delta_h\mu(\tilde{a})^*\sim \sum_{m=o}^{\infty} h^{(-2
+ \frac{m}{2})} {\cal R}_{2 -\frac{m}{2}} $\\ \hline
Intertwiner to SC normal form &
$W_h:=\Pi_{k=1}^{\infty}W_{h\frac{k}{2}}\mu(r_{\alpha}),\;\;\;\;\;$
$\tilde{W}_{h \frac{k}{2}}:= exp(ih^{\frac{k}{2}} Q_{\frac{k}{2}})$
\\ \hline
S.C. Normal form & $W_h^*{\cal R}_h W_h|_o \sim h^{-2} 
+ \sum_{j=0}^{\infty} h^{k} f_k (I_1, \cdots, I_n)$ \\ \hline

\end{tabular}
\bigskip

\noindent{\bf NI.2: Twisted model objects}
\bigskip

\begin{tabular}{r||l}  Twisted model object & $\cdot$ \\ \hline \hline
Functions & $\tau_L f = \mu(r_{\alpha}(L))f$\\ \hline
Operators & $\tau_L A \tau_L^* = \mu(r_{\alpha}(L)) A \mu(r_{\alpha}(L))^*$  \\ \hline
Maximal abelian algebra & ${\cal A}= <D_s, I_1,\dots,I_n>$ \\ \hline
Distinguished element &$D_s$ \\ \hline
Harmonic Oscillator &
$H_{\alpha}:= \frac{1}{2}\sum_{k=1}^n \alpha_k I_k $\\ \hline
Gaussians/Hermites  & $\gamma_o(y) =  \gamma_{iI}(q) =
e^{-\half |y|^2}\;\;\;\;\;\;\;\;\gamma_q = C_q A_1^{*q_1}\dots A_n^{*q_n}$   \\ \hline
${\cal A}$-Eigenfunctions& $ e^{ir_{kq}s}\otimes \gamma_q(y) $\\ \hline
$D_s$-Eigenvalues &$D_s e^{ir_{kq}s}\otimes \gamma_q(y) 
= r_{kq} e^{ir_{kq}s}\otimes \gamma_q(y)$ \\ \hline
Scaled Laplacian & 
${\cal D}_h:= \mu(a)\Delta_h \mu(a)^* \sim \sum_{m=o}^{\infty} h^{(-2 + \frac{m}{2})}
{\cal D}_{2 -\frac{m}{2}} $ \\ \hline
Intertwiner to SC normal form &
$\tilde{W}_h:=\Pi_{k=1}^{\infty}\tilde{W}_{h\frac{k}{2}}\mu(r_{\alpha}),\;\;\;\;\;$
$\tilde{W}_{h \frac{k}{2}}:= exp(ih^{\frac{k}{2}} \tilde{Q}_{\frac{k}{2}})$
\\ \hline 
S.C. Normal form &   $\tilde{W}_h^*{\cal D}_h \tilde{W}_h|_o \sim h^{-2}
  + \sum_{j=0}^{\infty} h^{k} f_k (I_1, \cdots, I_n)$ \\ \hline

\end{tabular}
\bigskip

\noindent{\bf NI.3: Adapted model objects}
\bigskip

\begin{tabular}{r||l} Adapted model object & $\cdot$ \\ \hline \hline
 Functions & $\tau_L f = \mu(T)f$\\ \hline
Operators & $\tau_L A \tau_L^* = \mu(T) A \mu(T)^*$\\ \hline
Maximal abelian algebra & ${\cal A}= <{\cal L}, \Lambda_1\Lambda_1^* ,\dots,
\Lambda_n
\Lambda_n^*>$ \\
\hline Distinguished element & ${\cal L} = D_s - \frac{1}{2}(\sum_{j=1}^n D_{u_j}^2 +
 \sum_{ij=1}^n K_{ij}(s)
u_iu_j).$ \\ \hline
Harmonic Oscillator & $H_{\alpha}:= \frac{1}{2}\sum_{k=1}^n \alpha_k \Lambda_k
\Lambda_k^* $\\ \hline
Gaussians/Hermites & $\mu(a^{-1}) \gamma_o (s,u):=U_o(s,u)  = 
(det Y(s))^{-1/2} exp (\frac{i}{2} <\Gamma(s) u,u>,$
 $\;\;\;\Gamma(s) := \frac{dY}{ds} Y^{-1};$\\ 
$\cdot$ & $\mu(a^{-1})\gamma_q:=U_q =C_q \Lambda_1^{q_1}...\Lambda_n^{q_n} U_o. $\\
\hline
${\cal A}$-Eigenfunctions & $\mu(\tilde{a}_s)(\phi^o_{kq}) =\phi_{kq}:=e^{i r_{kq}
s}U_q(s,u) $
 \\ \hline
Scaled Laplacian & $\Delta_h \sim \sum_{m=0}^{\infty} h^{(-2 +m/2)}{\cal L}_{2-m/2}$
 \\ \hline
\end{tabular}
\bigskip

\noindent{\bf NI.4 Intertwining operators}

1. $\mu(r_{\alpha}):= \int_{S^1}^{\oplus} \mu(r_{\alpha}(s))ds$\\
2. $a_s:= \left( \begin{array}{ll} Im\dot{Y}(s)^* \;\;\;&ImY(s)^*\\
Re\dot{Y}(s)^*\;\;\;&Re Y(s)^*  \end{array} \right) = {\cal A}^*.$\\
3. ${\cal A}_L(s):= \left( \begin{array}{ll} L^{\half} Im\dot{Y}(s) \;\;\;&L^{-\half} Re\dot{Y}(s)\\
L^{-\half} Im Y(s)\;\;\;&L^{-\half} Re Y(s)  \end{array} \right)$
\bigskip

\begin{tabular}{r|l} Intertwining operator $A$ & $\tau_L A \tau_L^*$ \\ \hline 
$\mu(a):$ Adapted Model $\rightarrow$ Twisted Model &
$\tau_L \mu(a) \tau_L^* =  \mu(r_{\alpha})\mu(a)$ \\ \hline
$\mu(r_{\alpha})$ : Model $\rightarrow$ Twisted
Model & $\tau_L \mu(r_{\alpha}) \tau_L^* = \mu(r_{\alpha}(L))
\mu(r_{\alpha})$\\ \hline
With $\tilde{a}_s:=a_s r_{\alpha}(s), \mu(\tilde{a})$: Adpated Model $\rightarrow$
Model & $\tau_L \mu(\tilde{a})\tau_L^* = \mu(\tilde{a})$  \\ \hline

\end{tabular}

\end{document}